\documentclass[12pt, leqno]{article}
\usepackage{amsfonts}
\usepackage{amsmath}
\usepackage{amssymb}
\usepackage{amscd}

\begin{document}

\begin{center}
{\bf \Large Arithmetic topology in Ihara theory}
\end{center} 

 \vskip 2pt

\begin{center}
{\em To the memory of Professor Akito Nomura}
\end{center} 

\begin{center}
Hisatoshi Kodani, Masanori Morishita  and Yuji Terashima 
\end{center}

\footnote[0]{
$\;\;$ Key words: Ihara representation, $l$-adic Milnor invariants, Pro-$l$ Johnson homomorphisms, Pro-$l$ Magnus-Gassner cocycles, $l$-adic Alexander invariants, Ihara power seires\\
$\;\; $ H.K. is  partly supported by Grants-in-Aid for JSPS Fellows 14J12303. \\
$\;\;$ Y.T. is partly supported by Grants-in-Aid for Scientific Research (C)  25400083 and by \\
$\;\;$ CREST, JST.
}

{\small
{\bf Abstract.} Ihara initiated to study a certain Galois representation which may be seen as an arithmetic analogue of the Artin representation of a pure braid group. We pursue the analogies in Ihara theory further and give foundational results, following after some issues and their inter-relations in the theory of braids and links such as Milnor invariants, Johnson homomorphisms, Magnus-Gassner cocycles and Alexander invariants, and study relations with arithmetic in Ihara theory. 
}

\vspace{.8cm}

\begin{center}
{\bf Introduction}
\end{center}

Let $l$ be a prime number. Let $S$ be a set of ordered $r+1$ ($r \geq 2$) distinct $\overline{\mathbb{Q}}$-rational points on the projective line $\mathbb{P}^1$ over the rational number field  $\mathbb{Q}$, where $\overline{\mathbb{Q}}$ is an algebraic closure of $\mathbb{Q}$. Let $k := \mathbb{Q}(S \setminus \{ \infty \})$, the finite algebraic number field generated by coordinates of points in $S \setminus \{ \infty \}$.  Note that the absolute Galois group ${\rm Gal}_{k} := {\rm Gal}(\overline{\mathbb{Q}}/k)$  is the \'{e}tale fundamental group of ${\rm Spec}\, k$ so that it acts on the geometric fiber $\mathbb{P}^1_{\overline{\mathbb{Q}}}\, \setminus S$ of the fibration $\mathbb{P}^1_{k} \setminus S \rightarrow {\rm Spec}\, k$ and hence on the pro-$l$ \'{e}tale  fundamental group $\pi_1^{\scriptsize \mbox{pro}-l}(\mathbb{P}^1_{\overline{\mathbb{Q}}}\,  \setminus S) \simeq \frak{F}_r$, where $\frak{F}_r$ denotes the free pro-$l$ group on $r$ generators $x_1, \dots , x_r$. In [Ih1], Ihara initiated to study this monodromy Galois representation 
$$ {\rm Ih}_S :  {\rm Gal}_{k} \longrightarrow {\rm Aut}(\frak{F}_r), \leqno{(0.1)}$$
particularly for the case $S = \{ 0, 1, \infty\}$ and $k = \mathbb{Q}$, in connection with deep arithmetic such as Iwasawa theory on cyclotomy and complex multiplications of Fermat Jacobians.
 We note that the image of ${\rm Ih}_S$ is contained in the subgroup consisting of $\varphi \in {\rm Aut}(\frak{F}_r)$ such that $ \varphi(x_i) \sim x_i^{\alpha}$ (conjugate) for $1\leq i \leq r$  and $\varphi(x_1\cdots x_r) = (x_1\cdots x_r)^{\alpha}$ for some $\alpha \in \mathbb{Z}_l^{\times}.$

As explained in [Ih3], the Ihara representation (0.1) may be regarded as an arithmetic analogue of the Artin representation of a pure braid group ([Ar]). Let $P_r$ be the pure braid group with $r$ strings $(r \geq 2)$. Note that $P_r$ is the topological fundamental group of the configuration space $D^r \setminus \Delta$ of  ordered $r$ points on a 2-dimensional disk $D$, where $\Delta$ denotes the hyperdiagonal of $D^r$. For $Q = (z_1, \dots, z_r) \in D^r \setminus \Delta$, we also write the same $Q$ for the subset $\{ z_1,\dots , z_r\}$ of $D$. Then $P_r$ acts, as the monodoromy, on the fiber $D \setminus Q$ of the universal bundle over a point $Q \in D^r \setminus \Delta$ and hence on the topological fundamental group $\pi_1(D \setminus Q) \simeq F_r$, where $F_r$ denotes the free group on $r$ generators $x_1, \dots , x_r$. Thus we have the Artin representation
$$ {\rm Ar}_Q : P_r \longrightarrow {\rm Aut}(F_r), \leqno{(0.2)} $$
which is in fact isomorphic onto the subgroup $\varphi \in {\rm Aut}(F_r)$ such that $\varphi(x_i) \sim x_i$ for $1 \leq i \leq r$ and $\varphi(x_1\cdots x_r) = x_1\cdots x_r$. 

We may see the following analogy between the Ihara representation (0.1) and the Artin representation (0.2):
$${\small
\mbox{
\begin{tabular}{ | c | c |} 
\hline
 absolute Galois group &  pure braid group \\
		${\rm Gal}_{k}$ &  $P_r$\\   
\hline
$\mathbb{P}^1_{k} \, \setminus S \rightarrow {\rm Spec} \, k$ &  universal bundle over $D^r \setminus \Delta$\\
with geometric fiber $\mathbb{P}^1_{\overline{\mathbb{Q}}} \, \setminus S$ & with fibers $D \setminus Q$\\
\hline
Ihara representation of   ${\rm Gal}_{k}$  &  Artin representation of $P_r$  \\ 
on $\pi_1^{\scriptsize \mbox{pro}-l}(\mathbb{P}^1_{\overline{\mathbb{Q}}}\,  \setminus S) = \frak{F}_r$  & on $\pi_1(D \setminus Q) = F_r$\\
\hline
\end{tabular}
}
}
\leqno{(0.3)}
$$
\vskip 4pt

The aim of this paper is, based on the above analogy (0.3),  to give foundational results obtained by pursuing pro-$l$ analogues for the Ihara representation of various objects derived from the Artin representation. To be precise, we shall investigate arithmetic (pro-$l$) analogues in Ihara theory of the following issues (I) $\sim$ (IV) and their inter-relations in the theory of braids and links:
\vskip 5pt 
\noindent
 (I) Milnor invariants of links,\\
 (II) Johnson homomorphisms for the pure braid group $P_r$, \\
 (III)  Magnus-Gassner representations of $P_r$, \\
 (IV) Alexander invariants of links. 
\vspace{.25cm}\\
Milnor invariants are higher order linking numbers of a link introduced by Milnor in [Mi]. For a pure braid link, they are defined as follows. For $b \in P_r$ and each $i$ $(1\leq i \leq r)$, we can write ${\rm Ar}_{Q}(b)(x_i) = y_i(b)x_i y_i(b)^{-1}$ for the unique $y_i(b) \in F_r$, where the sum of exponents of $x_i$ in the word $y_i(b)$ is $0$. The Milnor number $\mu(b; i_1\cdots i_ni) \in \mathbb{Z}$ is then defined to be the coefficient of  $X_{i_1}\cdots X_{i_n}$ in  the Magnus expansion of $y_i(b)$:
$$ y_i(b) = \sum_{1\leq i_1,\dots , i_n \leq r} \mu(b; i_1\cdots i_ni)X_{i_1}\cdots X_{i_n} \;\; (x_j = 1 + X_j).$$
The Milnor invariant $\overline{\mu}(\hat{b}; i_1\cdots i_ni)$ is defined by taking modulo a certain indeterminacy $\Delta(b; i_1\cdots i_ni)$: $\overline{\mu}(\hat{b}; i_1\cdots i_ni) := \mu(b; i_1\cdots i_ni)$ mod $\Delta(b; i_1\cdots i_ni)$. It turns out that it is an invariant of the link $\hat{b}$ obtained by closing $b$.  The Milnor invariants are  also interpreted in terms of Massey products in the cohomology of the link group ([Ki], [T]). Johnson homomorphisms are useful means to study the structure of the mapping class group of a surface ([J1],  [J2], [Mt1], [Mt2]). The main tools are algebraic and applicable to the study of the automorphism group ${\rm Aut}(F_r)$ of a free group $F_r$ ([Ka], [Sa]). Johnson homomorphisms describe the action of a certain filtration of ${\rm Aut}(F_r)$ on the nilpotent quotients $F_r/F_r(n)$ for $n \geq 1$, where $F_r(n)$ is the $n$-th term of the lower central series of $F_r$. Since the pure braid group $P_r$ is a subgroup of  the mapping class group of $r$ punctured disk, the theory of Johnson homomorphisms can also be  applied   to  $P_r$. It was shown in [Ko1], [Ko3; Chapter 1] that the Johnson homomorphisms are described by Milnor invariants of pure braid links. Magnus cocycles are crossed homomorphisms of $P_r$  defined by using the Fox free derivation ([B; 3.1, 3.2], [F]). The Gassner representation ${\rm Gass}$ of $P_r$ is a particular case of Magnus cocycles over the Laurent polynomial ring of $r$ variables and the determinant $\det({\rm Gass}(b) - I)$ gives the Alexander invariant which is  a polynomial invariant of the link $\hat{b}$ ([B; 3.3]). The relations of the Gassner representations with Johnson homomorphisms and Milnor invariants were studied in  [Ko2], [Ko3; Chapter 2]. \\

In this paper, based on the analogy (0.3), we shall study arithmetic analogues in Ihara theory of these issues  (I) $\sim$ (IV). The contents of this paper are organized as follows. In Section 1, we recall the construction of the Ihara representations and some basic results. In Section 2, we define $l$-adic Milnor numbers for each element in ${\rm Gal}_{k}$  and $l$-adic Milnor invariants for certain primes of $k(\zeta_{l^{\infty}})$, the field obtained by adjoining all $l$-powerth roots of unity to $k$. We introduce the pro-$l$ link group of each element of ${\rm Gal}_k$ and give a cohomological interpretation of $l$-adic Milnor invariants in terms of Massey products in the cohomology of the pro-$l$ link group.  In Section 3,  we present a general theory of the pro-$l$ Johnson map and pro-$l$ Johnson homomorphisms for the absolute Galois group ${\rm Gal}_{k}$. A similar theory has been developed in the context of non-abelian Iwasawa theory ([MT]). Among other things, we describe the pro-$l$ Johnson homomorphisms in terms of $l$-adic Milnor numbers. Sections 2 and 3 may be regarded as an arithmetic counterpart of [Ko1], [Ko3; Chapter 1]. In Section 4, we introduce the pro-$l$ Magnus cocycles   of ${\rm Gal}_{k}$ by using the pro-$l$ Fox free calculus, and give a relation with pro-$l$ Johnson homomorphisms. We consider the pro-$l$ (reduced) Gassner cocycle of ${\rm Gal}_k$  as a special case and express it by $l$-adic Milnor numbers. Section 4 may be regarded as an arithmetic counterpart of [Ko2], [Ko3; Chapter 2]. We note that Oda's unpublished notes ([O1], [O2]) also concern some issues related to Sections 3, 4. In Section 5, we introduce the pro-$l$ link module and $l$-adic Alexander invariants. In Section 6, we consider the case that $S = \{ 0,1,\infty \}$. We show that the Ihara power series $F_g(u_1,u_2)$ ($g \in {\rm Gal}_{\mathbb{Q}}$) introduced in [Ih1] coincides with our pro-$l$ reduced Gassner cocycle, and give a formula which expresses  $F_g(u_1,u_2)$ in terms of $l$-adic Milnor numbers.  Accordingly, using our formula and Ihara's formula, we express the Jacobi sum in $\mathbb{Q}(\zeta_{l^n})$ as a $(\zeta_{l^n}-1)$-adic expansion with coefficients $l$-adic Milnor numbers. Finally, combining our formula  and the result by Ihara, Kaneko and Yukinari ([IKY]), we give some formula relating Soul\'{e} characters ([So]) with $l$-adic Milnor numbers.  \\

This paper forms a (part of) elementary and group-theoretical foundation of {\it arithmetic topology} in Ihara theory. In the forthcoming papers, we shall study some connections of $l$-adic Milnor invariants and pro-$l$ Johnson homomorphisms in this paper with arithmetic of multiple power residue symbols in [Am], [Ms1], [Ms2; Chapter 8] and the works of Wojtkowiak on $l$-adic iterated integrals and $l$-adic polylogarithms ([NW], [W1] $\sim$ [W4] etc). See Remark 2.2.12. We shall also study arithmetic analogues of some issues in {\it quantum topology} such as Habegger-Masbaum's theorem on the relation between Milnor invariants and Kontsevich integrals ([HM]). \\
\\
{\it Notation.} We denote by $\mathbb{Z}$, $\mathbb{Q}$ and $\mathbb{C}$ the ring of rational integers, the field of rational numbers and the field of complex numbers, respectively.\\
Throughout this paper, $l$ denotes a fixed prime number.  We denote by $\mathbb{Z}_l$ and $\mathbb{Q}_l$ the ring of $l$-adic integers and the field of $l$-adic numbers, respectively.\\
For $a, b$ in a group $G$, $a \sim b$ means that $a$ is conjugate to $b$ in $G$. For subgroups $A, B$ of a topological group $G$, $[A,B]$ stands for the closed subgroup of $G$ generated by commutators  $[a,b] := aba^{-1}b^{-1}$ for all $a \in A, b \in B$. \\
For a positive integer $n$ and a ring $R$ with identity, ${\rm M}(n; R)$ denotes the ring of $n \times n$ matrices whose entries are in $R$ and ${\rm GL}(n; R)$ denotes the group of invertible elements of ${\rm M}(n;R)$.
\\

\begin{center}
{\bf 1.  The Ihara representation}
\end{center}

In this section,  we recall the set-up and some results on the Galois representation introduced by Ihara in [Ih1]. \\
\\
{\bf 1.1. The outer Galois representation.} Let $x_1, \dots, x_r$ be $r$ letters $( r\geq 2)$ and let $F_r$ denote the free group of rank $r$ on $x_1, \dots , x_r$. Let $x_{r+1}$ be the element of $F_r$ defined by $x_1\cdots x_r x_{r+1} = 1$ so that $F_r$ has the presentation $F_r = \langle x_1, \dots , x_r, x_{r+1} | x_1\cdots x_r x_{r+1} = 1 \rangle$. Let $\frak{F}_r$ denote the pro-$l$ completion of $F_r$. 
Let ${\rm Aut}(\frak{F}_r)$ (resp. ${\rm Int}(\frak{F}_r)$) denote the group of topological automorphisms (resp. inner-automorphisms) of $\frak{F}_r$ with compact-open topology. We note that any abstract automorphism of $\frak{F}_r$ is bicontinuous ([DDMS; Corollary 1.22]) and that ${\rm Aut}(\frak{F}_r)$ is virtually a pro-$l$ group ([DDMS; Theorem 5.6]). Let $H$ be the abelianization of $\frak{F}_r$, $H := \frak{F}_r/[\frak{F}_r, \frak{F}_r]$, and let $\pi : \frak{F}_r \rightarrow H$ be the abelianization homomorphism. For $f \in \frak{F}_r$, we let $[f] : = \pi(f)$. We set $X_i := [x_i]$ $(1\leq i \leq r+1)$ for simplicity so that $H$ is the free $\mathbb{Z}_l$-module with basis $X_1, \dots , X_r$ and we have $X_1 + \cdots + X_r + X_{r+1} = 0$. Each $\varphi \in {\rm Aut}(\frak{F}_r)$ induces an automorphism of the $\mathbb{Z}_l$-module $H$ which is denoted by $[\varphi] \in {\rm GL}(H)$.

Let $\overline{\mathbb{Q}}$ be the field of algebraic numbers in $\mathbb{C}$. Let $S$ be a given set of ordered $r+1$ $\overline{\mathbb{Q}}$-rational points  $P_1, \dots, P_{r+1}$ on the projective line $\mathbb{P}^1_{\mathbb{Q}}$  and we suppose that $P_1 =0, P_2 = 1$ and $P_{r+1} = \infty$. Let $k := \mathbb{Q}(S \setminus \{ \infty \})$, the finite algebraic number field generated over $\mathbb{Q}$ by coordinates of $P_1, \dots , P_r$, so that all $P_i$'s are $k$-rational points of $\mathbb{P}^1$.  Let ${\rm Gal}_{k} := {\rm Gal}(\overline{\mathbb{Q}}/k)$ be the absolute Galois group  of $k$ equipped with the Krull topology. Note that ${\rm Gal}_{k}$ is the \'{e}tale fundamental group $\pi_1^{\scriptsize \mbox{\'{e}t}}({\rm Spec} \, k)$ with the base point ${\rm Spec} \, \overline{\mathbb{Q}} \rightarrow {\rm Spec} \, k$.  Let $\pi_1^{\scriptsize \mbox{pro}-l}(\mathbb{P}^1_{\overline{\mathbb{Q}}}\,  \setminus S)$ denote the maximal pro-$l$ quotient of the \'{e}tale fundamental group of $\mathbb{P}^1_{\overline{\mathbb{Q}}}\,  \setminus S$ with a base point ${\rm Spec} \, \overline{\mathbb{Q}} \rightarrow \mathbb{P}^1_{\overline{\mathbb{Q}}}\,  \setminus S$ which lifts ${\rm Spec} \, \overline{\mathbb{Q}} \rightarrow {\rm Spec} \, k$. By [G; XII, Corollaire 5.2], $\pi_1^{\scriptsize \mbox{pro}-l}(\mathbb{P}^1_{\overline{\mathbb{Q}}}\,  \setminus S)$ is the pro-$l$ completion of the topological fundamental group $\pi_1(\mathbb{P}^1(\mathbb{C}) \setminus S)$. We fix once and for all an identification of $F_r$   with $\pi_1(\mathbb{P}^1(\mathbb{C}) \setminus S)$ obtained by associating to each $x_i$  the homotopy class of a small loop around $P_i$ and hence an identification  of $\pi_1^{\scriptsize \mbox{pro}-l}(\mathbb{P}^1_{\overline{\mathbb{Q}}}\,  \setminus S)$ with $\frak{F}_r$. 

The absolute Galois group  ${\rm Gal}_{k} = \pi_1^{\scriptsize \mbox{\'{e}t}}({\rm Spec} \, k)$  acts, as the monodromy, on the geometric fiber $\mathbb{P}^1_{\overline{\mathbb{Q}}}\, \setminus S$ of the fibration $\mathbb{P}^1_{k} \setminus S \rightarrow {\rm Spec}\, k$ and hence acts continuously on the pro-$l$ fundamental group $\pi_1^{\footnotesize{ \mbox{pro}-l}}(\mathbb{P}^1_{\overline{\mathbb{Q}}}\,  \setminus S) = \frak{F}_r$. The effect of changing a base point of $\mathbb{P}^1_{\overline{\mathbb{Q}}}\,  \setminus S$ is given as an inner automorphism of $\frak{F}_r$. Thus we have the continuous outer Galois representation 
$$ \Phi_S : {\rm Gal}_{k} \longrightarrow {\rm Out}(\frak{F}_r) := {\rm Aut}(\frak{F}_r)/{\rm Int}(\frak{F}_r). \leqno{(1.1.1)}$$

In terms of the field extensions, the representation $\Phi_S$ is described as follows. Let $t$ be a variable over $k$. We regard $\mathbb{P}^1$ as the $t$-line and so the function field $K$ of $\mathbb{P}^1_{\overline{\mathbb{Q}}}$ is the rational function field $\overline{\mathbb{Q}}(t)$. The $k$-rational points $P_i$ are identified with places of $K/\overline{\mathbb{Q}}$. Let $M$ be the maximal pro-$l$ extension of $K$ unramfied outside $P_i$ ($1\leq i \leq r+1$). We fix once and for all an  identification of $\frak{F}_r$  with ${\rm Gal}(M/K)$ obtained by associating to each $x_i$ a topological generator of the inertia group of an extension $P_i^M$ of $P_i$ to a place of $M$. Since $P_i$'s are $k$-rational, $M/k(t)$ is a Galois extension and so we have the exact sequence
$$ 1 \rightarrow \frak{F}_r =  {\rm Gal}(M/K) \rightarrow {\rm Gal}(M/k(t)) \rightarrow {\rm Gal}(K/k(t)) = {\rm Gal}_{k} \rightarrow 1.$$
For $g \in {\rm Gal}_{k}$, choose $\tilde{g} \in {\rm Gal}(M/k(t))$ which lifts $g$. Consider the action of ${\rm Gal}_{k}$ on ${\rm Gal}(M/K)$ defined by $f \mapsto \tilde{g} f \tilde{g}^{-1}$ and regard it as an automorphism of $\frak{F}_r$ via the isomorphism $\iota$. The effect of changing  a lift $\tilde{g}$ is given as an inner automorphism of $\frak{F}_r$. Thus we obtain the representation $\Phi_S$. Note further that  $g \circ P_i^M \circ \tilde{g}^{-1}$ is a place of $M$ which coincides with $P_i^M$ on $K$ ($1\leq i \leq r+1$). So we have $g \circ P_i^M \circ \tilde{g}^{-1} \circ h = P_i^M$ for some $h \in {\rm Gal}(M/K)$ so that $h^{-1}\tilde{g} x_i \tilde{g}^{-1} h$ is a topological generator of the inertia group of $P_i^M$. Hence $\tilde{g} x_i \tilde{g}^{-1} \sim x_i^{c_i}$ for some $c_i$ in $\mathbb{Z}_l$, the ring of $l$-adic integers. We pass to the abelianization $H$. Applying the conjugate  by $\tilde{g}$ on the equality $X_1 + \cdots  + X_{r+1} = 0$ in $H$, we have $c_1 X_1 + \cdots + c_{r+1} X_{r+1} = 0$. From these equations, we have $c_1 = \cdots = c_{r+1}$. Therefore the action  of ${\rm Gal}_{k}$ on $\frak{F}_r$ gives an element of the subgroup $\tilde{P}(\frak{F}_r)$ of ${\rm Aut}(\frak{F}_r)$  defined by
$$ \tilde{P}(\frak{F}_r) := \{ \varphi \in {\rm Aut}(\frak{F}_r) \, | \, \varphi(x_i) \sim x_i^{N(\varphi)} \;  (1\leq i \leq r+1) \; \mbox{for some}\; N(\varphi) \in \mathbb{Z}_l^{\times} \}.$$
Here the exponent $N(\varphi)$, called the {\em norm} of $\varphi$, gives a homomorphism $N : {\rm Aut}(\frak{F}_r) \rightarrow \mathbb{Z}_l^{\times}$. So each $\varphi \in \tilde{P}(\frak{F}_r)$ acts on the abelianization $H$ by the multiplication by $N(\varphi)$, $[\varphi](X_i) = N(\varphi)X_i$ for $1\leq i \leq r$.  It is easy to see ${\rm Int}(\frak{F}_r) \subset \tilde{P}(\frak{F}_r)$. Thus we have the outer Galois representation (1.1.1)
$$ \Phi_S : {\rm Gal}_{k} \longrightarrow \tilde{P}(\frak{F}_r)/{\rm Int}(\frak{F}_r). \leqno{(1.1.2)}$$
\\
{\bf 1.2. The Ihara representation.} We will lift $\Phi_S$ to a representation in ${\rm Aut}(\frak{F}_r)$. For this, consider the subgroup $P(\frak{F}_r)$ of $\tilde{P}(\frak{F}_r)$ defined by
$$P(\frak{F}_r) := \left\{ \varphi \in {\rm Aut}(\frak{F}_r) \, \Big| \, \begin{array}{l} \varphi(x_i) \sim x_i^{N(\varphi)} \;  (1\leq i \leq r-1) \;\varphi(x_r) \approx x_r^{N(\varphi)}, \\
                                                                                              \varphi(x_{r+1}) = x_{r+1}^{N(\varphi)}\; \mbox{for some}\; N(\varphi) \in \mathbb{Z}_l^{\times} \end{array} \right\},   \leqno{(1.2.1)}$$
where $\approx$ denotes conjugacy by an element of the subgroup $\frak{K}$ of $\frak{F}_r$ generated by $[\frak{F}_r,\frak{F}_r]$ and $x_1, \dots , x_{r-2}$. We denote by $P^1(\frak{F}_r)$ the kernel of $N|_{P(\frak{F}_r)}$:
$$  P^1(\frak{F}_r) := \left\{ \varphi \in {\rm Aut}(\frak{F}_r) \, \Big| \, \begin{array}{l} \varphi(x_i) \sim x_i \;  (1\leq i \leq r-1) \;\varphi(x_r) \approx x_r, \\
                                                                                               \;\; \varphi(x_{r+1}) = x_{r+1} \end{array} \right\}.  $$

The following proposition was proved in [Ih1; Proposition 3, page 55] for the case $r = 2$ and stated in [Ih3; page 252] for the general case.\\
 \\
{\bf Proposition 1.2.2.} {\em The natural homomorphism ${\rm Aut}(\frak{F}_r) \rightarrow {\rm Aut}(\frak{F}_r)/{\rm Int}(\frak{F}_r)$ induces the isomorphism $P(\frak{F}_r) \simeq \tilde{P}(\frak{F}_r)/{\rm Int}(\frak{F}_r)$.}  The representatives in $P(\frak{F}_r)$ of $\tilde{P}(\frak{F}_r)/{\rm Int}(\frak{F}_r)$ are called {\em Bely\u\i's lifts}. \\
\\
{\em Proof.}  Although the proof is similar to that for $r=2$, we give a concise proof for the sake of readers. First, we note that the centralizer of $x_i$ in $\frak{F}_r$ is $\langle x_i \rangle = x_i^{\mathbb{Z}_l}$ for $1\leq i \leq r+1$. \\
Injectivity:  Suppose $\varphi \in P(\frak{F}_r)$ and $\varphi = {\rm Int}(f)$ with $f \in \frak{F}_r$. Then $f x_{r+1} f^{-1} = x_{r+1}^{N(\varphi)}$. Passing to $H$, we see $N(\varphi) = 1$ and so $f$ is in the centralizer of $x_{r+1}$. Hence $f = x_{r+1}^a$ for some $a \in \mathbb{Z}_l$. Since $\varphi \in P(\frak{F}_r)$, $f x_r f^{-1} = \varphi(x_r) = g x_r g^{-1}$ for some $g \in \frak{K}$ and hence $g^{-1}f = x_r^b$ for some $b \in \mathbb{Z}_l$. Passing to the abelianization $H$, we find $ a = b = 0$. Hence  $f =1$ and $\varphi = 1$. \\
Surjectivity: Take $\varphi \in  \tilde{P}(\frak{F}_r)$. Multiplying $\varphi$ by an element of ${\rm Int}(\frak{F}_r)$, we may assume $\varphi(x_{r+1}) = x_{r+1}^{N(\psi)}$. Set $\varphi(x_r) = g x_r g^{-1}$ with $g \in \frak{F}_r$. Write $[g] = c_1 X_1 + \cdots c_{r}X_{r}$ in $H$ $(c_i \in \mathbb{Z}_l)$ and let $\varphi_1 :={\rm Int}(x_{r-1}^{-c_{r-1}}x_{r}^{-c_{r}}) \circ \varphi$. Then $\varphi_1(x_r) = g_1 x_r g_1^{-1}$ and $g_1 := x_{r-1}^{-c_{r-1}}x_{r}^{-c_{r}}g \in \frak{K}$. Hence $\varphi_1 \in P(\frak{F}_r)$ and $\varphi \equiv \varphi_1$ mod ${\rm Int}(\frak{F}_r)$. $\;\; \Box$\\
\\
By Proposition 1.2.2, we can lift $\Phi_S$ of (1.1.2) to the representation in $P(\frak{F}_r)$, denoted by  ${\rm Ih}_S$:
$$ {\rm Ih}_S : {\rm Gal}_{k} \longrightarrow P(\frak{F}_r), \leqno{(1.2.3)}$$
which we call the {\em Ihara representation} associated to $S$. Let $\Omega_S$ denote the subfield of $\overline{\mathbb{Q}}$ corresponding to the kernel of ${\rm Ih}_S$ so that ${\rm Ih}_S$ factors through the Galois group ${\rm Gal}(\Omega_S/k)$:
$$ {\rm Ih}_S :  {\rm Gal}(\Omega_S/k) \longrightarrow P(\frak{F}_r). \leqno{(1.2.4)}$$

We recall some arithmetic properties on the ramification in the Galois extension $\Omega_S/k$. For this, let us prepare some notations. Let $\zeta_{l^n}$ be a primitive $l^n$-th root of unity for a positive integer $n$ such that $(\zeta_{l^{n+1}})^l = \zeta_{l^n}$ for $n \geq 1$. We set $k(\zeta_{l^{\infty}}) := \cup_{n \geq 1}  k(\zeta_{l^n})$. The $l$-{\em cyclotomic character} $\chi_l : {\rm Gal}_{k} \rightarrow \mathbb{Z}_l^{\times}$ is defined by $g(\zeta_{l^n}) = \zeta_{l^n}^{\chi_l(g)}$ for $g \in {\rm Gal}_{k}$. Finally we define the set $R_S$ of finite primes of $k$ associated to $S$ as follows:  Let $s_i$ be the coordinate of $P_i$ for $1\leq i \leq r$, and let ${\cal O}_S$ be the integral closure of  $ \mathbb{Z}[ l^{-1}, \, (s_i - s_j)^{-1}  (1\leq i \neq j \leq r) ]$ in $k$. We then define $R_S$ by the maximal spectrum
$$ R_S := {\rm Spm} \, {\cal O}_S. \leqno{(1.2.5)}$$
\\
{\bf Theorem 1.2.6.} {\em Notations being as above, the following assertions hold}:\\
(1) ([Ih1; Proposition 2, page 53]). {\em $N \circ {\rm Ih}_S : {\rm Gal}_{k} \rightarrow \mathbb{Z}_l^{\times}$ coincides with $\chi_l$.  In particular, the restriction of $\varphi_S$ to ${\rm Gal}_{k(\zeta_{l^{\infty}})} := {\rm Gal}(\overline{\mathbb{Q}}/k(\zeta_{l^{\infty}}))$, denoted by ${\rm Ih}_S^1$, gives the representation
$${\rm Ih}_S^1  : {\rm Gal}_{k(\zeta_{l^{\infty}})} \longrightarrow P^1(\frak{F}_r)$$
and we  have $k(\zeta_{l^{\infty}}) \subset \Omega_S$.}\\
(2) ([AI; Proposition 2.5.2, Theorem 3]). {\em The Galois extension $\Omega_S/k$ is unramified over $R_S$ and $\Omega_S/k(\zeta_l)$ is a pro-$l$ extension.}\\
\\
{\bf  Remark 1.2.7} (cf. [Ih2]). By Artin's theorem ([Ar], [B; Theorem 1.9]), the Artin representation ${\rm Ar}_Q$ of the pure braid group $P_r$ in Introduction induces the isomorphism
$$ {\rm Ar}_Q : P_r \stackrel{\sim}{\longrightarrow} \{ \varphi \in {\rm Aut}(F_r) \, | \, \varphi(x_i) \sim x_i \;  (1\leq i \leq r), \; \varphi(x_1\cdots x_r) = x_1\cdots x_r \}.$$
So the  representation ${\rm Ih}_S^1 : {\rm Gal}_{k(\zeta_{l^{\infty}})} \rightarrow P^1(\frak{F}_r)$ (resp. ${\rm Ih}_S : {\rm Gal}_{k}  \rightarrow P(\frak{F}_r)$) may be seen as an (resp. extended) arithmetic analogue of the Artin representation ${\rm Ar}_Q$.\\

\begin{center}
{\bf 2.  $l$-adic Milnor invariants and pro-$l$ link groups}
\end{center}
\vspace{.2cm}
{\bf 2.1. Pro-$l$ Magnus expansions.}  Let $\{  \frak{F}_r(n) \}_{n \geq 1}$ be the lower central series of $\frak{F}_r$ defined by
$$\frak{F}_r(1) := \frak{F}_r, \;\; \frak{F}_r(n+1) := [\frak{F}_r(n), \frak{F}_r] \;\; (n \geq 1).$$
Note that each $\frak{F}_r(n)$ is a closed normal subgroup of $\frak{F}_r$ so that $\frak{F}_r(n)/\frak{F}_r(n+1)$ is central in $\frak{F}_r/\frak{F}_r(n+1)$, and that each $\frak{F}_r(n)$ is a finitely generated pro-$l$ group ([DDMS; 1.7, 1.14]). As in Section 1, let $H$ denote the abelianization of $\frak{F}_r$:
$$H := {\rm gr}_1(\frak{F}_r) = H_1(\frak{F}_r, \mathbb{Z}_l),$$
which is the free $\mathbb{Z}_l$-module with basis $X_1, \dots , X_r$, where $X_i$ is the image of $x_i$ in $H$. Let $T(H)$ be the
tensor algebra of $H$ over $\mathbb{Z}_l$ defined by $\bigoplus_{n \geq 0} H^{\otimes n}$, where $H^{\otimes 0} := \mathbb{Z}_l$ and $H^{\otimes n} := H \otimes_{\mathbb{Z}_l} \cdots \otimes_{\mathbb{Z}_l} H$ ($n$ times) for $n \geq 1$. It is nothing but the non-commutative polynomial algebra $\mathbb{Z}_l \langle X_1, \dots , X_r  \rangle$ over $\mathbb{Z}_l$ with variables $X_1, \dots , X_r$:
$$ T(H) = \bigoplus_{n \geq 0} H^{\otimes n} = \mathbb{Z}_l \langle X_1, \dots , X_r \rangle.$$
 Let $\widehat{T}(H)$ be the completion of $T(H)$ with respect to the $\frak{m}_T$-adic topology, where $\frak{m}_T$ is the maximal two-sided ideal of $T(H)$ generated by $X_1, \dots , X_r$ and $l$. It is the infinite product $\prod_{n \geq 0} H^{\otimes n}$, which is nothing but the {\em Magnus algebra}  $\mathbb{Z}_l \langle \langle X_1, \dots , X_r \rangle \rangle$ over $\mathbb{Z}_l$, namely,  the algebra of non-commutative formal power series (called {\em Magnus power series}) over $\mathbb{Z}_l$ with variables $X_1, \dots , X_r$:
$$ \widehat{T}(H)= \prod_{n \geq 0} H^{\otimes n} = \mathbb{Z}_l\langle \langle X_1, \dots , X_r \rangle \rangle.$$
For $n \geq 0$, we set $\widehat{T}(n) := \prod_{m \geq n} H^{\otimes m}$. The {\em degree} of a Magnus power series $\Phi$, denoted by ${\rm deg}(\Phi)$, is defined to be the minimum $n$ such that $\Phi \in \widehat{T}(n)$. We note that  $H^{\otimes n}$ is the free $\mathbb{Z}_l$-module on monomials $X_{i_1}\cdots X_{i_n}$ $(1\leq i_1,\dots , i_n \leq r)$ of degree $n$ and $\widehat{T}(n)$  consists of Magnus power series of degree $\geq n$.  

Let $\mathbb{Z}_l[[\frak{F}_r]]$ be the complete group algebra of $\frak{F}_r$ over $\mathbb{Z}_l$ and let $\epsilon_{\mathbb{Z}_l[[\frak{F}_r]]} : 
\mathbb{Z}_l[[\frak{F}_r]] \rightarrow \mathbb{Z}_l$ be the augmentation homomorphism with the augmentation ideal $I_{\mathbb{Z}_l[[\frak{F}_r]]} := {\rm Ker}(\epsilon_{\mathbb{Z}_l[[\frak{F}_r]]})$.
 The correspondence $x_i \mapsto 1 +X_i$ $(1\leq i \leq r)$ gives rise to the isomorphism of topological $\mathbb{Z}_l$-algebras
$$ \Theta : \mathbb{Z}_l[[\frak{F}_r]]\; \stackrel{\sim}{\longrightarrow}  \widehat{T}(H),  \leqno{(2.1.1)}$$
which we call the {\em pro-$l$  Magnus isomorphism}. Here $I_{\mathbb{Z}_l[[\frak{F}_r]]}^n$ corresponds, under $\Theta$,  to $\widehat{T}(n)$ for $n \geq 0$. For $\alpha \in \mathbb{Z}_l[[\frak{F}_r]]$, $\Theta(\alpha)$ is called the {\it pro-$l$ Magnus expansion} of $\alpha$. In the following, for a multi-index $I = (i_1 \cdots i_n)$,  $1\leq i_1,\dots , i_n \leq r$, we set
$$ |I| := n \; \mbox{and} \; X_I  := X_{i_1} \cdots X_{i_n}.$$
We call the coefficient of $X_I$ in $\Theta(\alpha)$ the {\em $l$-adic Magnus coefficient} of $\alpha$ for $I$  and denote it by $\mu(I;\alpha)$. So we have
$$ \Theta(\alpha) = \epsilon_{\mathbb{Z}_l[[\frak{F}_r]]}(\alpha) + \sum_{|I|\geq 1} \mu(I;\alpha) X_I.  \leqno{(2.1.2)}$$
Restricting $\Theta$ to $\frak{F}_r$, we have an injective group homomorphism, denoted by the same $\Theta$, 
$$ \Theta : \frak{F}_r \; \hookrightarrow \; 1 + \widehat{T}(1),    \leqno{(2.1.3)}$$
which we call the {\it pro-$l$ Magnus embedding} of $\frak{F}_r$ into $1 + \widehat{T}(1)$. 

Here are some basic properties of $l$-adic Magnus coefficients:\\
(2.1.4) For $\alpha, \beta \in \mathbb{Z}_l[[\frak{F}_r]]$ and a multi-index $I$, we have
$$ \mu(I;\alpha \beta) = \sum_{I = AB} \mu(A;\alpha)\mu(B;\beta),$$
where the sum ranges over all pairs $(A,B)$ of multi-indices such that $AB = I$, and we understand that $\mu(A;\alpha) =  \epsilon_{\mathbb{Z}_l[[\frak{F}_r]]}(\alpha)$ (resp. $\mu(B;\beta) =  \epsilon_{\mathbb{Z}_l[[\frak{F}_r]]}(\beta)$) if $|A| = 0$ (resp. $|B| = 0$). \\
(2.1.5) ({\it Shuffle relation}) For $ f \in \frak{F}_r$ and multi-indices $I, J$ with $|I|, |J| \geq 1$, we have
$$ \mu(I;f)\mu(J;f) = \sum_{A \in {\rm Sh}(I,J)} \mu(A;f),$$
where ${\rm Sh}(I,J)$ denotes the set of the results of all shuffles of $I$ and $J$ ([CFL]).\\
(2.1.6) For $f \in \frak{F}_r$ and $d \geq 2$, we have
$$ \begin{array}{ll} \mu(I;f) = 0 \; \mbox{for} \; |I| < d \;\; \mbox{i.e.,}\, {\rm deg}(\Theta(f-1)) \geq d  & \Longleftrightarrow \; f \in \frak{F}_r(d) \\
                                                                 &  \Longleftrightarrow \; f -1 \in  I_{\mathbb{Z}_l[[\frak{F}_r]]}^d. 
\end{array}$$

 An automorphism $\varphi$ of the topological $\mathbb{Z}_l$-algebra $\mathbb{Z}_l[[\frak{F}_r]]$ (resp. $\widehat{T}(H)$) is said to be {\em filtration-preserving} if $\varphi( I_{\mathbb{Z}_l[[\frak{F}_r]]}^n) = I_{\mathbb{Z}_l[[\frak{F}_r]]}^n$ (resp. $\varphi(\widehat{T}(n)) = \widehat{T}(n)$) for all $n \geq 1$. Let ${\rm Aut}^{\rm fil}(\mathbb{Z}_l[[\frak{F}_r]])$ (resp. ${\rm Aut}^{\rm fil}(\widehat{T}(H))$) be the group of filtration-preserving automorphisms of the topological $\mathbb{Z}_l$-algebras $\mathbb{Z}_l[[\frak{F}_r]]$ (resp. $\widehat{T}(H)$). The pro-$l$ Magnus isomorphism $\Theta$ in (2.1.1) induces the isomorphism
 $$ {\rm Aut}^{\rm fil}(\mathbb{Z}_l[[\frak{F}_r]]) \stackrel{\sim}{\longrightarrow} {\rm Aut}^{\rm fil}(\widehat{T}(H)); \;\; \varphi \mapsto \Theta \circ \varphi \circ \Theta^{-1}. \leqno{(2.1.7)}$$
 In the following we set
 $$ \varphi^* := \Theta \circ \varphi \circ \Theta^{-1}. \leqno{(2.1.8)}$$
 We note by (2.1.6)  that any automorphism $\varphi$ of $\frak{F}_r$ can be extended uniquely to a filtration-preserving topological automorphism of $\mathbb{Z}_l[[\frak{F}_r]]$, which is also denoted by $\varphi$. It is easy to see  by (2.1.8) that
 for $\varphi \in {\rm Aut}^{\rm fil}(\mathbb{Z}_l[[\frak{F}_r]]), \alpha \in \mathbb{Z}_l[[\frak{F}_r]]$, we have
 $$ \Theta(\varphi(\alpha)) = \varphi^*(\Theta(\alpha)). \leqno{(2.1.9)}$$
 \\
{\bf 2.2. $l$-adic Milnor invariants.} Let ${\rm Ih}_S : {\rm Gal}_{k} \rightarrow P(\frak{F}_r)$ be the Ihara representation associated to $S$ in (1.2.3). \\
\\
{\bf Lemma 2.2.1.} {\it  Let $g \in {\rm Gal}_{k}$.  For each $1\leq i \leq r$, there exists uniquely $y_i(g) \in \frak{F}_r$ satisfying the following properties}:\\
(1)   ${\rm Ih}_S(g)(x_i) = y_i(g) x_i^{\chi_l(g)} y_i(g)^{-1}$,  {\it where $\chi_l$ is the $l$-cyclotomic character}, \\
(2)   {\it In the expression $[y_i(g)] = c_1^{(i)}X_1 +  \cdots + c_r^{(i)} X_r$ $(c_j^{(i)} \in \mathbb{Z}_l)$, $c_i^{(i)} = 0.$}\\
\\
{\it Proof.}  Although the proof is standard, we give it for the sake of readers because this  lemma is 
basic in the theory of Milnor invariants.\\
Existence:  By the definition (1.2.1) of $P(\frak{F}_r)$ and Theorem 1.2.6 (1), there is $z_i \in \frak{F}_r$ such that ${\rm Ih}_S(g)(x_i) = z_i x_i^{\chi_l(g)} z_i^{-1}$ for each $i$. Let $[z_i] = a_1^{(i)}X_1 + \cdots + a_r^{(i)} X_r \, (a_j^{(i)} \in \mathbb{Z}_l)$. We set $y_i := z_i x_i^{- a_i^{(i)}}$. Then the conditions (1) and (2) are satisfied for $y_i$. \\
Uniqueness: Suppose that $y_i$ and $z_i$ in $\frak{F}_r$ satisfy the conditions (1) and (2). Since $z_i^{-1}y_i$ is in the centralizer of $x_i^{\chi_l(g)}$, $z_i^{-1}y_i = x_i^{b_i}$ for some $b_i \in \mathbb{Z}_l$. Comparing the coefficients of $X_i$ in $[z_i^{-1}y_i]$ and $[x_i^{b_i}]$, we have $b_i = 0$ and hence $y_i = z_i$. $\;\; \Box$\\
\\
We call $y_i(g) \in \frak{F}_l$ in Lemma 2.2.1  the $i$-th ({\em preferred}) {\em longitude} of $g \in {\rm Gal}_{k}$ for $S$. By Lemma 2.2.1, ${\rm Ih}_S(g)$ for $ g \in {\rm Gal}_{k}$ is determined by the $l$-cyclotomic value $\chi_l(g)$  and  the $r$-tuple ${\bf y}(g) := (y_1(g), \dots, y_r(g))$ of longitudes of $g$ for $S$. We note that ${\rm Ih}_S(g)$ acts on the abelianization $H$ of $\frak{F}_r$ by the multiplication by $\chi_l(g)$, $[{\rm Ih}_S(g)](X_i) = \chi_l(g) X_i$ for $1\leq i \leq r$. We also note that $y_i : {\rm Gal}_{k} \rightarrow \frak{F}_r$ is continuous, since ${\rm Ih}_S$ is continuous.

Following the case for pure braids ([Ko1], [Ko3; Chapter 1], [MK; Chapter 6, 4]), we will define the $l$-adic Milnor numbers of $g \in {\rm Gal}_{k}$ by the $l$-adic Magnus coefficients of the $i$-th longitude $y_i(g)$: Let $I = (i_1\cdots i_n)$ be a multi-index, where  $1\leq i_1, \dots , i_n \leq r$ and $|I| = n \geq 1$. The {\it $l$-adic Milnor number} of $g \in {\rm Gal}_{k}$ for $I$, denoted by $\mu(g; I) = \mu(g; i_1\cdots i_n)$, is defined by the $l$-adic Magnus coefficient of $y_{i_n}(g)$ for $I' := (i_1 \cdots i_{n-1})$:
$$  \mu(g; I) := \mu(I'; y_{i_n}(g)).  \leqno{(2.2.2)} $$
Here we set $\mu(g; I) := 0$ if $|I| = 1$. We note that the map $\mu(\; ; I) : {\rm Gal}_{k} \rightarrow \mathbb{Z}_l$ is continuous for each $I$, since $y_i :  {\rm Gal}_{k} \rightarrow \frak{F}_r$ is continuous. We define $\frak{a}(g)$ to be the ideal of $\mathbb{Z}_l$ generated by $\chi_l(g)-1$. Note that $\frak{a}(g) = 0$ when $g \in {\rm Gal}_{k(\zeta_{l^{\infty}})}$. We then define the {\em indeterminacy} $\Delta(g; I)$  by 
$$ \begin{array}{ll} \Delta(g; I) & := \mbox{the ideal of  $\mathbb{Z}_l$  generated by $\frak{a}(g)$ and $\mu(J;y_j(g))$, where $J$} \\
  & \;\; \;\;\; \; \mbox{ ranges over proper subsequence $I'$ and $j = i_n$ or $j$ is in $J$}
  \end{array} \leqno{(2.2.3)}
$$
We also write $\Delta(I'; y_{i_n}(g))$ for $\Delta(g; I)$ for the convenience later. We then set
$$ \overline{\mu}(g; I) := \mu(g; I) \; \mbox{mod}\; \Delta(g; I), \leqno{(2.2.4)} $$
which we call the {\it $l$-adic Milnor invariant} of $g \in {\rm Gal}_{k}$ for $I$. \\

We will show that the $l$-adic Milnor invariant $\overline{\mu}(g; I)$ for $g \in {\rm Gal}_{k}$ is unchanged when $g$ is replaced by its conjugate $hgh^{-1}$ for $h \in {\rm Gal}_{k(\zeta_{l^{\infty}})}$. To prove this, we prepare some lemmas. The formulas (1) and (2) of the next lemma was  proved by Wojtkowiak in terms of torsors of paths. See [W1; Proposition 1.0.7, Corollary 1.0.8 and Proposition 2.2.1].\\
\\
{\bf Lemma 2.2.5.} {\it For $g, h \in {\rm Gal}_{k}$ and $1\leq i \leq r$, we have}\\
(1)@$y_i(h^{-1}) = {\rm Ih}_S(h^{-1})(y_i(h)^{-1})$, \\
(2)@$y_i(hg) = {\rm Ih}_S(h)(y_i(g)) y_i(h)\;$ ({\it cocycle property}), \\
(3) @$y_i(hg h^{-1})   = {\rm Ih}_S(hg )(y_i(h^{-1})) {\rm Ih}_S(h)(y_i(g))y_i(h). $\\
\\
{\it Proof.} (1)  By Lemma 2.2.1, we have 
$$ \begin{array}{ll}
 x_i & = {\rm Ih}_S(h^{-1}) {\rm Ih}_S(h)(x_i) \\
& = {\rm Ih}_S(h^{-1})(y_i(h)x_i^{\chi_l(h)}y_i(h)^{-1}) \\  
&  = {\rm Ih}_S(h^{-1})(y_i(h))y_i(h^{-1})x_iy_i(h^{-1})^{-1} {\rm Ih}_S(h^{-1})(y_i(h)^{-1}),
\end{array}$$
from which we find ${\rm Ih}_S(h^{-1})(y_i(h))y_i(h^{-1}) = x_i^{a_i}$ for some $a_i \in \mathbb{Z}_l$. Passing to the abelianization $H$ of $\frak{F}_r$ and comparing the coefficients of $X_i$, we find $a_i = 0$ and hence  we obtain (1). \\
(2)  By Lemma 2.2.1, we have
$$ {\rm Ih}_S(hg)(x_i) = y_i(hg)x_i^{\chi_l(hg)}y_i(hg)^{-1}. \leqno{(2.2.5.1)}$$
On the other hand, we have
$$ \begin{array}{ll}
{\rm Ih}_S(h g)(x_i)  &
= {\rm Ih}_S(h) {\rm Ih}_S(g)(x_i) \\
&= {\rm Ih}_S(h)(y_i(g)x_i^{\chi_l(g)} y_i(g)^{-1})\\ 
&= {\rm Ih}_S(h)(y_i(g)) {\rm Ih}_S(h)(x_i^{\chi_l(g)}) {\rm Ih}_S(h)(y_i(g)^{-1}) \\
&= {\rm Ih}_S(h)(y_i(g))y_i(h)  x_i^{\chi_l(h g)}y_i(h)^{-1} {\rm Ih}_S(h)(y_i(g)^{-1}).
\end{array} \leqno{(2.2.5.2)}
$$
Comparing (2.2.5.1) and (2.2.5.2), we have $y_i(hg)^{-1} {\rm Ih}_S(h)(y_i(g))y_i(h) = x_i^{b_i}$ for some $b_i \in \mathbb{Z}_l$. Passing to the abelianization and comparing the coefficients of $X_i$, we find $b_i = 0$ and hence we obtain (2).\\
(3) By (2), we have
$$ y_i(hg h^{-1})   = {\rm Ih}_S(hg)(y_i(h^{-1}))y_i(hg) = {\rm Ih}_S(hg )(y_i(h^{-1})) {\rm Ih}_S(h)(y_i(g))y_i(h).  \;\; \Box$$
\\
For $\rho \in {\rm Gal}_k$ and a multi-index $J$ with $|J| \geq 1$, we define $\Theta_J(\rho)$ by
$$ \Theta_J(\rho) := {\rm Ih}_S(\rho)^*(X_J) - \chi_l(\rho)^{|J|} X_J. \leqno{(2.2.6)}$$
Since ${\rm Ih}_S(\rho)^*$ is  filtration-preserving, we note ${\rm deg}(\Theta_J(\rho)) \geq |J|$.\\
\\
{\bf Lemma 2.2.7.} {\em Notations being as above, the following assertions hold.}\\
(1) {\em $\Theta_J(\rho)$ is a Magnus power series $\sum_{|A|\geq |J|} m_A(J;\rho) X_A$ satisfying the following properties: }\\
$\;$ (i) {\it if $m_A(J;\rho) \neq 0$, then $A$ contains $J$ as a proper subsequence. So we may } \\
$\;\;$  {\em write} $\Theta_J(\rho) = \sum_{J \subsetneqq A} m_A(J;\rho) X_A.$\\
$\;$ (ii) {\it any coefficient $m_A(J;\rho)$ is a multiple of  $\mu(B;y_j(\rho))$ by an $l$-adic integer, \\
$\;\;$ where $B$ is some proper subsequence of $A$ and $j$ is in $J$.}\\
(2) {\em For $y \in \frak{F}_r$, we have}
$$ \begin{array}{ll} \Theta({\rm Ih}_S(\rho)(y))  & = \displaystyle{1 + \sum_{|J|\geq 1} \chi_l(\rho)^{|J|}\mu(J; y)X_J + \sum_{|J|\geq 1} \mu(J;y) \Theta_J(\rho)}\\
  &  \equiv \displaystyle{ \Theta(y) + \sum_{|J|\geq 1} \mu(J;y) \Theta_J(\rho)} \;\; \mbox{mod} \;\; \frak{a}(\rho).
\end{array}$$
\\
{\em Proof.}  (1)  Let $1\leq j \leq r$ and write $\Theta(y_j(\rho)) = 1 + Y_j(\rho)$. By (2.1.9) and Lemma 2.2.1, we have

$$ \begin{array}{ll}
 {\rm Ih}_S(\rho)^*(X_j) & = {\rm Ih}_S(\rho)^*(\Theta(x_j -1))\\
                                & = \Theta({\rm Ih}_S(\rho)(x_j -1))\\
                                & = \Theta(y_j(\rho)x_j^{\chi_l(\rho)} y_j(\rho)^{-1}) - 1 \\
                                & = \Theta(y_j(\rho)) \Theta(x_j)^{\chi_l(\rho)} \Theta(y_j(\rho)^{-1}) - 1\\
                                & = (1 + Y_j(\rho))(1+X_j)^{\chi_l(\rho)}(1 - Y_j(\rho) + Y_j(\rho)^2 - \cdots) - 1\\
                                & = \chi_l(\rho) X_j + \Theta_j(\rho),
                                \end{array} \leqno{(2.2.7.1)}
$$
where $\Theta_j(\rho)$ is the sum of terms of the form $u Y_j (\rho)^a X_j^b Y_j(\rho)^c$ for some $a,c \geq 0$ with $a+c\geq 1$, $b \geq 1$ and $u \in \mathbb{Z}_l$. Write $\Theta_j(\rho) = \sum_{|A| \geq 2} m_A(j;\rho) X_A$. 
Then it is easy to see that if $m_A(j; \rho) \neq 0$, then $A$ must contain $j$. Moreover, since $Y_j(\rho) = \sum_{|B| \geq 1} \mu(B; y_j(\rho)) X_B$, $m_A(j;\rho)$ is a multiple of $\mu(B; y_j(\rho))$ by an $l$-adic integer, where $B$ is some proper subsequence of $A$.  Let $J = (j_1\cdots j_n)$. By (2.2.7.1), we have
$$ \begin{array}{ll}
\sum_{|A| \geq |J|} m_A(J;\rho) X_A  & := \Theta_J(\rho)  \\
&  := {\rm Ih}_S(\rho)^*(X_J) - \chi_l(\rho)^{|J|}X_J \\
& = {\rm Ih}_S(\rho)^*(X_{j_1})\cdots {\rm Ih}_S(\rho)^*(X_{j_n}) - \chi_l(\rho)^{|J|} X_J\\
& = (\chi_l(\rho)X_{j_1} + \Theta_{j_1}(\rho))\cdots (\chi_l(\rho)X_{j_n} + \Theta_{j_n}(\rho)) - \chi_l(\rho)^{|J|} X_J\\
& = \Phi_{j_1}(\rho)\cdots \Phi_{j_n}(\rho),
\end{array}
$$
where $\Phi_j(\rho)$ is $\chi_l(\rho)X_j$ or $\Theta_j(\rho)$ and at least one $\Theta_j(\rho)$ is involved for some $j$. Hence, by the properties of coefficients of $\Theta_j(\rho) = \sum_{|A| \geq 2} m_A(j;\rho) X_A$ proved above, we obtain the properties (i) and (ii).\\
(2)  By  (2.1.9) and  (2.2.6), we have
$$ \begin{array}{ll}
\Theta({\rm Ih}_S(\rho)(y)) & = {\rm Ih}_S(\rho)^*(\Theta(y))\\
                                             & = \displaystyle{ {\rm Ih}_S(\rho)^*( 1 + \sum_{|J| \geq 1} \mu(J;y)X_J) } \\
                                             & = 1 + \displaystyle{ \sum_{|J| \geq 1} \mu(J;y) {\rm Ih}_S(\rho)^*(X_J) }\\
                                             & =  1 + \displaystyle{ \sum_{|J| \geq 1} \mu(J;y) (\chi_l(\rho)^{|J|}X_J + \Theta_J(\rho))  }\\
                                             & =  1 + \displaystyle{ \sum_{|J| \geq 1} \chi_l(\rho)^{|J|} \mu(J;y) X_J + \sum_{|J|\geq 1} \mu(J;y)\Theta_J(\rho)  } \\
                                              & \equiv \displaystyle{ \Theta(y) + \sum_{|J|\geq 1} \mu(J;y)\Theta_J(\rho)  } \;\; \mbox{mod} \;\; \frak{a}(\rho). \;\;\;\; \Box
\end{array}
$$
\\
We are ready to prove the following\\
\\
{\bf Theorem 2.2.8.} {\em For a multi-index $I$, the $l$-adic Milnor invariant $\overline{\mu}(g;I)$ for $g \in {\rm Gal}_k$ is unchanged when $g$ is replaced with its conjugate by an element of ${\rm Gal}_{k(\zeta_{l^{\infty}})}$. To be precise, let  $I$ be a multi-index with $|I| \geq 1$. Let $g \in {\rm Gal}_k$ and $h \in  {\rm Gal}_{k(\zeta_{l^{\infty}})}$. Then we have $\Delta(hgh^{-1}; I) = \Delta(g;I)$ and}
$$\overline{\mu}(hgh^{-1};I) = \overline{\mu}(g;I). $$
\\
{\it Proof.}  Let $I$ be a  multi-index with $|I| \geq 1$ and $1\leq i \leq r$. For $g, h \in {\rm Gal}_{k}$, 
we will show  
$$ \mu(I; y_i(hgh^{-1})) \equiv \mu(I; y_i(g)) \; \mbox{mod}\; \Delta(I; y_i(g)). \leqno{(2.2.8.1)}$$
By Lemma 2.2.5 (3), we have
$$ \Theta(y_i(hg h^{-1}))   = \Theta({\rm Ih}_S(hg )(y_i(h^{-1}))) \Theta({\rm Ih}_S(h)(y_i(g))) \Theta(y_i(h)). \leqno{(2.2.8.2)}$$
For simplicity, we set, for a multi-index $J$ with $|J| \geq 1$, 
$$a_J := \mu(J; {\rm Ih}_S(hg )(y_i(h^{-1}))), b_J := \mu(J; {\rm Ih}_S(h)(y_i(g))), c_J := \mu(J; y_i(h)).$$
Then, from (2.2.8.2) or  (2.1.4), we have
$$\begin{array}{l} \mu(I; y_i(hgh^{-1})) \\
\;\;= \displaystyle{a_I + b_I + c_I + \sum_{AB=I} a_A b_B  + \sum_{BC=I} b_B c_C + \sum_{AC=I} a_A c_C + \sum_{ABC=I}a_A b_B c_C, }
\end{array}\leqno{(2.2.8.3)}$$
where $A, B, C$ are multi-indices with $|A|, |B|, |C| \geq 1$. 

First, we  look at $b_B$ for a subsequence $B$ of $I$.  By Lemma 2.2.7 (1), (2) and as $h \in {\rm Gal}_{k(\zeta_{l^{\infty}})}$,  we have
$$ b_B  = \mu(B; y_i(g)) + \mu(J;y_i(g))(\mbox{an $l$-adic integer})  $$
for some proper subsequence $J$ of $B$. Therefore,  by (2.2.8.3) and the definition of $\Delta(I;y_i(g))$, we have
$$  \mu(I; y_i(hgh^{-1})) - \mu(I;y_i(g)) \equiv a_I +  c_I + \sum_{AC=I} a_A c_C \;\;\mbox{mod}\;  \Delta(I; y_i(g)).   \leqno{(2.2.8.4)}$$
Here we note that the right hand side of (2.2.8.4) is the coefficient of $X_I$ of $\Theta({\rm Ih}_S(hg )(y_i(h^{-1})))\Theta(y_i(h))$.  

So, next, we look at $\Theta({\rm Ih}_S(hg )(y_i(h^{-1})))\Theta(y_i(h))$. By (2.1.9), Lemma 2.2.5 (1) and Lemma 2.2.7 (2), we have
$$ \begin{array}{ll}
\Theta({\rm Ih}_S(hg )(y_i(h^{-1}))) & = {\rm Ih}_S(hg)^*(\Theta(y_i(h^{-1})) )\\
                                                & = {\rm Ih}_S(h)^* {\rm Ih}_S(g)^*(\Theta(y_i(h^{-1})))\\
                                                & \equiv \displaystyle{ {\rm Ih}_S(h)^*(\Theta(y_i(h^{-1})) + \sum_{|J| \geq 1} \mu(J; y_i(h^{-1})) \Theta_J(g))}\;\; (\mbox{mod}\; \frak{a}(g))\\
                                                & = \displaystyle{\Theta({\rm Ih}_S(h)(y_i(h^{-1}))) + \sum_{|J| \geq 1} \mu(J; y_i(h^{-1})) {\rm Ih}_S(h)^*(\Theta_J(g))}\\
                                                  & =   \displaystyle{\Theta(y_i(h)^{-1}) +   \sum_{|J| \geq 1} \mu(J; y_i(h^{-1})){\rm Ih}_S(h)^*(\Theta_J(g))}.                
\end{array}
\leqno{(2.2.8.5)}
$$ 
Here let us write $\Theta_J(g) = \sum_{J \subsetneqq A} m_A(J;g) X_A$ as in Lemma 2.2.7 (1). Then we have, as $h \in {\rm Gal}_{k(\zeta_{l^{\infty}})}$, 
$$ \begin{array}{ll} {\rm Ih}_S(h)^*(\Theta_J(g)) & = \displaystyle{\sum_{J \subsetneqq A}m_A(J;g) {\rm Ih}_S(h)^*(X_A) }\\
 & = \displaystyle{\sum_{J \subsetneqq A}m_A(J;g) (X_A + \Theta_A(h))} \;\; (\mbox{mod}\; \frak{a}(g))\\
& = \displaystyle{\sum_{J \subsetneqq A}m_A(J;g) (X_A + \sum_{A \subsetneqq A'} m_{A'}(A;h)X_{A'}).} 
\end{array} \leqno{(2.2.8.6)}$$
By (2.2.8.5) and (2.2.8.6), we have
$$ \begin{array}{l} \Theta( {\rm Ih}_S(hg )(y_i(h^{-1}))) \\
\equiv  \displaystyle{ \Theta(y_i(h)^{-1}) +   \sum_{|J| \geq 1} \sum_{J \subsetneqq A} \mu(J; y_i(h^{-1})) m_A(J;g) (X_A + \sum_{A \subsetneqq A'} m_{A'}(A;h)X_{A'})} \\
 \hspace{11cm} \mbox{mod}\; \frak{a}(g) \end{array} $$
and hence
$$\begin{array}{l}
\Theta( {\rm Ih}_S(hg )(y_i(h^{-1}))\Theta(y_i(h)) \\
 \equiv \displaystyle{ 1 + \sum_{|J| \geq 1}\sum_{J \subsetneqq A} \mu(J; y_i(h^{-1}))m_A(J;g)(X_A +\sum_{A \subsetneqq A'} m_{A'}(A;h)X_{A'})\Theta(y_i(h)) } \\   \hspace{11cm} \mbox{mod}\;\; \frak{a}(g). \end{array} \leqno{(2.2.8.7)}$$
Here we note by Lemma 2.2.7 (2) that $m_A(J;g)$ is a multiple of $\mu(B;y_j(g))$ by an $l$-adic integer for some proper subsequence $B$ of $A$ and $j$ in $J$. By the definition (2.2.3) of $\Delta(I;y_i(g))$,  the coefficient of $X_I$ in the right hand side of (2.2.8.7) must be congruent to 0 mod $\Delta(I; y_i(g))$. By (2.2.8.4), we obtain (2.2.8.1). 

Finally, we show that $\Delta(I;y_i(hgh^{-1})) = \Delta(I;y_i(g))$ by  induction on $|I|$. When $|I| = 1$, this is obviously true (both sides are $\frak{a}(g)=\frak{a}(hgh^{-1})$) by the definition. Assume that
$\Delta(I;y_i(hgh^{-1})) = \Delta(I;y_i(g)) $ for all $I$ with $|I|\leq n$ ($n \geq 1$).  Then, by (2.2.8.1), we have, for all $I$ with $|I| \leq n$ and $1\leq i\leq r$, 
$$ \mu(I; y_i(hgh^{-1})) \equiv \mu(I; y_i(g))  \; \mbox{mod} \; \Delta(I; y_i(g)) (=\Delta(I;y_i(hgh^{-1})). \leqno{(2.2.8.8)}$$ 
Using (2.2.8.8) and the definition (2.2.3) of $\Delta(I;y_i(\rho))$ for $\rho = hgh^{-1}, g$, we have $\Delta(I;y_i(hgh^{-1})) = \Delta(I;y_i(g))$ for $I$ with $|I| = n+1$. $\;\; \Box$
\\
\\
{\bf Remark 2.2.9.} It is known that a braid $\beta$ and its conjugate $\gamma \beta \gamma^{-1}$ give rise to the same link as their closures ($\beta \mapsto \gamma \beta \gamma^{-1}$ is one of Markov's transforms. cf.  [B; 2.2], [MK; Chapter 9]). In particular, they have the same Milnor invariants. So Theorem 2.2.8 may be seen as an arithmetic analogue of this known fact for braids.\\
\\
As a property of $l$-adic Milnor invariants, we have the following shuffle relation.\\
\\
{\bf Proposition 2.2.10.} {\em Let $g \in {\rm Gal}_{k}$. For multi-indices $I, J$ with $|I|, |J| \geq 1$ and $1\leq i \leq r$, we have}
$$ \sum_{H \in {\rm PSh}(I,J)} \overline{\mu}(g; Hi)  \equiv 0 \; \mbox{mod} \; \mbox{g.c.d} \{ \Delta(Hi) \; | \; H \in {\rm PSh}(I,J)\},$$
{\it where ${\rm PSh}(I,J)$ denotes the set of results of all proper shuffles of $I$ and $J$} ([CFL]). \\
\\
{\em Proof.} By (2.1.5), we have
$$ \mu(g; Ii)\mu(g;Ji) = \sum_{A \in {\rm Sh}(I,J)} \mu(g;Ai).$$
Taking mod $\mbox{g.c.d} \{ \Delta(Hi) \; | \; H \in {\rm PSh}(I,J)\}$, the left hand side is congruent to $0$ and any term $\mu(g;Ai)$ with $A \notin {\rm PSh}(I,J)$ is also congruent to $0$. So the assertion follows. $\;\; \Box$\\

Let $R_{S}^{\infty}$ be the set of primes of $k(\zeta_{l^{\infty}})$ lying over $R_S$ in (1.2.5).  For $\frak{p}_{\infty} \in R_S^{\infty}$, choose  a prime $\frak{P}$ of $\Omega_S$ lying over $\frak{p}_{\infty}$. Since $\frak{P}$ is unramified in the Galois extension $\Omega_S/k$ by Theorem 1.2.6 (2), we have the Frobenius automorphism $\sigma_{\frak{P}} \in {\rm Gal}(\Omega_S/k)$ of $\frak{P}$. By Theorem 2.2.8,
$\overline{\mu}(\sigma_{\frak{P}}; I)$ is independent of the choice of $\frak{P}$ lying over $\frak{p}_{\infty}$. So we define the {\it $l$-adic Milnor invariant} of $\frak{p}_{\infty}$ for a multi-index $I$ by
$$ \overline{\mu}(\frak{p}_{\infty}; I) := \overline{\mu}(\sigma_{\frak{P}}; I). \leqno{(2.2.11)}$$
We also set $\Delta(\frak{p}_{\infty}; I) := \Delta(\sigma_{\frak{P}}; I)$ so that $\overline{\mu}(\frak{p}_{\infty}; I) \in \mathbb{Z}_l/\Delta(\frak{p}_{\infty};I)$. 
Let $\frak{p}$ be the prime of $k$  lying below $\frak{p}_{\infty}$. Since $\chi_l(\sigma_{\frak{P}}) = {\rm N}\frak{p}$ (the norm of $\frak{p}$), in order to have $\mathbb{Z}_l/\Delta(\frak{p}_{\infty};I) \neq 0$, it is necessary that primes $\frak{p}_{\infty}$ in $R_S^{\infty}$ lie over 
$$ R_{S}^1 := \{ \frak{p} \in R_S \; | \; {\rm N}\frak{p} \equiv 1\; \mbox{mod} \;  l \}.$$
For $\frak{p} \in R_{S}^1$, let $e(\frak{p})$  denote the maximal integer such that 
$${\rm N}\frak{p} \equiv 1\; \mbox{mod} \;  l^{e(\frak{p})}.$$ 
It means that $\frak{p}$ is completely decomposed in $k(\zeta_{l^{e(\frak{p})}})/k$ and any prime of $k(\zeta_{l^{e(\frak{p})}})$ lying over $\frak{p}$ is inert in $k(\zeta_{l^{\infty}})/k(\zeta_{l^{e(\frak{p})}})$. Hence $\sigma_{\frak{P}} \in {\rm Gal}(\Omega_S/k(\zeta_{l^{e(\frak{p})}}))$. Then the indeterminacy $\Delta(\frak{p}_{\infty}; I)$  is an ideal of $\mathbb{Z}/l^{e(\frak{p})}\mathbb{Z}$.  We note that if $\mu(\sigma_{\frak{P}}; I) \equiv 0$ mod $l^{e(\frak{p})}$ for all $|I| \leq n$, then $ \overline{\mu}(\frak{p}_{\infty}; I)$ is well defined in $\mathbb{Z}/l^{e(\frak{p})}\mathbb{Z}$ for $|I| = n+1$. \\
\\
{\bf Remark 2.2.12.} In [Ms1] and [Ms2; Chapter 8], the arithmetic Milnor invariants for certain primes of a number field were introduced as multiple generalizations of power residue symbols and the R\'{e}dei triple symbol ([R]). See also [Am]. They are arithmetic analogues for primes of Milnor invariants of links. It is known ([Ko1], [Ko3; Chaper 1]) that Milnor invariants for a pure braid  coincide with those for the link obtained by closing the pure braid. Recently, we found a relation between $l$-adic Milnor invariants, Wojtkowiak's $l$-adic iterated integrals and $l$-adic polylogarithms ([NW], [W1]$\sim$[W4]) and multiple power residue symbols (in particular, R\'{e}dei symbols), which will be discussed in the forthcoming paper.\\

Finally, we introduce a filtration on ${\rm Gal}_{k}$ using $l$-adic Milnor numbers. We set ${\rm Gal}_{k}^{\rm Mil}[0] := {\rm Gal}_{k}$. For each integer $n \geq 1$, we define a subset ${\rm Gal}_{k}^{\rm Mil}[n]$ of 
${\rm Gal}_{k}$ by
$$ \begin{array}{ll}
{\rm Gal}_{k}^{\rm Mil}[n] & := \{ g \in {\rm Gal}_{k(\zeta_{l^{\infty}})} \, | \, \mu(g ; I) = 0 \; \mbox{for} \; |I| \leq  n \}\\
                                       & = \{  g \in {\rm Gal}_{k(\zeta_{l^{\infty}})} \, | \, {\rm deg}(\Theta(y_i(g)) -1) \geq n \; \mbox{for}\; 1\leq i \leq r \}.
\end{array}                                       
\leqno{(2.2.13)}$$
We then have the descending series 
$$  {\rm Gal}_{k} =  {\rm Gal}_{k}^{\rm Mil}[0]  \supset  {\rm Gal}_{k}^{\rm Mil}[1] \supset \cdots \supset  {\rm Gal}_{k}^{\rm Mil}[n]  \supset \cdots $$
and we call it  the {\it Milnor filtration} of  ${\rm Gal}_{k}$. \\
\\
{\bf Proposition 2.2.14.} {\it For $n \geq 0$, ${\rm Gal}_{k}^{\rm Mil}[n]$ is a closed normal subgroup of ${\rm Gal}_{k}$.}\\
\\
{\it Proof.}  This proposition is an immediate consequence of the coincidence of the Milnor filtration and the Johnson filtration which will be proved in Proposition 3.3.3 in Section 3.
So we give herewith a direct and brief proof.\\
We may assume $n \geq 1$. Since $\mu(\; ;I) : {\rm Gal}_{k} \rightarrow \mathbb{Z}_l$ is continuous for each $I$ and ${\rm Gal}_{k}^{\rm Mil}[n] = \bigcap_{|I|\leq n} \mu(\; ; I)^{-1}(0)$, ${\rm Gal}_{k}^{\rm Mil}[n]$ is closed in ${\rm Gal}_{k}$. Let $g, h \in {\rm Gal}_{k}^{\rm Mil}[n]$ and so ${\rm deg}(\Theta(y_i(\rho)) - 1) \geq n$ for $\rho = g, h$ and each $1\leq i \leq r$. Then we can show easily ${\rm deg}(\Theta(y_i(g^{-1})) - 1) \geq n$,  ${\rm deg}(\Theta(y_i(gh)) - 1) \geq n$ and ${\rm deg}(\Theta(y_i(hgh^{=1})) - 1) \geq n$ by using Lemma 2.2.5 (1), (2) and (3), respectively. $\;\; \Box$\\
\\
{\bf 2.3. Pro-$l$ link groups and Massey products.} Following the analogy with the link group of a pure braid link ([Ar],[B; Theorem 2.2]), we define the {\em pro-$l$ link group} of each Galois element $g \in {\rm Gal}_{k}$ associated to ${\rm Ih}_S$ by
$$ \begin{array}{ll} \Pi_S(g) & := \langle x_1, \dots , x_r \, | \, y_1(g)x_1^{\chi_l(g)}y_1(g)^{-1}=x_1, \cdots , y_r(g)x_r^{\chi_l(g)}y_r(g)^{-1}= x_r  \rangle \\
                                       & =   \langle x_1, \dots , x_r \, | \,  x_1^{1-\chi_l(g)}[x_1^{-1},y_1(g)^{-1}] = \cdots = x_r^{1-\chi_l(g)}[x_r^{-1},y_r(g)^{-1}] = 1 \rangle \\  
                                      & := \frak{F}_r/\frak{N}_S(g), 
\end{array} \leqno{(2.3.1)}$$
where $\frak{N}_S(g)$ denotes the closed subgroup of $\frak{F}_r$ generated normally by the pro-$l$ words $x_1^{1-\chi_l(g)}[x_1^{-1},y_1(g)^{-1}], \dots , x_r^{1-\chi_l(g)}[x_r^{-1},y_r(g)^{-1}] $. We will give a cohomological interpretation of  $l$-adic Milnor invariants of $g \in {\rm Gal}_{k}$ by Massey products in the cohomology of the pro-$l$ link group $\Pi_S(g)$. In the following, we let $g \in {\rm Gal}_{k}$ and $\frak{a}$ an ideal of $\mathbb{Z}_l$  such that  $\frak{a} \neq \mathbb{Z}_l$ and $\chi_l(g) \equiv 1$ mod $\frak{a}$. We may write $\frak{a} = l^a\mathbb{Z}_l$ for some $1 \leq  a \leq \infty$ ($l^a := 0$ if $a = \infty$).  When $g \in {\rm Gal}_{k(\zeta_l^{\infty})}$, we have $a = \infty$ and $\frak{a} = 0$.

Let $C^i(\Pi_S(g), \mathbb{Z}_l/\frak{a})$ be the $\mathbb{Z}_l/\frak{a}$-module of continuous $i$-cochains $(i \geq 0)$ of $\Pi_S(g)$ with coefficients in $\mathbb{Z}_l/\frak{a}$, where $\Pi_S(g)$ acts on $\mathbb{Z}_l/\frak{a}$ trivially.  We consider the differential graded algebra $(C^{\bullet}(\Pi_S(g), \mathbb{Z}_l/\frak{a}), d)$, where the product on $C^{\bullet}(\Pi_S(g), \mathbb{Z}_l/\frak{a}) = \bigoplus_{i \geq 0} C^i(\Pi_S(g), \mathbb{Z}_l/\frak{a})$ is given by the cup product and the differential $d$ is the coboundary operator. Then we have the cohomology ring $H^*(\Pi_S(g), \mathbb{Z}_l/\frak{a}) := \bigoplus_{i \geq 0} H^i(C^{\bullet}(\Pi_S(g), \mathbb{Z}_l/\frak{a}))$ of the pro-$l$ group $\Pi_S(g)$ with coefficients in $\mathbb{Z}_l/\frak{a}$. In the following, we deal with only one and two dimensional cohomology groups. For the sign convention, we follow [Dw]. For $c_1, \dots , c_n \in H^1(\Pi_S(g), \mathbb{Z}_l/\frak{a})$, an $n$-th {\it Massey product}  $\langle c_1,\dots, c_n \rangle$ is said to be {\it defined} if there is an array 
$$ W = \{ w_{ij} \in C^1(\Pi_S(g), \mathbb{Z}_l/\frak{a}) \; | \; 1 \leq i < j \leq n+1, (i,j) \neq (1,n+1) \}$$
such that
$$\left\{ 
\begin{array}{l}
[\omega_{i,i+1}] = c_i \;\; (1\leq i \leq n),\\
 \displaystyle{dw_{ij} = \sum_{a=i+1}^{j-1} w_{ia}\cup w_{aj}} \;\; (j \neq i+1).

\end{array}
\right.
$$
Such an array $W$ is called a {\it defining system} for $\langle c_1,\dots, c_n \rangle$. The value of  $\langle c_1,\dots, c_n \rangle$ relative to $W$ is defined by the cohomology class represented by the $2$-cocycle 
$$  \sum_{a=2}^n w_{1a} \cup w_{a,n+1},$$
and denoted by $\langle c_1,\dots, c_n \rangle_{W}$.  A Massey product $\langle c_1,\dots, c_n \rangle $ itself is taken to be the subset of $H^2(\Pi_S(g), \mathbb{Z}_l/\frak{a})$ consisting of elements $\langle c_1,\dots, c_n \rangle_{W}$ for some defining system $W$. By convention, $\langle c \rangle = 0$. The following lemma is a baisc fact ([Kr]). \\
\\
{\bf Lemma 2.3.2.} {\em We have  $\langle c_1, c_2 \rangle = c_1 \cup c_2$. For $n \geq 3$, $\langle c_1,\dots, c_n \rangle$ is defined and consists of a single element if $\langle c_{j_1},\dots, c_{j_a} \rangle = 0$ for all proper subsets $\{j_1,\dots, j_a \}$ $(a \geq 2)$ of $\{ 1,\dots ,n\}$.}\\

Next, we recall a relation between Massey products and the Magnus coefficients for our situation. Let $\psi : \frak{F}_r \rightarrow \Pi_S(g) = \frak{F}_r/\frak{N}_S(g)$ be the natural homomorphism. We denote by $\gamma_i$ the image of $x_i$ under $\psi$, $\gamma_i := x_i$ mod $\frak{N}_S(g)$, for $1\leq i \leq r$. By the definition (2.3.1) of $\Pi_S(g)$ and our assumption, $\pi$ induces the isomorphism $\frak{F}_r/\frak{F}_r^{l^a}\frak{F}_r(2) \stackrel{\sim}{\rightarrow} \Pi_S(g)/\Pi_S(g)^{l^a}[\Pi_S(g),\Pi_S(g)] \simeq (\mathbb{Z}_l/\frak{a})^{\oplus r}$ and so we have the isomorphism $H^1(\Pi_S(g), \mathbb{Z}_l/\frak{a}) \simeq H^1(\frak{F}_r,\mathbb{Z}_l/\frak{a})$. Therefore the Hochschild-Serre spectral sequence yields the isomorphism
$${\rm tg} : H^1(\frak{N}_S(g), \mathbb{Z}_l/\frak{a})^{\Pi_S(g)} \rightarrow H^2(\Pi_S(g),\mathbb{Z}_l/\frak{a}).$$
Here ${\rm tg}$ is  the transgression defined as follows. For $a \in H^1(\frak{N}_S(g), \mathbb{Z}_l/\frak{a})^{\Pi_S(g)}$, choose a 1-cochain $b \in C^1(\frak{F}_r, \mathbb{Z}_l/\frak{a})$ such that $b|_{\frak{N}_S(g)} = a$. Since the value $db(f_1,f_2)$, $f_i \in \frak{F}_r$, depends only on the cosets $f_i$ mod $\frak{N}_S(g)$, there is a 2-cocyle $c \in Z^2(\Pi_S(g), \mathbb{Z}_l/\frak{a})$ such that $\psi^*(c) = db$. Then ${\rm tg}(a)$ is defined to be the class of $c$. The dual to ${\rm tg}$ is  called the Hopf isomorphism:
$${\rm tg}^{\vee} : H_2(\Pi_S(g), \mathbb{Z}_l/\frak{a}) \stackrel{\sim}{\rightarrow} H_1(\frak{N}_S(g), \mathbb{Z}_l/\frak{a})_{\Pi_S(g)} = \frak{N}_S(g)/\frak{N}_S(g)^{l^a}[\frak{N}_S(g), \frak{F}_r].$$
 Then we have the following proposition (cf. [St; Lemma 1.5], [Ms1, Theorem 2.2.2]). \\
\\
{\bf Proposition 2.3.3.} {\it Notations being as above, let $c,\dots, c_n \in H^1(\Pi_S(g), \mathbb{Z}_l/\frak{a})$ and  $W = (w_{ij})$ a defining system for the Massey product $\langle c_1,\dots, c_n \rangle$. Let $ f \in \frak{N}_S(g)$ and set $\frak{r} := ({\rm tg}^{\vee})^{-1}(f \; {\rm mod}\;  \frak{N}_S(g)^{l^a}[\frak{N}_S(g),\frak{F}_r])$. Then we have} 
$$ \begin{array}{l} \langle c_1,\dots, c_n \rangle_{W}(\frak{r})\\
 = \displaystyle{\sum_{j=1}^n (-1)^{ j+1} \sum_{e_1+\cdots + e_j=n} \sum_{I =(i_1\cdots i_j)} w_{1,1+e_1}(\gamma_{i_1})\cdots w_{n+1-e_j, n+1}(\gamma_{i_j})\mu(I; f)_{\frak{a}},}
\end{array}$$
{\it where  $e_1,\dots,e_j$ run over positive integers satisfying $e_1+\cdots + e_j = n$ and  $\mu(I; f)_{\frak{a}} := \mu(I;f)$ mod $\frak{a}$. }\\

 Now, let $\gamma_1^*, \dots , \gamma_{r}^* \in H^1(\Pi_S(g), \mathbb{Z}_l/\frak{a})$ be the Kronecker dual to $\gamma_1, \dots , \gamma_{r}$, namely, $\gamma_i^*(\gamma_j) = \delta_{ij}$ for $1\leq i, j \leq r$.  Let $\frak{r}_i := ({\rm tg}^{\vee})^{-1}(x_i^{1-\chi_l(g)}[x_i^{-1}, y_i(g)^{-1}])$ mod  $[\frak{N}_S(g),\frak{F}_r])$ for $1\leq i \leq r$.
 Let $I = (i_1\cdots i_n)$ be a multi-index such that  $|I| = n \geq 2$. Let $g \in {\rm Gal}_{k}$. We assume the following conditions:
 $$ \left\{
\begin{array}{l}
(1) \;  \mu((j_1\cdots j_a); x_i^{1-\chi_l(g)}) \equiv 0 \; \mbox{mod} \; \frak{a} \; \mbox{for any subset}\\
\;\;\;\; \; \{ j_1, \dots , j_a\} \; \mbox{of} \; \{i_1,\dots, i_n\} \; \mbox{and} \; 1\leq i \leq r, \\
(2) \; i_1, \dots, i_n\; \mbox{are distinct each other, and}\\
\;\;\;\;\;   \mu(g; (j_1\cdots j_a)) \equiv  0 \; \mbox{mod}\; \frak{a} \; \mbox{for any proper subset}\\
\;\;\;\;  \; \{ j_1, \dots , j_a\} \; \mbox{of} \; \{i_1,\dots, i_n\}.\\
\end{array}\right.  \leqno{(2.3.4)}
$$
We note that the condition (1) is unnecessary when $g \in {\rm Gal}_{k(\zeta_{l^\infty})}$. The following theorem gives a cohomological interpretation of $\mu(g;I)_{\frak{a}} := \mu(g;I)$ mod $\frak{a}$ by the Massey product in the cohomology of $\Pi_S(g)$.\\
\\
{\bf Theorem 2.3.5.}  {\em Notations and assumtions being as above, the Massey product $\langle \gamma_{i_1}^*, \dots , \gamma_{i_n}^* \rangle$ in $H^2(\Pi_S(g),\mathbb{Z}_l/\frak{a})$ is uniquely defined and we have}
$$  \mu(g;I)_{\frak{a}} = (-1)^n \langle \gamma_{i_1}^*,\dots , \gamma_{i_n}^* \rangle(\frak{r}_{i_n}).$$
{\em Proof.}  First, we compute $\mu(J; x_i^{1-\chi_l(g)}[x_i^{-1},y_i(g)^{-1}])$ for a multi-index $J = (j_1\cdots j_a)$, where $\{j_1, \dots , j_a\}$ is a subset of $\{ i_1, \dots , i_n \}$. We note that
 $$\begin{array}{l}
\Theta(x_i^{1-\chi_l(g)}[x_i^{-1},y_i(g)^{-1}]) \\
\;\;\;\;\; = \Theta(x_i^{1-\chi_l(g)})(1 + \Theta(x_i^{-1})\Theta(y_i(g)^{-1})(\Theta(x_iy_i(g)) - \Theta(y_i(g)x_i))).
\end{array}$$
By our assumption (2.3.4) (1), we have
$$\begin{array}{ll} \mu(J; x_i^{1-\chi_l(g)}[x_i^{-1},y_i(g)^{-1}]) & \equiv \mu(J; x_iy_i(g)) - \mu(J; y_i(g)x_i) \\
 & \;\;\;\;\; \; \displaystyle{ + \sum_{A}(\mu(A; x_iy_i(g)) - \mu(A; y_i(g)x_i))\nu_A} \; \mbox{mod}\; \frak{a}, \end{array}  \leqno{(2.3.5.1)}$$
where $A$ runs over some proper subsequences of $J$ and $\nu_A \in \mathbb{Z}_l$. By the straightforward computation, we have
$$ \mu(J; x_iy_i(g)) = \left\{ \begin{array}{ll} 
\mu(g; (Ji)) \;\; &  (i \neq j_1),\\
\mu(g; (Jj_1)) + \mu(g; (j_2\cdots j_a j_1)) \;\; & (i = j_1),
\end{array} \right.
$$
and
$$
\mu(J; y_i(g)x_i) = \left\{ \begin{array}{ll} 
\mu(g; (Ji)) \;\; & (i \neq j_a),\\
\mu(g; (Jj_a)) + \mu(g; J) \;\; & (i = j_a).
\end{array} \right.
$$
Hence we have
$$ \mu(J; x_iy_i(g)) - \mu(J; y_i(g)x_i) = 
\left\{ \begin{array}{ll}
\mu(g;( j_2\cdots j_aj_1)) - \delta_{j_1, j_a}\mu(g; J) \; & ( i = j_1),\\
\mu(g; (j_2\cdots j_a j_1))\delta_{j_1,j_a} - \mu(g; J) \;&  (i = j_a),\\
0 \; & \mbox{(otherwise)}.
\end{array}\right.
\leqno{(2.3.5.2)}
$$
Now, let $n = 2$. Then we have $\langle \gamma_{i_1}^*, \gamma_{i_2}^* \rangle = \gamma_{i_1}^* \cup \gamma_{i_2}^*$. By Proposition 2.3.3, (2.3.4) (2), (2.3.5.1) and (2.3.5.2), we have 
$$\langle \gamma_{i_1}^*, \gamma_{i_2}^* \rangle (\frak{r}_{i_2}) = -\mu(I; [x_{i_2},y_{i_2}(g)])_{\frak{a}} = \mu(g;I)_{\frak{a}}. $$
Suppose $ n \geq 3$ and let $\{j_1, \dots, j_a \}$  be a proper subset of $\{i_1, \dots , i_n \}$. Then, by our assumption (2.3.4) (2), (2.3.5.1) and (2.3.5.2), we have
$$ \mu(J; x_i^{1-\chi_l(g)}[x_i^{-1},y_i(g)^{-1}]) \equiv  0 \; \mbox{mod} \; \frak{a}$$
for $J = (j_1\cdots j_a)$ and $1\leq i \leq r$. So,  by Proposition 2.3.4, we have
$$ \langle \gamma_{j_1}^*, \dots , \gamma_{j_a}^* \rangle (\frak{r}_i) = 0$$
for $1\leq i \leq r$. Since $H_2(\Pi(g),\mathbb{Z}_l/\frak{a})$ is generated by $x_i^{1-\chi_l(g)}[x_i,y_i(g)]$ for $1\leq i \leq r$, we have
$$ \langle c_{j_1}, \dots , c_{j_a} \rangle = 0.$$
Therefore, by Lemma 2.3.2, the Massey product $\langle c_{i_1}, \dots , c_{i_n} \rangle $ is uniquely defined. By Proposition 2.3.3, (2.3.4) (2), (2.3.5.1) and (2.3.5.2) again, we have
$$\langle \gamma_{i_1}^*, \dots , \gamma_{i_n}^* \rangle (\frak{r}_{i_n}) =(-1)^{n+1}\mu(I;x_n^{1-\chi_l(g)}[x_{i_n},y_{i_n}(g)])_{\frak{a}}= (-1)^n \mu(g; I)_{\frak{a}}.\;\; \Box
$$
\vspace{.2cm}

\begin{center}
{\bf 3.  Pro-$l$ Johnson homomorphisms}
\end{center}
\vspace{.2cm}
{\bf 3.1. Some algebras associated to lower central series.} For each integer $n \geq 1$, we let
$${\rm gr}_n(\frak{F}_r) := \frak{F}_r(n)/\frak{F}_r(n+1),$$
which is a free $\mathbb{Z}_l$-module whose rank $\ell_r(n)$ is given by the Witt formula ([MKS; 5.6, Theorem 5.11], [Se; Ch. IV, 4, 6]):
$$ \ell_r(n) = \frac{1}{n} \sum_{d|n} \mu(d) r^{n/d},$$
where $\mu(d)$ is the M\"{o}bius function.  The graded $\mathbb{Z}_l$-module
$$ {\rm gr}(\frak{F}_r) := \bigoplus_{n\geq 1} {\rm gr}_n(\frak{F}_r)$$
has the structure of a graded free Lie algebra over $\mathbb{Z}_l$: For $a = s$ mod $\frak{F}_r(m+1) \in {\rm gr}_m(\frak{F}_r)$ and $b = t$ mod $\frak{F}_r(n+1) \in {\rm gr}_n(\frak{F}_r)$ $(s \in  \frak{F}_r(m), t \in \frak{F}_r(n))$, the Lie bracket on ${\rm gr}(\frak{F}_r)$ is defined by
$$ [a, b]  := [s,t] \; {\rm mod} \; \frak{F}_r(m+n+1).$$

We consider the graded associative algebra over $\mathbb{Z}_l$ defined by
$$ {\rm gr}(\mathbb{Z}_l[[\frak{F}_r]]) := \bigoplus_{n \geq 0} {\rm gr}_n(\mathbb{Z}_l[[\frak{F}_r]]), \;\; {\rm gr}_n(\mathbb{Z}_l[[\frak{F}_r]]) := I_{\mathbb{Z}_l[[\frak{F}_r]]}^n/I_{\mathbb{Z}_l[[\frak{F}_r]]}^{n+1}.$$
The map $f \mapsto f-1$ $(f \in \frak{F}_r(n))$ defines an injective  $\mathbb{Z}_l$-linear map
$$ {\rm gr}_n(\frak{F}_r) \hookrightarrow {\rm gr}_n(\mathbb{Z}_l[[\frak{F}_r]])  \leqno{(3.1.1)} $$
for $n \geq 1$ and the injective Lie algebra homomorphism over $\mathbb{Z}_l$
$$  {\rm gr}(\frak{F}_r) \hookrightarrow {\rm gr}(\mathbb{Z}_l[[\frak{F}_r]]),$$
where ${\rm gr}(\mathbb{Z}_l[[\frak{F}_r]])$ is shown to be the universal enveloping algebra of the Lie algebra ${\rm gr}(\frak{F}_r)$. 
Moreover, by the correspondence $x_i -1$ mod $I_{\mathbb{Z}_l[[\frak{F}_r]]}^2 \in {\rm gr}_1(\mathbb{Z}_l[[\frak{F}_r]]) \mapsto X_i \in H$,   
 we have the isomorphism of $\mathbb{Z}_l$-modules
$$ \Theta_n : \;  {\rm gr}_n(\mathbb{Z}_l[[\frak{F}_r]]) \; \simeq  \; H^{\otimes n} \leqno{(3.1.2)}$$
for each $n \geq 0$ and so ${\rm gr}(\mathbb{Z}_l[[\frak{F}_r]])$ is identified with the tensor algebra $T(H)$:
$${\rm gr}(\mathbb{Z}_l[[\frak{F}_r]])  =  T(H)= \mathbb{Z}_l \langle X_1,\dots , X_r \rangle.$$
The composition of the map of (3.1.1) with $\Theta_n$ in (3.1.2), denoted also by $\Theta_n : {\rm gr}_n(\frak{F}_r) \hookrightarrow H^{\otimes n}$, is the degree $n$ part of the pro-$l$ Magnus embedding in (2.1.3):
$$ \Theta_n = (\Theta -1)|_{\frak{F}_r(n)} \; \; \mbox{mod} \; \widehat{T}(n+1).  \leqno{(3.1.3)}$$
Here we may note that $\Theta$ is multiplicative, $\Theta(f_1 f_2) = \Theta(f_1)\Theta(f_2)$ for $f_1, f_2 \in \frak{F}_r$, while $\Theta_n$ is additive, $\Theta_n(
[f_1 f_2]) = \Theta_n([f_1]+[f_2]) = \Theta_n([f_1]) + \Theta_n([f_2])$, where $[\,\cdot\,]$ stands for the class mod $\frak{F}_r(n+1)$. 

Let $S(H)$ be the symmetric algebra of $H$ over $\mathbb{Z}_l$ and let $q : T(H) \rightarrow S(H)$ be the natural map. We let $S^m(H) := q(H^{\otimes m})$ and $u_i := q(X_i)$ for $1\leq i \leq r$ so that $S(H)$ is the graded algebra $\bigoplus_{m\geq 0} S^m(H)$ which is noting but the commutative polynomial algebra over $\mathbb{Z}_l$ of variables $u_1, \dots , u_r$:
$$ S^m(H) = \bigoplus_{m\geq 0} S^m(H) = \mathbb{Z}_l[u_1,\dots, u_r].$$
\\
{\bf 3.2. The pro-$l$ Johnson map.}  This subsection concerns the pro-$l$ Johnson map associated to the Ihara representation, which is a pro-$l$ analogue of the Johnson map
 introduced by Kawazumi ([Ka]). Overall, we follow Kazazumi's arguments in [ibid] in our pro-$l$ setting.

For $\varphi \in {\rm Aut}^{\rm fil}(\widehat{T}(H))$, we denote by $[\varphi]$ the induced $\mathbb{Z}_l$-endomorphism of $H = \widehat{T}(1)/\widehat{T}(2) = \mathbb{Z}_l^{\oplus r}$.  \\
\\
{\bf Lemma 3.2.1.} {\it A $\mathbb{Z}_l$-algebra endomorphism $\varphi$ of $\widehat{T}(H)$ is a  filtration-preserving automorphism of $\widehat{T}(H)$, $\varphi \in {\rm Aut}^{\rm fil}(\widehat{T}(H))$, if and only if the following conditions are satisfied}:
\\
(1)  {\it $\varphi(\widehat{T}(n)) \subset \widehat{T}(n)$ for all $n \geq 0$.}\\
(2)  {\it the induced $\mathbb{Z}_l$-endomorphism $[\varphi]$ on $\widehat{T}(1)/\widehat{T}(2) = H$  is an isomorphism.}\\
\\
{\it Proof.}  Suppose $\varphi \in {\rm Aut}^{\rm fil}(\widehat{T}(H))$. Since $\varphi$ is filtration-preserving, the condition (1) holds. To show the condition (2), consider the following commutative diagram for vector spaces over $\mathbb{Z}_l$ with exact rows:
$$ \begin{array}{ccccccccc}
0 & \longrightarrow & \widehat{T}(2) & \longrightarrow & \widehat{T}(1) & \longrightarrow & H & \longrightarrow & 0 \\
  &                          & \;\;\;\;\;\;\; \downarrow {\footnotesize \varphi|_{\widehat{T}(2)}} & & \;\;\;\;\;\;\; \downarrow  {\footnotesize \varphi|_{\widehat{T}(1)}} &  & \;\;\;\;\; \downarrow {\footnotesize [\varphi]} & & \\
  0 & \longrightarrow & \widehat{T}(2) & \longrightarrow & \widehat{T}(1) & \longrightarrow & H & \longrightarrow & 0. 
  \end{array}
  $$
  Since $\varphi(\widehat{T}(n)) = \widehat{T}(n)$ for all $n \geq 0$, we have ${\rm Coker}(\varphi|_{\widehat{T}(i)}) = 0$ for $i=1,2$, in particular. Since $\varphi$ is an automorphism, we have ${\rm Ker}(\varphi) = 0$, in particular, 
  ${\rm Ker}(\varphi|_{\widehat{T}(i)}) = 0$ for $i=1,2$. By snake lemma applied to the above diagram, we obtain ${\rm Ker}([\varphi]) = 0$ and ${\rm Coker}([\varphi]) = 0$, hence the condition (2).

Suppose  that  a $\mathbb{Z}_l$-algebra endomorphism $\varphi$ of $\widehat{T}$ satisfies the conditions (1) and (2). Let $z = (z_m)$ be any element of $\widehat{T}$ with $z_m \in H^{\otimes m}$ for $m \geq 0$. To show that $\varphi$ is an automorphism, we have only to prove that there exists uniquely $y = (y_m) \in \widehat{T}$ such that
$$ z = \varphi(y). \leqno{(3.2.1.1)}$$
 Note by the condition (1) and (2) that $\varphi$ induces a $\mathbb{Z}_l$-linear automorphism of $\widehat{T}(m)/\widehat{T}(m+1) = H^{\otimes m}$, which is nothing but $[\varphi]^{\otimes m}$. Then, writing $\varphi(y_i)_j$ for the component
 of $\varphi(y_i)$ in $H^{\otimes j}$ for $i < j$, the equation (3.2.1.1) is equivalent to the following system of equations:
 $$\left\{ 
\begin{array}{l}
  z_0 = \varphi(y_0) = y_0,\\
 z_1 = [\varphi](y_1),\\
 z_2 = [\varphi]^{\otimes 2}(y_2) + \varphi(y_1)_2,\\
 \cdots \\
 z_m = [\varphi]^{\otimes m}(y_m) + \varphi(y_1)_m + \cdots + \varphi(y_{m-1})_m ,\\
 \cdots 
 \end{array} 
\right.
  \leqno{(3.2.1.2)}
 $$
Since $[\varphi]^{\otimes m}$ is an automorphism, we can find the unique solution $y = (y_m)$ of  (3.2.1.2)  from the lower degree. Therefore $\varphi$ is an $\mathbb{Z}_l$-algebra automorphism. Furthermore, we can see easily that if $z_0=\cdots = z_{n-1}=0$, then $y_0=\cdots = y_{n-1}=0$ for $n \geq 1$. This means that $\varphi^{-1}(\widehat{T}(n)) \subset \widehat{T}(n)$ and so $\varphi$ is filtration-preserving. $\;\; \Box$
\\

By Lemma 3.2.1, each $\varphi \in {\rm Aut}^{\rm fil}(\widehat{T}(H))$ induces a $\mathbb{Z}_l$-linear automorphism $[\varphi]$ of $H = \widehat{T}(1)/\widehat{T}(2)$ and so we have a group homomorphism 
$$ [\;\;] \; : \; {\rm Aut}^{\rm fil}(\widehat{T}(H)) \longrightarrow {\rm GL}(H),$$
where ${\rm GL}(H)$ denotes the group of $\mathbb{Z}_l$-linear automorphisms of $H$. We then define the {\it induced automorphism group} of $\widehat{T}$ by 
$$ \begin{array}{ll} {\rm IA}(\widehat{T}(H)) & := {\rm Ker}([\;\;])\\
 & = \{ \varphi \in {\rm Aut}(\widehat{T}(H))\; | \; \varphi(h) \equiv h \; \mbox{mod}\; \widehat{T}(2) \; \mbox{for any}\; h \in H \}.
\end{array}$$
We note that there is a natural splitting $s : {\rm GL}(H) \rightarrow {\rm Aut}^{\rm fil}(\widehat{T}(H))$ of $[\;\; ],$ which is defined by 
$$  s(P)((z_n)) := (P^{\otimes n}(z_n)) \;\; \mbox{for} \; P \in {\rm GL}(H).$$
In the following, we also regard $[P] \in {\rm GL}(H)$ as an element of ${\rm Aut}^{\rm fil}(\widehat{T})$ through the splitting $s$. 
Thus we have the following\\
\\
{\bf Lemma 3.2.2.}  {\it We have a semi-direct decomposition}
$$ {\rm Aut}^{\rm fil}(\widehat{T}(H)) \; = \; {\rm IA}(\widehat{T}(H)) \rtimes {\rm GL}(H); \;\; \varphi  =  (\varphi \circ [\varphi]^{-1}, [\varphi]). $$
\vskip 5pt

Let $\varphi \in {\rm IA}(\widehat{T}(H))$. Then we have $\varphi(h) - h \in \widehat{T}(2) \; \; \mbox{for any} \; h \in H$,  and so we have a map
$$ E \; : \; {\rm IA}(\widehat{T}) \longrightarrow {\rm Hom}_{\mathbb{Z}_l}(H,\widehat{T}(2)); \; \varphi \mapsto \varphi|_H - {\rm id}_H, \leqno{(3.2.3)}$$
where ${\rm Hom}_{\mathbb{Z}_l}(H,\widehat{T}(2))$ denotes the $\mathbb{Z}_l$-module of $\mathbb{Z}_l$-homomorphisms $H \rightarrow \widehat{T}(2)$. The following Proposition will play a key role in our discussion.\\
\\
{\bf Proposition 3.2.4.}  {\it The map $E$ is bijective.}\\
\\
{\it Proof.}  Injectivity: Suppose $E(\varphi) = E(\varphi')$ for $\varphi, \varphi' \in {\rm IA}(\widehat{T}(H))$. Then we have $\varphi|_H = \varphi'|_H$. Since an $\mathbb{Z}_l$-algebra endomorphism of $\widehat{T}(H)$ is determined by its restriction on $H$, we have $\varphi = \varphi'$. \\
Surjectivity: Take any $\phi \in {\rm Hom}_{\mathbb{Z}_l}(H,\widehat{T}(2))$.  We can extend $\phi + {\rm id}_H : H \rightarrow \widehat{T}(2)$ uniquely to a $\mathbb{Z}_l$-algebra endomorphism $\varphi$ of $\widehat{T}(H)$. Then we have obviously $\varphi(\widehat{T}(n)) \subset \widehat{T}(n)$ for all $n \geq 0$. Since $\widehat{T}(1)/\widehat{T}(2) = H$ and we see that 
$$ [\varphi](h \, {\rm mod} \, \widehat{T}(2)) = \varphi(h) \, {\rm mod}\, \widehat{T}(2) = h + \phi(h) \, {\rm mod} \, \widehat{T}(2) = h \, {\rm mod}\, \widehat{T}(2),$$
we have $[\varphi] = {\rm id}_{H}$. By Lemma 3.2.1, we have $\varphi \in {\rm IA}(\widehat{T})$ and $E(\varphi) = \phi$.  $\;\; \Box$\\
\\
By  Lemma 3.2.2 and Proposition 3.2.4, we have the following\\
\\
{\bf Corollary 3.2.5.}  {\it We have a bijection}
$$ \hat{E} \; : \; {\rm Aut}^{\rm fil}(\widehat{T}(H)) \simeq {\rm Hom}_{\mathbb{Z}_l}(H,\widehat{T}(2)) \times {\rm GL}(H)$$
 defined by $\hat{E}(\varphi) = (E(\varphi \circ [\varphi]^{-1}), [\varphi]).$
\\

Now, let ${\rm Ih}_S : {\rm Gal}_{k} \rightarrow P(\frak{F}_r)$ be the Ihara representation associated to $S$ in (1.2.3). 
We recall that the correspondence $\varphi \mapsto \varphi^* := \Theta \circ \varphi \circ \Theta^{-1}$ in (2.1.8) gives the injective homomorphism ${\rm Aut}(\frak{F}_r) \rightarrow {\rm Aut}^{\rm fil}(\widehat{T}(H))$ and hence the inclusion $P(\frak{F}_r) \hookrightarrow {\rm Aut}^{\rm fil}(\widehat{T}(H))$ which satisfies $[\varphi] = [\varphi^*]$ in ${\rm GL}(H)$. Composing ${\rm Ih}_S$ with this inclusion, we have  the homomorphism $\hat{\eta}_S : {\rm Gal}_{k} \rightarrow {\rm Aut}^{\rm fil}(\widehat{T}(H))$ defined by
$$ \hat{\eta}_S(g) := {\rm Ih}_S(g)^* = \Theta \circ {\rm Ih}_S(g) \circ \Theta^{-1}. $$
We then define the map $ \eta_S : {\rm Gal}_{k} \rightarrow {\rm IA}(\widehat{T}(H))$ by composing $\hat{\eta}_S$ with the projection on  ${\rm IA}(\widehat{T}(H))$:
$$ \eta_S(g) := \hat{\eta}_S(g) \circ [{\rm Ih}_S(g)]^{-1} = {\rm Ih}_S(g)^* \circ [{\rm Ih}_S(g)]^{-1} = \Theta \circ {\rm Ih}_S(g) \circ \Theta^{-1} \circ [{\rm Ih}_S(g)]^{-1}.    
 \leqno{(3.2.6)}$$
Thus,  we have $ \hat{\eta}_S(g) = (\eta_S(g), [{\rm Ih}_S(g)])$ for $g \in {\rm Gal}_{k}$ under the semi-direct decomposition ${\rm Aut}^{\rm fil}(\widehat{T}(H)) = {\rm IA}(\widehat{T}(H)) \rtimes {\rm GL}(H)$ of Lemma 3.2.2.\\

Now, we define the {\it pro-$l$ Johnson map}
$$ \tau_S \; : \; {\rm Gal}_{k} \longrightarrow {\rm Hom}_{\mathbb{Z}_l}(H,\widehat{T}(2)) $$
by the composing $\eta_S$ with $E$ in (3.2.3), and define the {\it extended  pro-$l$ Johnson map}
$$ \hat{\tau}_S \; : \; {\rm Gal}_{k} \longrightarrow {\rm Hom}_{\mathbb{Z}_l}(H,\widehat{T}(2)) \rtimes {\rm GL}(H) $$
by  composing $\hat{\eta}_S$ with $\hat{E}$ of Corollary 3.2.5. 
So we have, for $g \in {\rm Gal}_{k}$,
$$ \begin{array}{l}
\begin{array}{ll} 
\tau_S(g) & := E(\eta_S(g)) = \eta_S(g)|_H - {\rm id}_H \\
                & = {\rm Ih}_S(g)^* \circ [{\rm Ih}_S(g)]^{-1}|_H - {\rm id}|_H =  \Theta \circ {\rm Ih}_S(g) \circ \Theta^{-1}\circ [{\rm Ih}_S(g)]^{-1}|_H - {\rm id}|_H, 
                 \end{array}\\
\; \hat{\tau}_S(g) := (\tau_S(g), [{\rm Ih}_S(g)]).
\end{array}
\leqno{(3.2.7)}
$$
For $m \geq 1$, let ${\rm Hom}_{\mathbb{Z}_l}(H, H^{\otimes (m+1)})$ denote the $\mathbb{Z}_l$-module of $\mathbb{Z}_l$-homomorphisms $H \rightarrow H^{\otimes (m+1)}$, and we define the {\it $m$-th pro-$l$ Johnson map} 
$$ \tau_S^{(m)} \; : \; {\rm Gal}_{k} \longrightarrow {\rm Hom}_{\mathbb{Z}_l}(H, H^{\otimes (m+1)}) $$
by the $m$-th component of $\tau_S$:
$$ \tau_S(g) :=  \sum_{m\geq 1}\tau_{S}^{(m)}(g) \;\;\;\;  (g \in {\rm Gal}_{k}).  \leqno{(3.2.8)}$$

We note that the pro-$l$ Johnson map $\tau_S$ is no longer a homomorphism. In fact, we have the following\\
\\
{\bf Proposition 3.2.9.} {\it For $g_1, g_2 \in {\rm Gal}_k$, we have}
$$ \eta_S(g_1 g_2) = \eta_S(g_1)\circ [{\rm Ih}_S(g_1)]\circ \eta_S(g_2)\circ [{\rm Ih}_S(g_1)]^{-1}.$$
\\
{\it Proof.} By (3.2.6), we have
$$ \begin{array}{ll}
\eta_S(g_1 g_2) & = {\rm Ih}_S(g_1g_2)^* \circ [{\rm Ih}_S(g_1g_2)]^{-1}\\
                               & = \Theta \circ ({\rm Ih}_S(g_1 g_2) \circ \Theta^{-1} \circ [{\rm Ih}_S(g_1 g_2)]^{-1}\\
                              & = \Theta \circ {\rm Ih}_S(g_1) \circ {\rm Ih}_S(g_2) \circ \Theta^{-1} \circ [{\rm Ih}_S(g_2)]^{-1} \circ [{\rm Ih}_S(g_1)]^{-1}\\
                              & = \Theta \circ {\rm Ih}_S(g_1) \circ \Theta^{-1} \circ [{\rm Ih}_S(g_1)]^{-1} \circ [{\rm Ih}_S(g_1)] \circ \Theta \circ {\rm Ih}_S(g_2) \circ \Theta^{-1} \\
                              &    \;\;\;\;\; \circ [{\rm Ih}_S(g_2)]^{-1} \circ [{\rm Ih}_S(g_1)]^{-1}\\
                              & = \eta_S(g_1) \circ [{\rm Ih}_S(g_1)] \circ \eta_S(g_2) \circ [{\rm Ih}_S(g_1)]^{-1}.  \;\;\;\;\;\;\;\;\;\; \Box
\end{array}
$$
\\
Proposition 3.2.9 yields  coboundary relations among $\tau_S^{(m)}$. Here we give the formulas only for $\tau_S^{(1)}$ and $\tau_S^{(2)}$.\\
\\
{\bf Proposition 3.2.10.} {\em For $g_1, g_2 \in {\rm Gal}_k$, we have}
$$ \begin{array}{ll}
\tau_S^{(1)}(g_1 g_2) &= \tau_S^{(1)}(g_1) + [{\rm Ih}_S(g_1)]^{\otimes 2} \circ \tau_S^{(1)}(g_2) \circ [{\rm Ih}_S(g_1)]^{-1},\\
\tau_S^{(2)}(g_1 g_2) & = \tau_S^{(2)}(g_1) + (\tau_S^{(1)}(g_1)\otimes {\rm id}_H + {\rm id}_H \otimes \tau_S^{(1)}(g_1)) \circ [{\rm Ih}_S(g_1)]^{\otimes 2}\\
                                   & \;\;\;\;\; \circ \tau_S^{(1)}(g_2) \circ [{\rm Ih}_S(g_1)]^{-1} + [{\rm Ih}_S(g_1)]^{\otimes 3} \circ \tau_S^{(2)}(g_2) \circ [{\rm Ih}_S(g_1)]^{-1}.
\end{array}
$$
{\em Proof.} By definition (3.2.8), we have
$$ \tau_S(g_1 g_2) = \sum_{m \geq 1} \tau_S^{(m)}(g_1 g_2). \leqno{(3.2.10.1)}$$
On the other hand, by Proposition 3.2.9 and (3.2.7), we have, for $h \in H$,
$$ \begin{array}{ll}
\tau_S(g_1 g_2)(h) & = -h + \eta_S(g_1 g_2)(h)\\
           & = -h + (\eta_S(g_1)\circ [{\rm Ih}_S(g_1)] \circ \eta_S(g_2) \circ [{\rm Ih}_S(g_1)]^{-1})(h)\\
           & = -h + ( \eta_S(g_1)\circ [{\rm Ih}_S(g_1)] \circ ({\rm id}_H + \tau_S(g_2)))([{\rm Ih}_S(g_1)]^{-1}(h))\\
           & = \displaystyle{-h + (\eta_S(g_1) \circ [{\rm Ih}_S(g_1)]) \left( [{\rm Ih}_S(g_1)]^{-1}(h) + \sum_{m \geq 1} (\tau_S^{(m)}(g_2) \circ [{\rm Ih}_S(g_1)]^{-1})(h) \right)}\\
           & = \displaystyle{-h + \eta_S(g_1)\left( h + \sum_{m \geq 1} ([{\rm Ih}_S(g_1)]^{\otimes m+1} \circ \tau_S^{(m)}(g_2) \circ [{\rm Ih}_S(g_1)]^{-1})(h) \right) }\\
           & = -h + \eta_S(g_1)(h) \\
           & \;\;\;\;\; \;\;\; + \eta_S(g_1)( ([{\rm Ih}_S(g_1)]^{\otimes 2} \circ \tau_S^{(1)}(g_2) \circ [{\rm Ih}_S(g_1)]^{-1})(h))\\
           & \;\;\;\;\;  \;\;\; + \eta_S(g_1)( ([{\rm Ih}_S(g_1)]^{\otimes 3} \circ \tau_S^{(2)}(g_2) \circ [{\rm Ih}_S(g_1)]^{-1})(h)) \;\; {\rm mod}\;  \widehat{T}(4).
 \end{array}
 $$
We note that
$$ \eta_S(g)|_{H^{\otimes m}} = ({\rm id}_H + \tau_S(g))^{\otimes m} : H^{\otimes m} \longrightarrow H \times \widehat{T}(2m)$$
for any $g \in {\rm Gal}_{k}$ and so we have the following congruences mod $\widehat{T}(4)$:
$$ \begin{array}{l}
\eta_S(g_1)(h) \equiv h + \tau_S^{(1)}(g_1)(h) + \tau_S^{(2)}(g_1)(h),\\
\eta_S(g_1)(([{\rm Ih}_S(g_1)]^{\otimes 2} \circ \tau_S^{(1)}(g_2) \circ [{\rm Ih}_S(g_1)]^{-1})(h)) \\
 \;\;\;\;\;\;\;\; \;\; \equiv  ([{\rm Ih}_S(g_1)]^{\otimes 2} \circ \tau_S^{(1)}(g_2) \circ [{\rm Ih}_S(g_1)]^{-1})(h)\\
 \;\;\;\;\;\;\;\;\;\;\;\;                     +  ( (\tau_S^{(1)}(g_1)\otimes {\rm id}_H + {\rm id}_H \otimes \tau_S^{(1)}(g_1)) \circ [{\rm Ih}_S(g_1)]^{\otimes 2} \circ \tau_S^{(1)}(g_2) \circ [{\rm Ih}_S(g_1)]^{-1})(h),\\
 \eta_S(g_1)( ([{\rm Ih}_S(g_1)]^{\otimes 3} \circ \tau_S^{(1)}(g_2) \circ [{\rm Ih}_S(g_1)]^{-1})(h))  \equiv    ([{\rm Ih}_S(g_1)]^{\otimes 3} \circ \tau_S^{(2)}(g_2) \circ [{\rm Ih}_S(g_1)]^{-1})(h).
\end{array}
$$
Therefore we have
$$ \begin{array}{l} \tau_S(g_1 g_2)(h)   \\
\equiv  \tau_S^{(1)}(g_1)(h)  +   \tau_S^{(2)}(g_1)(h)\\
\;\; +  ([{\rm Ih}_S(g_1)]^{\otimes 2} \circ \tau_S^{(1)}(g_2) \circ [{\rm Ih}_S(g_1)]^{-1})(h) \\
 \;\; +   ( (\tau_S^{(1)}(g_1)\otimes {\rm id}_H + {\rm id}_H \otimes \tau_S^{(1)}(g_1)) \circ [{\rm Ih}_S(g_1)]^{\otimes 2} \circ \tau_S^{(1)}(g_2) \circ [{\rm Ih}_S(g_1)]^{-1})(h)\\
\;\;  +     ([{\rm Ih}_S(g_1)]^{\otimes 3} \circ \tau_S^{(2)}(g_2) \circ [{\rm Ih}_S(g_1)]^{-1})(h) \;\;\;\;  \mbox{mod} \; \widehat{T}(4).
\end{array}
\leqno{(3.2.10.2)}
$$                                                                                           
Comparing (3.2.10.1) and (3.2.10.2), we obtain the assertions. $\;\;\Box$\\
\\
{\bf 3.3. Pro-$l$ Johnson homomorphisms.}  For $n \geq 0$, let $\pi_n : \frak{F}_r \rightarrow \frak{F}_r/\frak{F}_r(n+1)$ be the natural homomorphism. Since each $\frak{F}_r(n)$ is a characteristic subgroup of $\frak{F}_r$, $\pi_n$ induces  the natural homomorphism
${\pi_n}_* : P(\frak{F}_r) \hookrightarrow {\rm Aut}(\frak{F}_r) \rightarrow {\rm Aut}(\frak{F}_r/\frak{F}_r(n+1))$.  Let ${\rm Ih}_{S}^{(n)}$ denote the composite of ${\rm Ih}_S$ with ${\pi_n}_*$:
$$ {\rm Ih}_{S}^{(n)} : {\rm Gal}_{k} \longrightarrow  {\rm Aut}(\frak{F}_r/\frak{F}_r(n+1)). $$
In particular, ${\rm Ih}_S^{(1)}(g) = [{\rm Ih}_S(g)]$ for $g \in {\rm Gal}_{k}$. Let ${\rm Gal}_{k}^{\rm Joh}[n]$ denote the kernel of ${\rm Ih}_{S}^{(n)}$:
$$ \begin{array}{ll} {\rm Gal}_{k}^{\rm Joh}[n]  & := {\rm Ker}({\rm Ih}_S^{(n)})\\
                                                                                & = \{ g \in   {\rm Gal}_{k} \, | \, {\rm Ih}_S(g)(f) f^{-1} \in \frak{F}_r(n+1) \; \mbox{for all}\; f \in \frak{F}_r \}.
\end{array} \leqno{(3.3.1)}$$
We then have the descending series of closed normal subgroups of ${\rm Gal}_{k}$:
$$ {\rm Gal}_{k} = {\rm Gal}_{k}^{\rm Joh}[0] \supset {\rm Gal}_{k}^{\rm Joh}[1] \supset \cdots \supset {\rm Gal}_{k}^{\rm Joh}
[n] \supset \cdots$$
and we call it the {\it Johnson filtration} of ${\rm Gal}_{k}$ associated to the Ihara representation $\varphi_S$ (cf. [Aa], [J1], [J2]). We note by Theorem 1.2.6 (1) 
$${\rm Gal}_{k}^{\rm Joh}[1] = {\rm Ker}({\rm Ih}_S^{(1)} :  {\rm Gal}_{k} \rightarrow {\rm GL}(H)) = {\rm Gal}_{k(\zeta_{l^{\infty}})}. \leqno{(3.3.2)}$$
The relation with the Milnor filtration defined in (2.2.13) is given as follows.\\
\\
{\bf Proposition 3.3.3.} {\em The Johnson filtration coincides with the Milnor filtration, namely, for each $n \geq 0$, we have}
$$   {\rm Gal}_{k}^{\rm Joh}[n] = {\rm Gal}_{k}^{\rm Mil}[n].$$
{\em Proof.} We may assume $n \geq 1$ and hence $g \in {\rm Gal}_{k(\zeta_{l^{\infty}})}$. Then we have
$$ \begin{array}{ll}
 g \in {\rm Gal}_{k}^{\rm Joh}[n]  & \Leftrightarrow  {\rm Ih}_S(g)(x_i) x_i^{-1} \in \frak{F}_r(n+1) \; \mbox{for all}\; 1\leq i \leq r\\
                                                                    & \Leftrightarrow  y_i(g)x_i y_i(g)^{-1} x_i^{-1} \in \frak{F}_r(n+1) \; \mbox{for all}\; 1\leq i \leq r\\
                                                                    & \Leftrightarrow  y_i(g) \in \frak{F}_r(n) \; \mbox{for all}\; 1\leq i  \leq r\\
                                                                    & \Leftrightarrow  {\rm deg}(\Theta(y_i(g) - 1)) \geq n \; \mbox{for all}\; 1\leq i  \leq r\\
                                                                    & \Leftrightarrow   g \in {\rm Gal}_{k}^{\rm Mil}[n]   \;\; \;\; \Box
\end{array}
$$
Note that Proposition 3.3.3 yields Proposition 2.2.14. In the following, we simply write ${\rm Gal}_{k}[n]$ for the $n$-th term of the Johnson (or Milnor) filtration for $n \geq 0$ and we denote by $k[n]$ the Galois subextension of $k$ in $\overline{\mathbb{Q}}$ corresponding to ${\rm Gal}_{k}[n]$. By (3.3.2), we have $k[1] = k(\zeta_{l^{\infty}})$. \\
 \\
We give some basic properties of the Johnson filtration. The following Lemma 3.3.4, Proposition 3.3.5 and Theorem 3.3.6 (2) were shown by Ihara for the case $r=2$. See [Ih1; Proposition 7, page 59] and also [O1]. We give herewith concise proofs for the sake of readers.\\
\\
{\bf Lemma 3.3.4.} {\it For $g \in {\rm Gal}_{k}[m]$ $(m\geq 0)$ and $f \in \frak{F}_r(n)$ $(n \geq 1)$, we have}
$$ {\rm Ih}_S(g)(f)f^{-1} \in  \frak{F}_r(m+n).$$
{\it Proof.}  We fix $m \geq 0$ and  $g \in {\rm Gal}_{k}[m]$. We prove the assertion by induction on $n$. For $n =1$, the assertion ${\rm Ih}_S(g)(f)f^{-1} \in \frak{F}_r(m+1)$ is true by the definition (3.3.1). Assume that
$$ {\rm Ih}_S(g)(f)f^{-1} \in \frak{F}_r(m+i) \; \mbox{if} \; f \in \frak{F}_r(i) \; \mbox{and} \; 1\leq i \leq n. \leqno{(3.3.4.1)}$$ 
Let $[\frak{F}_r(n),\frak{F}_r]_{\rm abst}$ denote the abstract group generated by $[a,b]$ ($a \in \frak{F}_r(n), b \in \frak{F}_r$). Since ${\rm Ih}_S(g)$ is continuous and $
[\frak{F}_r(n),\frak{F}_r]_{\rm abst}$ is dense in $\frak{F}_r(n+1)$, it suffices to show that
$${\rm Ih}_S(g)(f)f^{-1} \in \frak{F}_r(m+n+1) \; \mbox{for} \; f \in  [\frak{F}_r(n),\frak{F}_r]_{\rm abst}. 
$$
For this, we have only to show 
$$ {\rm Ih}_S(g)([b,c])[b,c]^{-1} \in \frak{F}_r(m+n+1) \; \mbox{if} \; b \in \frak{F}_r(n), c \in \frak{F}_r.  \leqno{(3.3.4.2)}
$$
For simplicity, we shall use the notation: $[\varphi, x] := \psi(x)x^{-1}$ and $[x,\varphi] := x \varphi(x)^{-1}$ for $x \in \frak{F}_r$ and $\varphi \in {\rm Aut}(\frak{F}_r)$. 
By  the ``three subgroup lemma" and the induction hypothesis (3.3.4.1), we have
$$ \begin{array}{ll}
{\rm Ih}_S(g)([b,c])[b,c]^{-1} & = [{\rm Ih}_S(g), [b,c]] \\
                                  & \in [{\rm Ih}_S(g), [\frak{F}_r(n), \frak{F}_r]]\\
                                  & \subset [[{\rm Ih}_S(g), \frak{F}_r(n)], \frak{F}_r][[\frak{F}_r, {\rm Ih}_S(g)], \frak{F}_r(n)]\\
                                  & \subset [\frak{F}_r(m+n), \frak{F}_r][\frak{F}_r(m+1), \frak{F}_r(n)]\\
                                  & = \frak{F}_r(m+n+1)
\end{array}
$$ 
and our claim (3.3.4.2) follows. $\Box$\\
\\                               
 Lemma 3.3.4 yields the following \\
 \\
 {\bf Proposition 3.3.5.}  {\it For $m, n \geq 0$, we have}
 $$  [{\rm Gal}_{k}[m], {\rm Gal}_{k}[n]] \subset {\rm Gal}_{k}[m+n] \;\; \mbox{for}\, m, n \geq 0.  $$
 {\em In particular, the Johnson $($or Milnor$)$ filtration is a central series.}\\
 \\
 {\it Proof.}  Using the same notation as in the proof of (3.3.4.2) and  Lemma 3.3.4, we have
 $$ \begin{array}{l}
 [[{\rm Gal}_{k}[n], \frak{F}_r],{\rm Gal}_{k}[m]] \subset [\frak{F}_r[n+1], {\rm Gal}_{k}[m]] \subset \frak{F}_r[m+n+1],\\
 
 [[\frak{F}_r,{\rm Gal}_{k}[m]], {\rm Gal}_{k}[n]] \subset [\frak{F}_r(m+1), {\rm Gal}_{k}[n]] \subset \frak{F}_r(m+n+1).
 \end{array}
 $$
 By the three subgroup lemma, we have
 $$ \begin{array}{ll} [[{\rm Gal}_{k}[m], {\rm Gal}_{k}[n]], \frak{F}_r] & \subset [{\rm Gal}_{k}[n], \frak{F}_r],{\rm Gal}_{k}[m]] [[\frak{F}_r, {\rm Gal}_{k}[m]], {\rm  Gal}_{k}[n]] \\
& \subset \frak{F}_r(m+n+1),
\end{array} $$
which yields the assertion by the definition (3.3.1).  $\;\; \Box$  \\
 
 For $n \geq 0$, let
 $$ {\rm gr}_n({\rm Gal}_k) := {\rm Gal}_k[n]/{\rm Gal}_k[n+1],$$
 which is a $\mathbb{Z}_l$-module.  Then, by Proposition 3.3.5, the graded $\mathbb{Z}_l$-module
 $$ {\rm gr}({\rm Gal}_k) := \bigoplus_{n \geq 0} {\rm gr}_n({\rm Gal}_k)$$
 has the structure of a graded Lie algebra over $\mathbb{Z}_l$, where the Lie bracket is defined by
 the commutator: For $ a = g$ mod ${\rm Gal}_k[m+1]$, $b = h$ mod ${\rm Gal}_k[n+1]$ ($g \in {\rm Gal}_k[m], h \in {\rm Gal}_k[n]$), 
 $$ [a,b] := [g,h] \; \mbox{mod} \; {\rm Gal}_k[m+n+1].$$
 
Now, for $m\geq 1$,  we let  $\tau_S^{[m]}$ denote the restriction of the $m$-th $l$-adic Johnson map $\tau_S^{(m)}$ in (3.2.8) to ${\rm Gal}_k[m]$:
$$ \tau_S^{[m]} := \tau_S^{(m)}|_{{\rm Gal}_k[m]} : {\rm Gal}_k[m] \longrightarrow {\rm Hom}_{\mathbb{Z}_l}(H, H^{\otimes (m+1)}). $$
The following theorem asserts that $\tau_S^{[m]}$ describes the action of  ${\rm Gal}_k[m]$ on $\frak{F}_r/\frak{F}_r(m+2)$. \\
 \\ 
{\bf Theorem 3.3.6.} {\em Notations being as above, the following assertions hold.}\\
(1) {\em For $g \in {\rm Gal}_{k}[m]$ and $f \in \frak{F}_r$, we have}
$$ \tau_S^{[m]}(g)([f]) = \Theta_{m+1}({\rm Ih}_S(g)(f)f^{-1}),$$
{\em where $\Theta_{m+1} : {\rm gr}_{m+1}(\frak{F}_r) \hookrightarrow H^{\otimes (m+1)}$ s the degree $(m+1)$-part of the Magnus embedding in $(3.1.3)$. }\\
(2) {\em The map $\tau_S^{[m]}$ is a $\mathbb{Z}_l$-homomorphism and ${\rm Ker}(\tau_S^{[m]}) = {\rm Gal}_{k}[m+1]$. Hence $\tau_S^{[m]}$ induces the injective $\mathbb{Z}_l$-homomorphism ${\rm gr}_m({\rm Gal}_k) \hookrightarrow {\rm Hom}_{\mathbb{Z}_l}(H, H^{\otimes (m+1)})$. In particular, we have}
$$ {\rm gr}_m({\rm Gal}_k) \simeq \mathbb{Z}_l^{\oplus r_m}$$
for some integer $r_m \geq 0$.\\
\\
{\em Proof.}  (1)  We need to show that for $g \in {\rm Gal}_k[m]$, 
$$ \tau_S^{(m)}(g)(X_i) = \Theta_{m+1}({\rm Ih}_S(g)(x_i)x_i^{-1}) \;\; 1\leq i \leq r.  \leqno{(3.3.6.1)}$$
By (3.2.7) and $[{\rm Ih}_S(g)] = {\rm id}_H$, we have
$$ \begin{array}{ll}
\tau_S(g)(X_i) & = (\Theta \circ {\rm Ih}_S(g) \circ \Theta^{-1})(\Theta(x_i)-1) - (\Theta(x_i)-1)\\
 & = \Theta({\rm Ih}_S(g)(x_i)) - \Theta(x_i).
 \end{array}
 $$
 Therefore, by (3.2.8),  we have
 $$\tau_S^{(m)}(g)(X_i) = \mbox{ the component in} \; H^{\otimes (m+1)} \; \mbox{of}\; \Theta({\rm Ih}_S(g)(x_i)) - \Theta(x_i).  \leqno{(3.3.6.2)}$$
 On the other hand, since ${\rm Ih}_S(g)(x_i)x_i^{-1} \in \frak{F}_r(m+1)$, we have
 $$ \Theta({\rm Ih}_S(g)(x_i)x_i^{-1}) \equiv 1 + \Theta_{m+1}({\rm Ih}_S(g)(x_i) x_i^{-1})  \;\; {\rm mod}\; \widehat{T}(m+2).$$
 Multiplying the above equation by  $\Theta(x_i)$ from right, we have
 $$ \Theta({\rm Ih}_S(g)(x_i)) \equiv \Theta(x_i) + \Theta_{m+1}({\rm Ih}_S(g)(x_i) x_i^{-1}) \;\; {\rm mod}\; \widehat{T}(m+2). \leqno{(3.3.6.3)}$$
By (3.3.6.2) and (3.3.6.3), we obtain (3.3.6.1).\\
(2)  By (1), for $g, h \in {\rm Gal}_k[m]$ and $f \in \frak{F}_r$, we have
$$ \begin{array}{ll}
\tau_S^{[m]}(gh)([f]) & = \Theta_{m+1}({\rm Ih}_S(gh)(f)f^{-1})\\
                                & =  \Theta_{m+1}({\rm Ih}_S(g)({\rm Ih}_S(h)(f)) f^{-1})\\
                                & =   \Theta_{m+1}({\rm Ih}_S(g)({\rm Ih}_S(h)(f)f^{-1}) {\rm Ih}_S(g)(f) f^{-1}).\\
\end{array}
$$
Since ${\rm Ih}_S(h)(f)f^{-1} \in \frak{F}_r(m+1)$, we have ${\rm Ih}_S(g)({\rm Ih}_S(h)(f)f^{-1}) \equiv {\rm Ih}_S(h)(f)f^{-1}$ mod $\frak{F}_r(2m+1) (\subset \frak{F}_r(m+2))$ by Lemma 3.3.4, and hence
$$
\begin{array}{ll}
\tau_S^{[m]}(gh)([f]) & =  \Theta_{m+1}({\rm Ih}_S(g)(f)f^{-1}) + \Theta_{m+1}({\rm Ih}_S(h)(f) f^{-1}).\\
 & = (\tau_S^{[m]}(g) + \tau_S^{[m]}(h))([f])
\end{array}
$$
for any $f \in \frak{F}_r$. Since ${\rm Ih}_S$ is continuous, we see that $\tau_S^{[m]}$ is a $\mathbb{Z}_l$-homomorphism. By (1) and (3.3.1), ${\rm Ker}(\tau_S^{[m]}) = {\rm Gal}_{k}[m+1]$, and hence $\tau_S^{[m]}$ induces the injective $\mathbb{Z}_l$-homomorphism ${\rm gr}_m({\rm Gal}_k) \hookrightarrow {\rm Hom}_{\mathbb{Z}_l}(H, H^{\otimes (m+1)})$. Since ${\rm Hom}_{\mathbb{Z}_l}(H, H^{\otimes (m+1)})$ is a free $\mathbb{Z}_l$-module of finite rank, the last assertion follows. $\;\; \Box$\\
\\                                               
 By Theorem 3.3.6 (1), $\tau_S^{[m]}$ factors through ${\rm Hom}_{\mathbb{Z}_l}(H, {\rm gr}_{m+1}(\frak{F}_r))$
$$ \tau_S^{[m]} : {\rm Gal}_k[m] \longrightarrow {\rm Hom}_{\mathbb{Z}_l}(H, {\rm gr}_{m+1}(\frak{F}_r)); \; g \mapsto ([f] \mapsto {\rm Ih}_S(g)(f) f^{-1})$$
followed by the map ${\rm Hom}_{\mathbb{Z}_l}(H, {\rm gr}_{m+1}(\frak{F}_r)) \rightarrow {\rm Hom}_{\mathbb{Z}_l}(H, H^{\otimes (m+1)})$ induced by $\Theta_{m+1}$. We call $\tau_S^{[m]} :   {\rm Gal}_{k}[m] \longrightarrow {\rm Hom}_{\mathbb{Z}_l}(H, H^{\otimes (m+1)})$ ($m \geq 1$) or the induced injective  $\mathbb{Z}_l$-homomorphism  ${\rm gr}_m({\rm Gal}_k) \hookrightarrow {\rm Hom}_{\mathbb{Z}_l}(H, H^{\otimes (m+1)})$, denoted by the same $\tau_S^{[m]}$, the {\em $m$-th  pro-$l$ Johnson homomorphism}. \\
\\
 A relation between the $m$-th pro-$l$ Johnson homomorphisms and $l$-adic Milnor numbers in Section 2 is given as follows.\\
\\
 {\bf Theorem 3.3.7.} {\em For $g \in {\rm Gal}_k[m]$ $(m\geq 1)$, we have}
 $$ \tau_S^{[m]}(g)(X_i) = - \sum_{|J| = m+1} \mu(J)X_J,$$
 where for $J = (j_1\cdots j_{m+1})$, 
$$\mu(J) = \left\{ \begin{array}{l}
 \mu(g; j_2\cdots j_{m+1}j_1) - \delta_{j_1, j_{m+1}}\mu(g; J) \; \;( i = j_1),\\
\mu(g; j_2\cdots j_{m+1} j_1)\delta_{j_1,j_{m+1}} - \mu(g; J) \;\;  (i = j_{m+1}),\\
0 \;\;  \mbox{(otherwise)}.
\end{array}\right.
$$
{\em Proof.} By Theorem 3.3.6 (1), we have
$$ \begin{array}{ll}
\tau_S^{[m]}(g)(X_i) & = \Theta_{m+1}({\rm Ih}_S(g)(x_i)x_i^{-1})\\
          & = \Theta_{m+1}(y_i(g)x_i y_i(g)^{-1}x_i^{-1}) \\
          & =  - \Theta_{m+1}([x_i,y_i(g)])\\
          & =  - \displaystyle{ \sum_{|J| = m+1} \mu(J; [x_i,y_i(g)])X_J.}
          \end{array}
          \leqno{(3.3.7.1)}
          $$
 By the computation in the proof of Theorem 2.3.3, we have, for $|J| = (j_1\cdots j_{m+1})$, 
 $$\begin{array}{ll} \mu(J; [x_i,y_i(g)])  & = \mu(J;x_iy_i(g)) - \mu(J;y_i(g)x_i) \\
 &
= \left\{ \begin{array}{l}
 \mu(g; j_2\cdots j_{m+1}j_1) - \delta_{j_1, j_{m+1}}\mu(g; J)\; \; ( i = j_1),\\
\mu(g; j_2\cdots j_{m+1} j_1)\delta_{j_1,j_{m+1}} - \mu(g; J) \;\;  (i = j_{m+1}).\\
0 \;\; \mbox{(otherwise)}.
\end{array}\right.
\end{array}
\leqno{(3.3.7.2)}
$$
 By (3.3.7.1) and (3.3.7.2), the assertion follows. $\;\; \Box$\\
 \\
 {\bf Remark 3.3.8.} A correspondence between Johnson invariants and Milnor invariants was given by Habegger in a topological framework ([Ha]). Our treatment
 in this paper is group-theoretical and similar to that given in [Ko1], [Ko3; Chapter 1] for pure braids. \\
 
We compute the pro-$l$ Johnson homomorphisms on commutators.\\
\\
{\bf Proposition 3.3.9.} {\it For $g \in {\rm Gal}_k[m], h \in {\rm Gal}_k[n]$ $(m, n \geq 0)$ and $f \in \frak{F}_r$, we have}
$$ \begin{array}{ll} \tau_S^{[m+n]}([g, h])([f]) & = \Theta_{m+n+1}({\rm Ih}_S(g)({\rm Ih}_S(h)(f) f^{-1})({\rm Ih}_S(h)(f) f^{-1})^{-1} \\
 \;\;\;\;\;\; \;\; \;\; &  -  {\rm Ih}_S(h)({\rm Ih}_S(g)(f)f^{-1})({\rm Ih}_S(g)(f)f^{-1})^{-1}).
\end{array}
$$
{\it Proof.} For simplicity, we set $\psi := {\rm Ih}_S(g), \phi := {\rm Ih}_S(h)$. By a straightforward computation using $[g, h] \in {\rm Gal}_k[m+n]$ (Proposition 3.3.5) and $\psi(f) f^{-1} \in \frak{F}_r(m+1)$ (Lemma 3.3.4), we obtain
$$  \begin{array}{l}
[\psi,\phi](f) f^{-1}\\
= [\psi,\phi]((\phi(f)f^{-1})^{-1}) \cdot (\psi \phi \psi^{-1})((\psi(f)f^{-1})^{-1})\cdot \psi(\phi(f)f^{-1}) \cdot \psi(f)f^{-1}\\
\equiv  (\phi(f) f^{-1})^{-1} \cdot \phi((\psi(f)(f^{-1})^{-1})\cdot \psi(\phi(f)f^{-1}) \cdot \psi(f)f^{-1} \;\; {\rm mod} \; \frak{F}_r(m+n+2).
\end{array}
$$
Since $\psi(f) f^{-1} \in \frak{F}_r(m+1), \phi(f)f^{-1} \in \frak{F}_r(n+1)$ and $[\frak{F}_r(m+1), \frak{F}_r(n+1)]\subset \frak{F}(m+n+2)$, we have
$$  \begin{array}{l}
[\psi,\phi](f) f^{-1}\\
\equiv  (\phi(f)f^{-1})^{-1} \cdot \psi(\phi(f)f^{-1}) \cdot \phi((\psi(f) f^{-1})^{-1})\cdot \psi(f)f^{-1} \;\; {\rm mod} \;  \frak{F}_r(m+n+2).
\end{array}
$$ 
Since we easily see that 
$$\left\{ 
\begin{array}{l}
(\phi(f)f^{-1})^{-1} \psi(\phi(f)f^{-1}) \equiv \psi(\phi(f)f^{-1}) (\phi(f)f^{-1})^{-1} \; {\rm mod} \; \frak{F}_r(m+n+2),\\
\phi((\psi(f) f^{-1})^{-1})\cdot \psi(f)f^{-1} \equiv (\phi(\psi(f) f^{-1})\cdot (\psi(f)f^{-1})^{-1})^{-1} \; {\rm mod}\;  \frak{F}_r(m+n+2),
\end{array}
\right.
$$
we obtain the assertion. $\;\; \Box$
\\
\\
By Proposition 3.3.9,  the direct sum of Johnson homomorphisms $\tau_S^{[m]}$ over all $m \geq 1$ defines a graded Lie algebra homomorphism from ${\rm gr}({\rm Gal}_k)$ to the derivation algebra of ${\rm gr}(\frak{F}_r)$ as follows. Recall that a $\mathbb{Z}_l$-linear endomorphism $\delta$ of ${\rm gr}(\frak{F}_r)$ is called a {\it derivation} on ${\rm gr}(\frak{F}_r)$ if it satisfies 
$$\delta([x,y]) = [\delta(x),y]+[x,\delta(y)] \;\;\;\; (x, y \in {\rm gr}(\frak{F}_r)).$$
 Let ${\rm Der}({\rm gr}(\frak{F}_r))$ denote the associative $\mathbb{Z}_l$-algebra of all derivations on ${\rm gr}(\frak{F}_r)$ which has  a Lie algebra structure over $\mathbb{Z}_l$ with the Lie bracket defined by $[\delta, \delta'] := \delta \circ \delta' - \delta' \circ \delta$ for $\delta, \delta' \in {\rm Der}({\rm gr}(\frak{F}_r))$.  For $m \geq 0$, we define the subspace ${\rm Der}_m({\rm gr}(\frak{F}_r))$ of ${\rm Der}({\rm gr}(\frak{F}_r))$, the degree $m$ part,  by
$$ {\rm Der}_m({\rm gr}(\frak{F}_r)):= \{ \delta \in {\rm Der}({\rm gr}(\frak{F}_r)) \, | \, \delta({\rm gr}_n(\frak{F}_r)) \subset {\rm gr}_{m+n}(\frak{F}_r) \; {\rm for}\; n \geq 1 \}$$
so that ${\rm Der}({\rm gr}(\frak{F}_r))$ is a graded Lie algebra over $\mathbb{Z}_l$:
$$ {\rm Der}({\rm gr}(\frak{F}_r)) = \bigoplus_{m\geq 0} {\rm Der}_m({\rm gr}(\frak{F}_r)).$$
A derivation $\delta \in {\rm Der}_m({\rm gr}(\frak{F}_r))$ is called a {\em special derivation} if there are $Y_i \in {\rm gr}_m(\frak{F}_r)$ such that $\delta(X_i) = [Y_i, X_i]$ $(1 \leq i \leq r)$ and moreover if the condition $\sum_{i=1}^r [Y_i,X_i] = 0$ is satisfied, a special derivation is said to be {\em normalized} ([Ih4; $\S 2$]). It is easy to see that normalized special derivations form a graded Lie subalgebra 
$${\rm Der}^{\rm n.s}({\rm gr}(\frak{F}_r)) = \bigoplus_{m\geq 0} {\rm Der}_m^{\rm n.s}({\rm gr}(\frak{F}_r))$$
 of ${\rm Der}({\rm gr}(\frak{F}_r))$. Since a derivation on ${\rm gr}(\frak{F}_r)$ is determined by its restriction on $H = {\rm gr}_1(\frak{F}_r)$, we have a natural inclusion, for each $m \geq 1$, 
$$ {\rm Der}_m({\rm gr}(\frak{F}_r)) \subset  {\rm Hom}_{\mathbb{Z}_l}(H,{\rm gr}_{m+1}(\frak{F}_r)); \;\; \delta \mapsto \delta|_H. $$
Hence we have the inclusions
$$ {\rm Der}^{\rm n.s}_{+}({\rm gr}(\frak{F}_r)) \subset {\rm Der}_{+}({\rm gr}(\frak{F}_r)) \subset \bigoplus_{m \geq 1} {\rm Hom}_{\mathbb{Z}_p}(H, {\rm gr}_{m+1}(\frak{F}_r)),$$
where ${\rm Der}_{+}({\rm gr}(\frak{F}_r))$ (resp. ${\rm Der}^{\rm n.s}_{+}({\rm gr}(\frak{F}_r))$) is the Lie subalgebra of ${\rm Der}({\rm gr}(\frak{F}_r))$  (resp. ${\rm Der}^{\rm n.s}({\rm gr}(\frak{F}_r))$) consisting of positive degree part. Although we make use of the arithmetic pro-$l$ Johnson homomorphisms, the following proposition was essentially proved by Ihara in [Ih4; $\S 2$].\\
\\
{\bf Proposition 3.3.10.} {\it The direct sum of $\tau_S^{[m]}$ over $m\geq 1$ defines the Lie algebra homomorphism}
$$ {\rm gr}(\tau) := \bigoplus_{m\geq 1} \tau_S^{[m]} \, : \, {\rm gr}({\rm Gal}_k) \longrightarrow {\rm Der}^{\rm n.s}_{+}({\rm gr}(\frak{F}_r)).$$
{\it Proof.} (cf. [Da; Proposition 3.18]) By Theorem 3.3.6 (1), it suffices to show that for $g \in {\rm Gal}_k[m]$, the map $f \mapsto {\rm Ih}_S(g)(f)f^{-1}$ is indeed a special derivation on ${\rm gr}(\frak{F}_r)$.  This was shown in [Ih4; $\S 2$] for the case $r=2$. We give herewith a proof for the sake of readers. Let  $g \in {\rm Gal}_k[m]$ $(m \geq 1)$ and $s \in \frak{F}_r(i), h \in \frak{F}_r(j)$. By using the commutator formulas
$$ [ab,c] = a[b,c]a^{-1}\cdot[a,c], \;\; [a,bc] = [a,b]\cdot  b[a,c]b^{-1} \;\; (a,b,c \in G),$$
we have
$$ \begin{array}{l}
{\rm Ih}_S(g)([s,t]) [s,t]^{-1} \\
 = [{\rm Ih}_S(g)(s), {\rm Ih}_S(g)(t)][s,t]^{-1}\\
 = [ss^{-1}{\rm Ih}_S(g)(s), {\rm Ih}_S(g)(t)t^{-1}t][s,t]^{-1}\\
 =  s [s^{-1} {\rm Ih}_S(g)(s), {\rm Ih}_S(g)(t)t^{-1}]\cdot ({\rm Ih}_S(g)(t)t^{-1})[ s^{-1}{\rm Ih}_S(g)(s),t]({\rm Ih}_S(g)(t)t^{-1})^{-1} s^{-1}\\
  \;\; \;\;\; \cdot [s,{\rm Ih}_S(g)(t)t^{-1}] ({\rm Ih}_S(g)(t)t^{-1})[s,t]({\rm Ih}_S(g)(t)t^{-1})^{-1}[s,t]^{-1}\\
 = s [s^{-1}{\rm Ih}_S(g)(s), {\rm Ih}_S(g)(t)t^{-1}]\cdot ({\rm Ih}_S(g)(t)t^{-1})[ s^{-1}{\rm Ih}_S(g)(s),t]({\rm Ih}_S(g)(t)t^{-1})^{-1} s^{-1}\\
 \;\; \;\;\; \cdot [s,{\rm Ih}_S(g)(t)t^{-1}] [{\rm Ih}_S(g)(t)t^{-1},[s,t]]. 
\end{array}
$$
Since $s^{-1}{\rm Ih}_S(g)(s) \in \frak{F}_r(i+m), {\rm Ih}_S(g)(t)t^{-1} \in \frak{F}_r(j+m)$ by Lemma 3.3.4, we have
$$ [s^{-1}{\rm Ih}_S(g)(s), {\rm Ih}_S(g)(t)t^{-1}] \in \frak{F}_r(i+j+2m), \;  [{\rm Ih}_S(g)(t)t^{-1},[s,t]] \in \frak{F}_r(i+2j+m).$$
By these claims together, we obtain
$$ \begin{array}{l}
{\rm Ih}_S(g)([s,t]) [s,t]^{-1} \\
\equiv s {\rm Ih}_S(g)(t)t^{-1} [ s^{-1}{\rm Ih}_S(g)(s),t](s {\rm Ih}_S(g)(t)t^{-1})^{-1}
 [s,{\rm Ih}_S(g)(t)t^{-1}] \;\; {\rm mod}\; \frak{F}_r(i+j+m+1).
\end{array}$$
Noting $x [s^{-1}{\rm Ih}_S(g)(s),t] x^{-1} \equiv [ s^{-1}{\rm Ih}_S(g)(s),t] \; {\rm mod} \; \frak{F}_r(i+j+m+1)$ for $x \in \frak{F}_r$, we proved that $f \mapsto {\rm Ih}_S(g)(f)f^{-1}$ is a derivation. That it is special and normalized follows from ${\rm Ih}_S(g)(x_i) = y_i(g)x_i y_i(g)^{-1}$ ($1\leq i \leq r$) and ${\rm Ih}_S(g) (x_1\cdots x_r) = x_1\cdots x_r$ for $g \in {\rm Gal}_k[m] \; (m\geq 1)$.
$\;\; \Box$\\
 
Finally we introduce an analogue of the Morita trace map ([Mt1; 6]).  For each $m \geq 1$, we identify 
${\rm Hom}_{\mathbb{Z}_l}(H, H^{\otimes (m+1)})$ with $H^* \otimes_{\mathbb{Z}_l} H^{\otimes (m+1)}$, where $H^* := {\rm Hom}_{\mathbb{Z}_l}(H,\mathbb{Z}_l)$ is the dual $\mathbb{Z}_l$-module, and let
$$ c_{m+1} : {\rm Hom}_{\mathbb{Z}_l}(H, H^{\otimes (m+1)}) = H^* \otimes_{\mathbb{Z}_l} H^{\otimes (m+1)} \longrightarrow H^{\otimes m}$$
be the contraction at $(m+1)$-component defined by
$$ c_{m+1}(\phi \otimes h_1 \otimes \cdots \otimes h_{m+1}) := \phi(h_{m+1}) h_1 \otimes \cdots \otimes h_m  \leqno{(3.3.11)}$$
for $\phi \in H^*, h_i \in H$. We then define the {\em $m$-th pro-$l$ Morita trace map}
$$ {\rm Tr}^{[m]} : {\rm Hom}_{\mathbb{Z}_l}(H, H^{\otimes (m+1)}) \longrightarrow S^m(H)  \leqno{(3.3.12)}$$
by the composite map $q \circ c_{m+1}$.\\

\begin{center}
{\bf 4.  Pro-$l$ Magnus-Gassner cocycles}
\end{center}
\vspace{.2cm}
{\bf 4.1. Pro-$l$ Fox free derivation.} The {\em pro-$l$ Fox free derivative} $\frac{\partial}{\partial x_j} : \mathbb{Z}_l[[\frak{F}_r]] \rightarrow \mathbb{Z}_l[[\frak{F}_r]]$ $(1\leq j \leq r)$ is a continuous $\mathbb{Z}_l$-linear map satisfying the following property: For any $\alpha \in \mathbb{Z}_l[[\frak{F}_r]]$, 
$$ \alpha =  \epsilon_{\mathbb{Z}_l[[\frak{F}_r]]}(\alpha) + \sum_{j=1}^r \frac{\partial \alpha}{\partial x_j}(x_j-1). \leqno{(4.1.1)}$$
We note by (4.1.1) that $\displaystyle{\frac{\partial \alpha}{\partial x_j}} \in I_{\mathbb{Z}_l[[\frak{F}_r]]}^{n-1}$ if $\alpha - \epsilon_{\mathbb{Z}_l[[\frak{F}_r]]}(\alpha) \in I_{\mathbb{Z}_l[[\frak{F}_r]]}^n$ for $n\geq 1$.

Here are  some basic rules for the pro-$l$ free calculus:
\vspace{.05cm}\\
(i) $\displaystyle{\frac{\partial x_i}{\partial x_j} = \delta_{ij}.}$\\
(ii) $\displaystyle{ \frac{\partial \alpha \beta}{\partial x_j} = \frac{\partial \alpha}{\partial x_j}\epsilon_{\mathbb{Z}_l[[\frak{F}_r]]}(\beta) + \alpha \frac{\partial \beta}{\partial x_j} }$  $\;\; (\alpha, \beta \in \mathbb{Z}_l[[\frak{F}_r]]).$\\
(iii) $\displaystyle{ \frac{\partial f^{-1}}{\partial x_j} = - f^{-1} \frac{\partial f}{\partial x_j}}$   $\;\; (f \in \frak{F}_r).$\\
(iv) $\displaystyle{ \frac{\partial f^{\alpha}}{\partial x_j} = \beta \frac{\partial f}{\partial x_j}}$  $\;(f \in \frak{F}_r, \alpha \in \mathbb{Z}_l)$, where $\beta$ is any element of  $\mathbb{Z}_l[[\frak{F}_r]]$ \\ $\;\;\;\;$ such that $\beta(f-1) = f^{\alpha} - 1$ if exists.\\
(v) $\displaystyle{ \frac{\partial \varphi(\alpha)}{\partial \varphi(x_j)} = \varphi (\frac{\partial \alpha}{\partial x_j}) }$  $\;\;(\varphi \in {\rm Aut}(\frak{F}_r), \alpha \in \mathbb{Z}_l[[\frak{F}_r]]).$  (Note that $\varphi(x_1), \dots , \varphi(x_r)$ \\ $\;\;\;\;$ are free generators of $\frak{F}_r$.)\\
(vi) If $\frak{F}'$ is an open free subgroup of $\frak{F}_r$ with free generators $y_1, \cdots , y_s$, we \\
 $\;\;\;\;$ have the chain rule: $\displaystyle{ \frac{\partial \alpha}{\partial x_j} = \sum_{i=1}^s \frac{\partial \alpha}{\partial y_i} \frac{\partial y_i}{\partial x_j } }$  $\;\; (\alpha \in \mathbb{Z}_l[[\frak{F}']]).$ 
 \vspace{.05cm}
 
The higher derivatives are defined inductively and the $l$-adic Magnus coefficient $\mu(I;\alpha)$ of $\alpha \in \mathbb{Z}_l[[\frak{F}_r]]$ for $I = (i_1\cdots i_n)$  is expressed by
$$ \mu(I;\alpha) = \epsilon_{\mathbb{Z}_l[[\frak{F}_r]]}\left(\frac{ \partial^n \alpha}{\partial x_{i_1} \cdots \partial x_{i_n}}\right) $$
so that the pro-$l$ Magnus expansion (2.1.4) is written as
$$ \Theta(\alpha) = \epsilon_{\mathbb{Z}_l[[\frak{F}_r]]}(\alpha) + \sum_{1\leq i_1,\dots , i_n \leq r} \epsilon_{\mathbb{Z}_l[[\frak{F}_r]]}\left(\frac{ \partial^n \alpha}{\partial x_{i_1} \cdots \partial x_{i_n}}\right) X_{i_1}\cdots X_{i_n}.$$
 \\
{\bf 4.2. Pro-$l$ Magnus cocycles.} Let ${\rm Ih}_S : {\rm Gal}_k \rightarrow P(\frak{F}_r) \subset {\rm Aut}(\frak{F}_r)$ be the Ihara representation associated to $S$ in (1.2.3).
Let $\bar{\;\;}  : \mathbb{Z}_l[[\frak{F}_r]] \rightarrow \mathbb{Z}_l[[\frak{F}_r]]$ denote the anti-automorphism induced by the involution $\frak{F}_r \ni f \mapsto f^{-1} \in \frak{F}_r$. We define the {\em pro-$l$ Magnus cocycle}  ${\rm M}_S : {\rm Gal}_k \rightarrow {\rm M}(r;\mathbb{Z}_l[[\frak{F}_r]])$ associated to ${\rm Ih}_S$ by 
$$ {\rm M}_{S}(g) := \left( \overline{\frac{ \partial {\rm Ih}_S(g)(x_j)}{\partial x_i}} \right) \leqno{(4.2.1)} $$
for $g \in {\rm Gal}_{k}$. In fact, we have the following \\
\\
{\bf Lemma 4.2.2.} {\em The map ${\rm M}_{S}$ is a $1$-cocycle of ${\rm Gal}_k$ with coefficients in ${\rm GL}(r; \mathbb{Z}_l[[\frak{F}_r]])$ with respect to the action ${\rm Ih}_S$. To be precise, for $g, h \in {\rm Gal}_k$, we have}
$$ {\rm M}_S(gh) = {\rm M}_S(g) {\rm Ih}_S(g)({\rm M}_S(h)),$$
{\em where ${\rm Ih}_S(g)({\rm M}_S(h))$ is the matrix obtained by applying ${\rm Ih}_S(g)$ to each entry of ${\rm M}_S(h)$.}\\
\\
{\em Proof.} Let $y_j := {\rm Ih}_S(h)(x_j)$ for $1\leq j \leq r$. Then we have
$$ \frac{\partial {\rm Ih}_S(gh)(x_j)}{\partial x_i} = \frac{\partial {\rm Ih}_S(g)(y_j)}{\partial x_i}. \leqno{(4.2.2.1)}$$
Using the basic rules (v), (vi) of the pro-$l$ Fox derivatives, we have
$$  \begin{array}{ll} \displaystyle{ \frac{\partial {\rm Ih}_S(g)(y_j)}{\partial x_i} }& = \displaystyle{ \sum_{a=1}^r \frac{\partial {\rm Ih}_S(g)(y_j)}{\partial  {\rm Ih}_S(g)(x_a)} \frac{\partial {\rm Ih}_S(g)(x_a)}{\partial x_i} }\\
 & = \displaystyle{ \sum_{a=1}^r {\rm Ih}_S(g) \left( \frac{\partial y_j}{\partial x_a} \right) \frac{\partial {\rm Ih}_S(g)(x_a)}{\partial x_i}. }
 \end{array}
 \leqno{(4.2.2.2)}
 $$
 By (4.2.2.1) and (4.2.2.2), we have
 $$ \overline{\frac{\partial {\rm Ih}_S(gh)(x_j)}{\partial x_i} } = \sum_{a=1}^r  \overline{\frac{\partial {\rm Ih}_S(g)(x_a)}{\partial x_i}}
  \cdot \overline{{\rm Ih}_S(g)\left(\frac{\partial y_j}{\partial x_a}\right)}.
 $$
 Since ${\rm Ih}_S(g)$ and $\bar{\;\;}$ are commutative operators, we obtain  the desired equality of the matrices. Taking $h = g^{-1}$, we see that ${\rm M}_S(g) \in {\rm GL}(r; \mathbb{Z}_l[[\frak{F}_r]])$ for $g \in {\rm Gal}_k$.
$\;\; \Box$.\\
\\
For $m \geq 1$, we let ${\rm M}_S^{[m]}$ be the composite of ${\rm M}_S$ restricted to ${\rm Gal}_k[m]$ with the natural homomorphism ${\rm GL}(r;\mathbb{Z}_l[[\frak{F}_r]]) \rightarrow {\rm GL}(r;\mathbb{Z}_l[[\frak{F}_r]]/I_{\mathbb{Z}_l[[\frak{F}_r]]}^{m+1})$
$$ {\rm M}_S^{[m]} : {\rm Gal}_k[m] \longrightarrow {\rm GL}(r;\mathbb{Z}_l[[\frak{F}_r]]/I_{\mathbb{Z}_l[[\frak{F}_r]]}^{m+1}).$$
A relation between ${\rm M}_S^{[m]}$ and the $m$-th pro-$l$ Johnson homomorphism is given as follows. First, recall the identification
 $\Theta_{n} : {\rm gr}_{n}(\frak{F}_r) \simeq H^{\otimes n}$ by the degree $n$ part of the  Magnus isomorphism in (3.1.2). We then have a matrix representation of ${\rm Hom}_{\mathbb{Z}_l}(H, H^{\oplus (m+1)})$ for $m \geq 1$
$$ ||\;\;|| :  {\rm Hom}_{\mathbb{Z}_l}(H, H^{\oplus (m+1)}) \longrightarrow {\rm M}(r; {\rm gr}_{m}(\mathbb{Z}_l[[\frak{F}_r]]))$$
by associating to each element $\tau \in  {\rm Hom}_{\mathbb{Z}_l}(H, H^{\oplus (m+1)}) $ the matrix
$$ ||\tau|| := \left(  \frac{ \partial (\Theta_{m+1}^{-1}\circ \tau)(X_j)}{\partial x_i} \right) \in {\rm M}(r; {\rm gr}_m(\mathbb{Z}_l[[\frak{F}_r]])). \leqno{(4.2.3)}$$
\\
{\bf Proposition 4.2.4.} {\em For $g \in {\rm Gal}_k[m]$, we have}
$$ {\rm M}_S^{[m]}(g) = I + \overline{||\tau_S^{[m]}(g)||}.$$
{\em Proof.} By Theorem 3.3.6, we have
$$ (\Theta_{m+1}^{-1}\circ \tau_S^{[m]}(g))(X_j) = {\rm Ih}_S(g)(x_j) x_j^{-1}$$
and so
$$ \begin{array}{ll}
\displaystyle{ \frac{\partial (\Theta_{m+1}^{-1}\circ \tau_S^{[m]})(X_j)}{\partial x_i}} & = \displaystyle{ \frac{\partial {\rm Ih}_S(g)(x_j) x_j^{-1}}{\partial x_i} }\\
              & = \displaystyle{ \frac{\partial {\rm Ih}_S(g)(x_j)}{\partial x_i} - {\rm Ih}_S(g)(x_j) x_j^{-1} \delta_{ij}.}
              \end{array}
$$
Since ${\rm Ih}_S(g)(x_j)x_j^{-1} \in \frak{F}_r(m+1)$, we have ${\rm Ih}_S(g)(x_j) x_j^{-1} \delta_{ij} \equiv \delta_{ij}$ mod $I_{\mathbb{Z}_l[[\frak{F}_r]]}^{m+1}$
and hence the assertion is proved. $\;\; \Box$ \\
\\
In terms of $||\cdot||$, the $m$-th pro-$l$ Morita trace ${\rm Tr}^{[m]}(\tau)$ in (3.3.12) is, in fact, written as the trace of the matrix $||\tau||$ .\\
\\
{\bf Proposition 4.2.5.} {\em For $m \geq 1$ and $,\tau \in {\rm Hom}_{\mathbb{Z}_l}(H, H^{\otimes (m+1)})$,  we have}
$$ {\rm Tr}^{[m]}(\tau) = q_m ({\rm tr}(\Theta_m(||\tau||))),$$
{\em where $q_m : H^{\otimes m} \rightarrow S^m(H)$ is the natural map.}\\
\\
{\em Proof.} We identify ${\rm Hom}_{\mathbb{Z}_l}(H, H^{\otimes (m+1)})$ with $H^* \otimes H^{\otimes m}$. Let $\tau = \phi \otimes X_{i_1} \otimes \cdots \otimes X_{i_{m+1}}$ ($\phi \in H^*$). By (4.2.3), we have
$$ \displaystyle{ {\rm tr}(||\tau||) = \sum_{i=1}^r \frac{\partial (\Theta_{m+1}^{-1} \circ \tau)(X_i)}{\partial x_i} = \sum_{i=1}^r \phi(X_i) \frac{\partial \Theta_{m+1}^{-1}(X_{i_1} \otimes \cdots \otimes X_{i_{m+1}})}{\partial x_i}.} \leqno{(4.2.5.1)}$$
We note that any element $Y$ of $H^{\otimes (m+1)}$ can be written uniquely as 
$$ Y = Y_1 \otimes X_1 + \cdots + Y_r\otimes X_r,\;\; Y_i \in H^{\otimes m}$$
and then we have, by (4.1.1),
$$ \frac{\partial \Theta_{m+1}^{-1}(Y)}{\partial x_i} = \Theta_m^{-1}(Y_i).$$
Therefore we have
$$ \frac{\partial \Theta_{m+1}^{-1}(X_{i_1} \otimes \cdots \otimes X_{i_{m+1}})}{\partial x_i} = \delta_{i, i_{m+1}} X_{i_1}\otimes \cdots \otimes X_{i_m} $$
and hence, by (4.2.5.1),
$${\rm tr}(\Theta_m(||\tau||)) = \phi(X_{i_{m+1}})X_{i_1}\otimes \cdots \otimes X_{i_m},$$
where the right hand side is $c_{m+1}(\tau)$ by (3.3.11). By the definition (3.3.12), the assertion is proved. $\;\; \Box$\\

Now, for some application later on, we extend the construction of the pro-$l$ Magnus cocycle to a relative situation. Let $\frak{G}$ be a pro-$l$ group and  let $\psi : \frak{F}_r \rightarrow \frak{G}$ be a continuous surjective homomorphism. We also denote by $\psi$ the induced surjective homomorphism $\mathbb{Z}_l[[\frak{F}_r]] \rightarrow \mathbb{Z}_l[[\frak{G}]]$ of complete group algebras over $\mathbb{Z}_l$. Let $\frak{N} := {\rm Ker}(\psi)$ so that $\frak{F}_r/\frak{N} \simeq \frak{G}$. We assume that $\frak{N}$ is stable under the action of ${\rm Gal}_k$ through ${\rm Ih}_S$, namely  ${\rm Ih}_S(g)(\frak{N}) \subset \frak{N}$ for all $g \in {\rm Gal}_k$ (This is certainly satisfied if $\frak{N}$ is a characteristic subgroup of $\frak{F}_r$). Then we have a homomorphism ${\rm Ih}_{S,\psi} : {\rm Gal}_k \rightarrow  {\rm Aut}(\mathbb{Z}_l[[\frak{G}]])$ defined by
$$ {\rm Ih}_{S,\psi}(g)(\psi(\alpha))  := \psi({\rm Ih}_S(g)(\alpha))  \;\; (\alpha \in \mathbb{Z}_l[[\frak{F}_r]] ). \leqno{(4.2.6)} $$ 
Let ${\rm Gal}_k[\psi]$ be the subgroup of ${\rm Gal}_k$ defined by 
 $$ \begin{array}{ll} {\rm Gal}_{k}[\psi] & := {\rm Ker}({\rm Ih}_{S, \psi})\\
                    & = \{ g \in {\rm Gal}_k \; | \; \psi \circ {\rm Ih}_S(g) = \psi \}\\
                    \end{array}
       \leqno{(4.2.7)}             $$
and let $k[\psi]$ denote the subfield of $\overline{\mathbb{Q}}/k$  corresponding to ${\rm Gal}_k[\psi]$. Now we  define the {\em pro-$l$ Magnus cocycle} ${\rm M}_{S, \psi} : {\rm Gal}_{k} \rightarrow {\rm GL}(r; \mathbb{Z}_l[[\frak{G}]])$ associated to ${\rm Ih}_S$ and $\psi$ by
$$ {\rm M}_{S, \psi}(g) := \psi ({\rm M}_S(g)) \;\; \;\; (g \in {\rm Gal}_k),$$
where the right hand side is the matrix obtained by applying $\psi$ to each entry of ${\rm M}_S(g)$. For $m \geq 1$, let ${\rm M}_{S,\psi}^{[m]}$ be the composite of ${\rm M}_S^{[m]}$ with the natural homomorphism ${\rm GL}(r;\mathbb{Z}_l[[\frak{F}_r]]/I_{\mathbb{Z}_l[[\frak{F}_r]]}^{m+1}) \rightarrow {\rm GL}(r;\mathbb{Z}_l[[\frak{G}]]/I_{\mathbb{Z}_l[[\frak{G}]]}^{m+1})$ induced by $\psi$.  
Lemma 4.2.2 and Proposition 4.2.4 are extended to the following \\
\\
{\bf Proposition 4.2.8.} {\em Notations being as above, the following assertions hold.}\\
(1) {\em For  $g, h \in {\rm Gal}_k$, we have}
$$ {\rm M}_{S,\psi}(gh) = {\rm M}_{S,\psi}(g) {\rm Ih}_{S,\psi}(g)({\rm M}_{S,\psi}(h)).$$
(2) {\em For $g \in  {\rm Gal}_k$, we have} 
$$ {\rm M}_{S,\psi}^{[m]}(g) = I + \psi(\overline{||\tau_S^{[m]}(g)||}). $$
(3) {\em The restriction of ${\rm M}_{S, \psi}$ to ${\rm Gal}_{k}[\psi]$, denoted by the same  ${\rm M}_{S, \psi}$,}
$${\rm M}_{S, \psi} : {\rm Gal}_k[\psi] \longrightarrow {\rm GL}(r;\mathbb{Z}_l[[\frak{G}]/I_{\mathbb{Z}_l[[\frak{G}]]}^{m+1})$$
{\em is a homomorphism and factors through the Galois group ${\rm Gal}(\Omega_S/k[\psi])$, where $\Omega_S$ is the subfield of $\overline{\mathbb{Q}}$ corresponding to ${\rm Ker}({\rm Ih}_S)$ as in $(1.2.4)$.} We call it the {\em pro-$l$ Magnus representation} of ${\rm Gal}_k[\psi]$ associated to ${\rm Ih}_S$ and $\psi$.\\
\\
{\em Proof.} (1) The formula is obtained by applying $\psi$ to the both sides of the formula in Lemma 4.2.2. (2)  This is also obtained by applying $\psi$ to the matrices of the both sides of the formula in Proposition 4.2.4.  (3) Suppose $g, h \in {\rm Gal}_{k,\psi}$. Since $\psi \circ {\rm Ih}_S(g) = \psi$, we have ${\rm Ih}_{S,\psi}(g)({\rm M}_{S,\psi}(h)) = {\rm M}_{S,\psi}(h)$ and so ${\rm M}_{S,\psi}(gh) = {\rm M}_{S,\psi}(g){\rm M}_{S,\psi}(h)$. Since ${\rm M}_{S,\psi}(g) = I$ for $g \in {\rm Ker}({\rm Ih}_S)$, we have ${\rm Ker}({\rm M}_{S,\psi}) \supset {\rm Ker}({\rm Ih}_S)$ and hence ${\rm M}_{S,\psi}$ factors through ${\rm Gal}(\Omega_S/k[\psi])$. $\;\; \Box$\\

For $n \geq 0$, let $\pi_n : \frak{F}_r \rightarrow \frak{F}_r/\frak{F}_r(n+1)$ be the natural homomorphism. We consider the case that $\psi = \pi_n$ and so ${\rm Ih}_{S,\psi} = {\rm Ih}_{S}^{(n)}$. By (3.3.1) and Lemma 3.3.3, we have 
$$ \begin{array}{ll} {\rm Gal}_k[\pi_n] & = \{ g \in {\rm Gal}_k \; | \; \pi_n \circ {\rm Ih}_S(g) = \pi_n \}\\
                                                                              & = \{ g \in {\rm Gal}_k \; | \; {\rm Ih}_S(g)(f) \equiv f \; \mbox{mod} \; \frak{F}_r(n+1) \; \mbox{for all}\; f \in \frak{F}_r \}\\
                                                                              & = {\rm Gal}_k[n].
                                                                              \end{array}
$$
 Then we have a family of pro-$l$ Magnus cocycles 
$$ {\rm M}_{S,\pi_n} : {\rm Gal}_k \longrightarrow {\rm GL}(r; \mathbb{Z}_l[[\frak{F}_r/\frak{F}_r(n+1)]]), \leqno{(4.2.9)}$$
and the pro-$l$ Magnus representation
$${\rm M}_{S,\pi_n} : {\rm Gal}_k[n] \longrightarrow {\rm GL}(r; \mathbb{Z}_l[[\frak{F}_r/\frak{F}_r(n+1)]]) \leqno{(4.2.10)}$$
associated to ${\rm Ih}_S$ and $\pi_n$ for $n \geq 0$.  \\
\\
{\bf 4.3.  Pro-$l$ Gassner cocycles.}  This subsection concerns the pro-$l$ (reduced) Gassner cocycles as special cases of the Magnus cocycles. For the construction of the pro-$l$ reduced Gassner cocycles, we follow Oda's arguments [O2]. We also refer to [N; II] for Magnus-Gassner matrices.

The pro-$l$ Gassner cocycle is  defined by ${\rm M}_{S,\pi_1}$ in (4.2.9). 
To be precise, let $\Lambda_r := \mathbb{Z}_l[[u_1,\dots , u_r]]$ denote the algebra of commutative formal power series over $\mathbb{Z}_l$ of variables $u_1, \dots , u_r$, called the {\em Iwasawa algebra} of $r$ variables. The correspondence $x_i \; \mbox{mod}\; \frak{F}_r(2) \mapsto 1+u_i$ ($1\leq i \leq r$) gives the abelianized pro-$l$ Magnus isomorphism
$$ \theta : \mathbb{Z}_l[[\frak{F}_r/\frak{F}_r(2)]] \stackrel{\sim}{\longrightarrow} \Lambda_r.$$
We let $\pi := \pi_1$ and 
$$\chi_{\Lambda_r} := {\rm Ih}_{S, \theta \circ \pi} : {\rm Gal}_k \rightarrow {\rm Aut}(\Lambda_r), \leqno{(4.3.1)}$$ 
which is defined by (4.2.4) with $\psi = \theta \circ \pi$. In fact, by Lemma 2.2.1, $\chi_{\Lambda_r}$ is given by
$$ \chi_{\Lambda_r}(g)(u_i) = (\theta \circ \pi)({\rm Ih}_S(g)(x_i -1)) = (1+u_i)^{\chi_l(g)} -1 \;\;\;\; (1\leq i \leq r). \leqno{(4.3.2)}$$
Then the {\em pro-$l$ Gassner cocycle} of ${\rm Gal}_k$ associated to ${\rm Ih}_S$ 
$$ {\rm Gass}_S : {\rm Gal}_k \longrightarrow {\rm GL}(r; \Lambda_r)$$
is defined by
$$ {\rm Gass}_S(g) := \displaystyle{\left( (\theta \circ \pi)\left(\frac{\partial {\rm Ih}_S(g)(x_j)}{\partial x_i}\right)\right)} \;\;\;\; (g \in {\rm Gal}_k), \leqno{(4.3.3)}$$
where we note that we do not need to take the anti-automorphism $\bar{\;\;}$ in (4.3.3) to obtain the 1-cocycle relation
$$ {\rm Gass}_S(gh) = {\rm Gass}_S(g) \chi_{\Lambda_r}(g)({\rm Gass}_S(h)) \;\;\;\; (g, h \in {\rm Gal}_k),$$ 
since $\Lambda_r$ is commutative.  Here $\chi_{\Lambda_r}(g)({\rm Gass}_S(h))$ is the matrix obtained by applying $\chi_{\Lambda_r}(g)$ to each entry of ${\rm Gass}_S(h)$. We can express ${\rm Gass}_S(g)$  in terms of $l$-adic Milnor numbers as follows.\\
\\
{\bf Proposition 4.3.4.} {\em The $(i,j)$-entry of ${\rm Gass}_S(g)$ $(g \in {\rm Gal}_k)$ is expressed by}
$$ {\rm Gass}_S(g)_{ij} = \left\{ 
\begin{array}{ll}
\displaystyle{ \frac{\chi_{\Lambda_r}(g)(u_i)}{u_i}\left( 1 + \sum_{n\geq 1} \sum_{{\scriptstyle 1\leq i_1,\dots , i_n\leq r}\atop{\scriptstyle  i_n \neq i }} \mu(g; i_1\cdots i_n i) u_{i_1} \cdots u_{i_n}\right)} &  (i=j),\\
\displaystyle{ - \chi_{\Lambda_r}(g)(u_j) \left( \mu(g;ij) +  \sum_{n\geq 1} \sum_{1\leq i_1,\dots , i_n\leq r} \mu(g;i_1\cdots i_n i j) u_{i_1}\cdots u_{i_n}\right) } &  (i \neq j).\\
\end{array}\right.
$$
{\em Proof.}  By Lemma 2.2.1 and a straightforward computation, we have
$$ \begin{array}{ll} \displaystyle{ \frac{\partial {\rm Ih}_S(g)(x_j)}{\partial x_i} }  & = \displaystyle{ \frac{\partial y_j(g)x_j^{\chi_l(g)}y_j(g)^{-1}}{\partial x_i} }\\
       & = \displaystyle{y_j(g)\frac{x_j^{\chi_l(g)}-1}{x_j-1} \delta_{ij} + (1-y_j(g)x_j^{\chi_l(g)}y_j(g)^{-1})\frac{\partial y_j(g)}{\partial x_i}   }
\end{array}
$$
and hence, by (4.3.2), 
$$ \begin{array}{l} \displaystyle{ (\theta \circ \pi)\left( \frac{\partial {\rm Ih}_S(g)(x_j)}{\partial x_i} \right) } \\
\displaystyle{ = (\theta \circ \pi)(y_j(g)) \frac{(1+u_j)^{\chi_l(g)}-1}{u_j}\delta_{ij} + (1-(1+u_j)^{\chi_l(g)}) (\theta \circ \pi)\left( \frac{\partial y_j(g)}{\partial x_i}\right)}\\
\displaystyle{ =   \frac{\chi_{\Lambda_r}(g)(u_j)}{u_j}(\theta \circ   \pi)(y_j(g))\delta_{ij} - \chi_{\Lambda_r}(g)(u_j)   (\theta \circ \pi)\left( \frac{\partial y_j(g)}{\partial x_i}\right).}                                                                                                                                                      
\end{array}
 \leqno{(4.3.4.1)} $$
Here we have 
$$\displaystyle{(\theta \circ \pi)(y_j(g)) = 1 + \sum_{|I|\geq 1}\mu(g; Ij)u_I,} \leqno{(4.3.4.2)}$$
where we set $u_I := u_{i_1}\cdots u_{i_n}$ for $I = (i_1\cdots i_n)$, and (4.1.1) yields 
$$\displaystyle{ (\theta \circ \pi)\left(\frac{\partial y_j(g)}{\partial x_i}\right) = \sum_{|I|\geq 0} \mu(g; Iij)u_I. } \leqno{(4.3.4.3)}$$
By (4.3.3), (4.3.4.1), (4.3.4.2) and (4.3.4.3), we have
$$ \begin{array}{l}
\displaystyle{ {\rm Gass}_S(g)= (\theta \circ \pi)\left( \frac{\partial {\rm Ih}_S(g)(x_j) }{\partial x_i} \right) }\\
\displaystyle{ = \delta_{ij} \frac{\chi_{\Lambda_r}(g)(u_j)}{u_j}\left(1+ \sum_{|I|\geq 1}\mu(g;Ij)u_I\right) - \chi_{\Lambda_r}(g)(u_j) \sum_{|I|\geq 0} \mu(g;Iij)u_I. } \\

\end{array}
$$
By $\mu(g;ii) = 0$ and a simple observation, we obtain the assertion.  $\;\; \Box$\\
\\
By (4.2.10), when ${\rm Gass}_S$ is restricted to ${\rm Gal}_k[1]$, we have a representation
$$ {\rm Gass}_S : {\rm Gal}_k[1] \longrightarrow {\rm GL}_r(\Lambda_r), $$
which we call the {\em pro-$l$ Gassner representation} of ${\rm Gal}_k[1]$ associated to ${\rm Ih}_S$.  It factors through the Galois group ${\rm Gal}(\Omega_S/k[1])$
by Theorem 4.2.8 (3).\\

In the following, for simplicity, we let 
$$\frak{F}_r' := \frak{F}_r(2), \frak{F}_r'' := [\frak{F}_r',\frak{F}_r'], \; \mbox{and}\; \frak{L}_r := \frak{F}_r'/\frak{F}_r'' = H_1(\frak{F}_r',\mathbb{Z}_l).$$
 We consider $\frak{L}$ as a $\Lambda_r$-module by conjugation: For $ f  \in \frak{F}_r$ and $f' \in \frak{F}_r'$, we set
$$ [f].(f' \; \mbox{mod}\; \frak{F}_r'') := f f' f^{-1}\; \mbox{mod} \; \frak{F}_r''.$$
and extend it by the $\mathbb{Z}_l$-linearity and continuity. The structure of the $\Lambda_r$-module $\frak{L}_r$ can be described by means of the pro-$l$ Crowell exact sequence ([Ms2; Chapter 9]).   Attached to the surjective homomorphism $\pi : \frak{F}_r \longrightarrow \frak{F}_r/\frak{F}_r'$,
the {\em pro-$l$ Crowell exact sequence} reads as the exact sequence of $\Lambda_r$-modules:
$$ 0 \longrightarrow \frak{L}_r \stackrel{\nu_1}{\longrightarrow} \Lambda_r^{\oplus r} \stackrel{\nu_2}{\longrightarrow} I_{\Lambda_r} \longrightarrow 0,$$
where $I_{\Lambda_r}$ is the (augmentation) ideal of $\Lambda_r$ generated by $u_1, \dots , u_r$ and $\nu_1, \nu_2$ are $\Lambda_r$-homomorphisms defined by
$$\displaystyle{ \nu_1(f' \; \mbox{mod}\; \frak{F}_r'') := ( (\theta \circ \pi)\left( \frac{\partial f'}{\partial x_i} \right)) \;\; (f' \in \frak{F}_r'); \;\; \;\;
 \nu_2((\lambda_i)) := \sum_{i=1}^r \lambda_i u_i.} \leqno{(4.3.5)}
$$
({\em Convention}: An element $(\lambda_i)$ of $\Lambda_r^{\oplus r}$ is understood as a column vector.) 
Hence we have the isomorphism of $\Lambda_r$-modules induced by $\nu_1$, called the {\em Blanchfield-Lyndon isomorphism}: 
$$ \nu_1 : \frak{L}_r \stackrel{\sim}{\longrightarrow}\{ (\lambda_i) \in \Lambda_r^{\oplus r} \; | \;  \sum_{i=1}^r \lambda_i u_i = 0 \}. \leqno{(4.3.6)}$$

We define the action ${\rm Meta}_S$ of ${\rm Gal}_k$ on $\frak{L}_r$ through the Ihara representation ${\rm Ih}_S$: For $g \in {\rm Gal}_k$ and $f' \in \frak{F}_r'$,
$$ {\rm Meta}_S(g)(f' \; \mbox{mod} \; \frak{F}_r'') := {\rm Ih}_S(g)(f') \; \mbox{mod} \; \frak{F}_r''.$$
It is easy to see that ${\rm Meta}_S(g)$ is a $\chi_{\Lambda_r}$-linear automorphism of $\frak{L}_r$, namely, a $\mathbb{Z}_l$-linear automorphism and satisfies
$$ {\rm Meta}_S(g)(\lambda.(f' \; \mbox{mod}\; \frak{F}_r'')) = \chi_{\Lambda_r}(g)(\lambda).(f' \; \mbox{mod}\; \frak{F}_r'')$$
for $\lambda \in \Lambda_r$ and $f' \in \frak{F}_r'$. When ${\rm Meta}_S$ is restricted to ${\rm Gal}_k[1]$, we have the representation, which we call the {\em pro-$l$ meta-abelian representation} of ${\rm Gal}_k[1]$ associated to ${\rm Ih}_S$, 
$$ {\rm Meta}_S : {\rm Gal}_k[1] \longrightarrow {\rm GL}_{\Lambda_r}(\frak{L}_r),$$
where ${\rm GL}_{\Lambda_r}(\frak{L}_r)$ is the group of $\Lambda_r$-module automorphisms of $\frak{L}_r$. Regarding $\frak{L}_r$ as a $\Lambda_r$-submodule of $\Lambda_r^{\oplus r}$ by the isomorphism (4.3.6),  ${\rm Meta}_S$ and ${\rm Gass}_S$ has the following relation.\\
\\
{\bf Proposition 4.3.7.}  {\em For $g \in {\rm Gal}_k$ and $f' \in \frak{F}_r'$, we have}
$$ (\nu_1 \circ {\rm Meta}_S(g))(f' \; \mbox{mod}\; \frak{F}_r'') = {\rm Gass}_S(g) (\chi_{\Lambda_r}(g) \circ \nu_1)(f' \; \mbox{mod}\; \frak{F}_r'').$$
{\em When ${\rm Meta}_S$ and ${\rm Gass}_S|_{\frak{L}_r}$ are restricted to ${\rm Gal}_k[1]$, they are  equivalent representations over $\Lambda_r$.}
\\
\\
{\em Proof.}  The first assertion follows from the direct computation: By (4.3.1), (4.3.3) and (4.3.5), we have, for any $g \in {\rm Gal}_k$ and $f' \in \frak{F}_r'$, 
$$ \begin{array}{ll} 
(\nu_1\circ {\rm Meta}_S(g))(f' \; \mbox{mod}\; \frak{F}_r'') & = \nu_1( {\rm Ih}_S(g)(f') \; \mbox{mod}\; \frak{F}_r'')\\
 & = \displaystyle{ ( (\theta \circ \pi) \left( \frac{\partial {\rm Ih}_S(g)(f')}{\partial x_i} \right) ) }\\
  & =  \displaystyle{ ( (\theta \circ \pi)\left( \sum_{a=1}^r \frac{\partial {\rm Ih}_S(g)(f')}{\partial {\rm Ih}_S(g)(x_a)} \frac{\partial {\rm Ih}_S(g)(x_a)}{\partial x_i} \right) )} \\
  & = \displaystyle{ ( \sum_{a=1}^r (\theta \circ \pi) \left( \frac{\partial {\rm Ih}_S(g)(x_a)}{\partial x_i} \right) (\theta \circ \pi  \circ {\rm Ih}_S(g)) \left( \frac{\partial f'}{\partial x_a} \right)) }
  \\
  & = {\rm Gass}_S(g) \chi_{\Lambda_r}(g)(\nu_1(f' \; \mbox{mod}\; \frak{F}_r'')). 
  \end{array} 
$$
When ${\rm Meta}_S$ and ${\rm Gass}_S$ are restricted to ${\rm Gal}_k[1]$, by the first assertion, we have the commutative diagram of $\Lambda_r$-modules for any $g \in {\rm Gal}_k[1]$:
$$ \begin{array}{rcccl}
                      & \frak{L}_r  & \stackrel{\nu_1}{\hookrightarrow} &  \Lambda_r^{\oplus r} & \\
{\rm Meta}_S(g) & \downarrow  &                                                   & \downarrow & {\rm Gass}_S(g) \\
& \frak{L}_r & \stackrel{\nu_1}{\hookrightarrow} & \Lambda_r^{\oplus r},& 
\end{array}
$$
from which the latter assertion follows. $\;\; \Box$\\

 Next, we introduce the pro-$l$ reduced Gassner cocycle associated to the Ihara representation ${\rm Ih}_S$. For this, we follow Oda's arguments ([O2]). We first define a certain $\Lambda_r$-submodule $\frak{L}_r^{\rm prim}$ of $\frak{L}_r$,  which Oda calls the primitive part of $\frak{L}$, as follows. 
For $1\leq i \leq r$, let $\frak{N}_i$ be the closed subgroup generated normally by $x_i$ and let $\frak{F}_{r}^{(i)} := \frak{F}_r/\frak{N}_i$. 
Let $\Lambda_r^{(i)} := \mathbb{Z}_l[[u_1,\dots , \hat{u_i},\dots , u_r]] \simeq \mathbb{Z}_l[[\frak{F}_r^{(i)}/(\frak{F}_r^{(i)})']]$ ($\hat{u_i}$ means deleting $u_i$) with augmentation ideal $I_{\Lambda_r^{(i)}}$, and let  $\delta_i : \Lambda_r \rightarrow \Lambda_r^{(i)}$ be the $\mathbb{Z}_l$-algebra homomorphism defined by $\delta_i(u_j) := u_j$ if $j \neq i$ and $\delta_i(u_i) := 0$. Note that any $\Lambda_r^{(i)}$-module is regarded as a $\Lambda_r$-module via $\delta_i$. Let $\frak{L}_r^{(i)} := (\frak{F}_r^{(i)})'/(\frak{F}_r^{(i)})''$ and let $\xi_i : \frak{L}_r \rightarrow \frak{L}_r^{(i)}$ be the $\Lambda_r$-homomorphism induced by the natural homomorphism $\frak{F}_r \rightarrow \frak{F}_r^{(i)}$.  Then the {\em primitive part} $\frak{L}_r^{\rm prim}$ of $\frak{L}_r$ is defined by 
$$ \displaystyle{ \frak{L}_r^{\rm prim} := \bigcap_{i=1}^r \;{\rm Ker}(\xi_i). }  \leqno{(4.3.8)}$$
We set $w := u_1\cdots u_r$. The following theorem and  the proof are due to Oda. \\
\\
{\bf Theorem 4.3.9} ([O2]). {\em Notations being as above, the following assertions hold.}\\
(1) {\em The Blanchfield-Lyndon isomorphism $\nu_1$ in $(4.3.6)$ restricted to $\frak{L}_r^{\rm prim}$ induces the following isomorphism of $\Lambda_r$-modules}
$$ \frak{L}_r^{\rm prim} \simeq \{ \displaystyle{(\lambda_j \frac{w}{u_j} ) \in \Lambda_r^{\oplus r}} \; | \; \lambda_j \in \Lambda_r, \; \sum_{j=1}^r \lambda_j = 0\}.$$
{\em In particular, $\frak{L}_r^{\rm prim}$ is the free $\Lambda_r$-module of rank $r-1$ on the basis }
$$\displaystyle{{\boldsymbol v}_1 := \, ^t(-\frac{w}{u_1}, \frac{w}{u_2}, 0, \dots , 0),} \dots, \displaystyle{{\boldsymbol v}_{r-1} := \, ^t(0,\dots , 0, -\frac{w}{u_{r-1}}, \frac{w}{u_r})}.$$\\
(2) {\em $\frak{L}_r^{\rm prim}$ is stable under the action of ${\rm Gal}_k$ through ${\rm Meta}_S$ and defines $1$-cocycle}
$$ {\rm Gass}_S^{\rm red} : {\rm Gal}_k \longrightarrow {\rm GL}_{r-1}(\Lambda_r)$$
{\em with respect to the basis ${\boldsymbol v}_1, \dots , {\boldsymbol v}_{r-1}$ and the action $\chi_{\Lambda_r}$ in $(4.3.1)$.} We call ${\rm Gass}_S^{\rm red}$ the {\em pro-$l$ reduced Gassner cocycle} of ${\rm Gal}_k$ associated to ${\rm Ih}_S$.\\
\\
{\em Proof.} (1)  We define the $\Lambda_r$-homomorphism $\tilde{\xi_i} :  \Lambda_r^{\oplus r} \rightarrow (\Lambda_r^{(i)})^{\oplus (r-1)}$ by
$$ \tilde{\xi_i}( \, ^t(\lambda_1,\dots,\lambda_r)) := \, ^t(\delta_i(\lambda_1),\dots, \delta_i(\lambda_{i-1}),\delta_i(\lambda_{i+1}),\dots,\delta_i(\lambda_r)).$$
Then we have $\xi_i = \tilde{\xi_i}|_{\frak{L}_r}$ for $1\leq i \leq r$ and the commutative diagram of $\Lambda_r$-modules:
$$ \begin{array}{ccccccccc}
0 & \longrightarrow & \frak{L}_r & \longrightarrow & \Lambda_r^{\oplus r} & \longrightarrow & I_{\Lambda_r} & \longrightarrow & 0\\
  &    & \downarrow \xi_i  &                                         & \downarrow \tilde{\xi_i} &                                     & \downarrow \delta_i&                                & \\
0 & \longrightarrow & \frak{L}_r^{(i)}& \longrightarrow & (\Lambda_r^{(i)})^{\oplus (r-1)} & \longrightarrow & I_{\Lambda_r^{(i)}} & \longrightarrow & 0,\\
\end{array}
$$ 
where two rows are the pro-$l$ Crowell exact sequences. It is easy to see that ${\rm Ker}(\tilde{\xi_i})$ is given by
$$   {\rm Ker}(\tilde{\xi}) = \{\, ^t(\lambda_1 u_i, \dots, \lambda_{i-1}u_i, \lambda_i, \lambda_{i+1}u_i,\dots, \lambda_r u_i) \; | \; \lambda_j \in \Lambda_r \;(1\leq j \leq r)\}$$
and hence,  by (4.3.6) and (4.3.8), we have
$$ \frak{L}_r^{\rm prim} = \{ (\lambda_j) \in \Lambda_r^{\oplus r} \; | \; \sum_{j=1}^r \lambda_j u_j = 0, \;\; \lambda_j \equiv 0 \; \mbox{mod}\; u_i \; \mbox{if}\; i \neq j   \}.$$
Since $\Lambda_r$ is a regular local ring, it is factorial. Therefore we have the first assertion
$$ \frak{L}_r^{\rm prim} = \{ (\lambda_j) \in \Lambda_r^{\oplus r} \; | \; \sum_{j=1}^r \lambda_j u_j = 0, \;\; \lambda_j \equiv 0 \; \mbox{mod}\; \frac{w}{u_j}  \; (1\leq j \leq r)\}.$$
The assertion for a basis of $\frak{L}_r^{\rm prim}$ is clear. \\
(2) Since ${\rm Ih}_S(g)(x_i)$ is conjugate to $x_i^{\chi_l(g)}$ for $g \in {\rm Gal}_k$ and $1\leq i \leq r$, the definition (4.3.8) implies that $\frak{L}_r^{\rm prim}$ is ${\rm Gal}_k$-stable under the action ${\rm Meta}_S$. So we may write, for $1\leq j \leq r-1$,
$$ {\rm Ih}_S(g)({\boldsymbol v}_j) = \sum_{i=1}^{r-1} {\rm Gass}_S^{\rm red}(g)_{ij} {\boldsymbol v}_i,  \leqno{(4.3.9.1)}$$
where ${\rm Gass}_S^{\rm red}(g)_{ij} \in \Lambda_r$ is the $(i,j)$-entry of the representation matrix of ${\rm Ih}_S(g)$ with respect to ${\boldsymbol v}_1, \dots , {\boldsymbol v}_{r-1}$. Then we have, for $g, h \in {\rm Gal}_k$,
$$ \begin{array}{ll}
{\rm Ih}_S(gh)({\boldsymbol v}_j) & = {\rm Ih}_S(g)({\rm Ih}_S(h)({\boldsymbol v}_j))\\
 & = \displaystyle{ {\rm Ih}_S(g) \left( \sum_{i=1}^{r-1} {\rm Gass}_S^{\rm red}(h)_{ij} {\boldsymbol v}_i  \right) }\\
 & =   \displaystyle{ \sum_{i=1}^{r-1} \chi_{\Lambda_r}({\rm Gass}_S^{\rm red}(h)_{ij}) {\rm Ih}_S(g)({\boldsymbol v}_i) } \;\; (\mbox{by} \; (4.3.1))\\
 & =  \displaystyle{ \sum_{t=1}^{r-1} \left( \sum_{i=1}^{r-1} {\rm Gass}_S^{\rm red}(g)_{ti} \chi_{\Lambda_r}(g)({\rm Gass}_S^{\rm red}(h)_{ij}) \right) {\boldsymbol v}_t, }
\end{array}
$$
which means the cocycle relation
$$ {\rm Gass}_S^{\rm red}(gh) = {\rm Gass}_S^{\rm red}(g) \chi_{\Lambda_r}(g)({\rm Gass}_S^{\rm red}(h)).$$
Hence the assertion is proved. $\;\; \Box$
\\
\\
When we restrict ${\rm Gass}_S^{\rm red}$ to ${\rm Gal}_k[1]$, we have a representation
$$ {\rm Gass}_S^{\rm red} : {\rm Gal}_k[1] \longrightarrow {\rm GL}(r-1; \Lambda_r), $$
which we call the {\em pro-$l$ reduced Gassner representation} of ${\rm Gal}[1]$ associated to ${\rm Ih}_S$.  \\

Let $\Gamma$ be a free pro-$l$ group of rank $1$ generated by $x$ so that $\mathbb{Z}_l[[\Gamma]]$ is identified with the Iwasawa algebra $\Lambda := \mathbb{Z}_l[[u]]$ ($x \leftrightarrow 1+u$). Let 
$\frak{z} : \frak{F}_r \rightarrow \Gamma$ be the homomorphism defined by $\frak{z}(x_i) := x$ for $1\leq i \leq r$. Let $\chi_{\Lambda}$ be the action of ${\rm Gal}_k$ on $\Lambda$ defined by $ \chi_{\Lambda}(g)(u) := (1+u)^{\chi_l(g)} -1$ for $g \in {\rm Gal}_k$
Then we have the pro-$l$ Magnus cocycle associated to ${\rm Ih}_S$ and $\frak{z}$
$$ {\rm Bur}_S : {\rm Gal}_k \longrightarrow {\rm GL}(r; \Lambda),$$
which we call the {\em pro-$l$ Burau cocycle} of ${\rm Gal}_k$ associated to ${\rm Ih}_S$. It is the 1-cocycle of ${\rm Gal}_k$ with coefficients in ${\rm GL}(r;\Lambda)$ with respect to the action $\chi_{\Lambda}$.  By definition, we have
$$ {\rm Bur}_S(g) = {\rm Gass}_S(g)|_{u_1=\cdots = u_r = u}.$$
Similarly, we have the {\em pro-$l$ reduced Burau cocycle} associated to ${\rm Ih}_S$
$$ {\rm Bur}_S^{\rm red} : {\rm Gal}_k \longrightarrow {\rm GL}(r-1;\Lambda)$$
defined by
$$ {\rm Bur}_S^{\rm red}(g) := {\rm Gass}_S^{\rm red}(g)|_{u_1=\cdots = u_r = u}.$$
Since $(\frak{z}\circ {\rm Ih}_S(g))(x_i) = \frak{z}(y_i(g)x_iy_i(g)^{-1}) = \frak{z}(x_i)$ for $g \in {\rm Gal}_k[1]$, we have
$$ \frak{z} \circ {\rm Ih}_S(g) = \frak{z}    \;\; (g \in {\rm Gal}_k[1]).$$
So, when we restrict ${\rm Bur}_S$ and ${\rm Bur}_S^{\rm red}$ to ${\rm Gal}_k[1]$, we have representations
$${\rm Bur}_S : {\rm Gal}_k[1] \rightarrow {\rm GL}_r(\Lambda), \;\; {\rm Bur}_S^{\rm red} : {\rm Gal}_k[1] \rightarrow {\rm GL}_{r-1}(\Lambda),$$
which are called the {\em pro-$l$ Burau representation} and the {\em pro-$l$ reduced Burau representation} of ${\rm Gal}_k[1]$ associated to ${\rm Ih}_S$, respectively.
\\

\begin{center}
{\bf 5.  $l$-adic Alexander invariants}
\end{center}
\vspace{.2cm}
{\bf 5.1. Pro-$l$ link modules.} Let $g \in {\rm Gal}_k$. As in (2.3.1), let $\Pi_S(g)$ be the pro-$l$ link group of $g$ associated to the Ihara representation ${\rm Ih}_S$:
$$ \begin{array}{ll} \Pi_S(g) & = \langle x_1, \dots , x_r \, | \, y_1(g)x_1^{\chi_l(g)}y_1(g)^{-1}x_1^{-1} = \cdots =  y_r(g)x_r^{\chi_l(g)}y_r(g)^{-1}x_r^{-1} = 1 \rangle \\
 & = \frak{F}_r/\frak{N}_S(g), 
\end{array}$$
where $\frak{N}_S(g)$ is the closed subgroup of $\frak{F}_r$ generated normally by the pro-$l$ words $y_1(g)x_1^{\chi_l(g)}y_1(g)^{-1}x_1^{-1}, \dots ,  y_r(g)x_r^{\chi_l(g)}y_r(g)^{-1}x_r^{-1}$. Let $\psi : \frak{F}_r \rightarrow \Pi_S(g)$ be the natural homomorphism and let $\gamma_i := \psi(x_i)$ ($1\leq i \leq r$). 
Recall that $\frak{a}(g)$ denotes the ideal of $\mathbb{Z}_l$ generated by $\chi_l(g)-1$.  Then we have
$$ \Pi_S(g)/\Pi_S(g)' = \mathbb{Z}_l/\frak{a}(g) [\gamma_1] \oplus \cdots \oplus \mathbb{Z}_l/\frak{a}(g) [\gamma_r] \simeq (\mathbb{Z}_l/\frak{a}(g))^{\oplus r},$$
 where $[\gamma_i] := \gamma_i$ mod $\Pi_S(g)'$ ($1\leq i \leq r$). The correspondence $\gamma_i \mapsto u_i$ induces the $\mathbb{Z}_l$-algebra isomorphism
 $$ \theta(g) : \mathbb{Z}_l[[\Pi_S(g)/\Pi_S(g)']] \simeq \Lambda_r/((1+u_1)^{\chi_l(g)-1}-1, \dots , (1+u_r)^{\chi_l(g)-1}-1).$$
 We denote the right hand side by $\Lambda_r(g)$:
 $$ \Lambda_r(g) := \Lambda_r/((1+u_1)^{\chi_l(g)-1}-1, \dots , (1+u_r)^{\chi_l(g)-1}-1),$$
 and by $I_{\Lambda_r(g)}$ the augmentation ideal of $\Lambda_r(g)$. 
 
 We define the {\em pro-$l$ link module} $\frak{L}_S(g)$ of $g$ associated to ${\rm Ih}_S$ by
$$ \frak{L}_S(g) := \Pi_S(g)'/\Pi_S(g)'',$$
which is considered as a $\Lambda_r(g) = \mathbb{Z}_l[[\Pi_S(g)/\Pi_S(g)']]$-module.  It may be seen as an analogue of the classical link module in link theory (cf. [Hi], [Ms2, Chapter 9]).

Let $\varpi : \Pi_S(g) \rightarrow \Pi_S(g)/\Pi_S(g)'$ be the abelianization map. We define the {\em pro-$l$ Alexander module} $\frak{A}_S(g)$ of $g$ associated to ${\rm Ih}_S$ by the {\em pro-$l$ differential module} associated to $\varpi$,  namely the quotient module of the free $\Lambda_r(g)$-module on symbols $d\gamma$ for $\gamma \in \Pi_S(g)$ by the $\Lambda_r(g)$-submodule generated by 
$d(\gamma_1\gamma_2) - d\gamma_1 - \varpi(\gamma_1)d\gamma_2$ for $\gamma_1, \gamma_2 \in \Pi_S(g)$ ([Ms; 9.3]):
$$ \frak{A}_S(g) := \bigoplus_{\gamma \in \Pi_S(g)} \Lambda_r(g) d\gamma / \langle d(\gamma_1\gamma_2) - d\gamma_1 - \varpi(\gamma_1)d\gamma_2 \; (\gamma_1, \gamma_2 \in \Pi_S(g) ) \rangle_{\Lambda_r(g)}.$$
We define the {\em $l$-adic Alexander matrix $Q_S(g)$}  by the Jacobian matrix of the relators of $\Pi_S(g)$:
$$ Q_S(g) := \left( (\theta(g) \circ \varpi \circ \psi )\left(  \frac{\partial\, y_j(g)x_j^{\chi_l(g)}y_j(g)^{-1} x_j^{-1}}{\partial x_i} \right) \right). \leqno{(5.1.1)} $$
\\
{\bf Proposition 5.1.2.} {\em Notations being as above, the following assertions hold.}\\
(1) {\em The correspondence $\displaystyle{ d\gamma \mapsto ( (\theta(g) \circ \varpi \circ \psi)\left( \frac{\partial f}{\partial x_i} \right)) }$ gives the isomorphism}
$$\frak{A}_S(g) \stackrel{\sim}{\longrightarrow} {\rm Coker}( Q_S(g) : \Lambda_r(g)^{\oplus r} \rightarrow \Lambda_r(g)^{\oplus r}),$$
{\em where $f$ is any element of $\frak{F}_r$ such that $\gamma = \psi(f)$.}\\
(2) (Pro-$l$ Crowell exact sequence) {\em We have the following exact sequence of $\Lambda_r(g)$-modules}:
$$ 0 \longrightarrow \frak{L}_S(g) \stackrel{\nu_1}{\longrightarrow} \frak{A}_S(g) \stackrel{\nu_2}{\longrightarrow} I_{\Lambda_r(g)} \longrightarrow 0,$$
{\em where $\nu_1, \nu_2$ are given by}
$$\displaystyle{ \nu_1(\gamma' \; \mbox{mod}\; \Pi_S(g)'') := d\gamma \; (\gamma' \in \Pi_S(g)'); \;\; 
 \nu_2(d\gamma) := (\theta(g) \circ \varpi)(\gamma) -1 \; (\gamma \in \Pi_S(g)).} $$
\\
{\em Proof.} We refer to [Ms2; Theorems 9.3.6, 9.4.2].\\
\\
Let $\phi_g : \Lambda_r \rightarrow \Lambda_r(g)$ be the natural $\mathbb{Z}_l$-algebra homomorphism.\\
\\
{\bf Proposition 5.1.3.} {\em We have}
$$ Q_S(g) = \phi_g({\rm Gass}_S(g) - I)$$
and its $(i,j)$-entry is given by
$$ Q_S(g)_{ij} = \left\{ 
\begin{array}{ll}
\displaystyle{\phi_g\left(\sum_{n\geq 1} \sum_{{\scriptstyle1\leq i_1,\dots , i_n\leq r}\atop{\scriptstyle  i_n \neq i }} \mu(g; i_1\cdots i_n i) u_{i_1} \cdots u_{i_n}\right) } & \; (i=j),\\
\displaystyle{\phi_g \left( - u_j \left(\mu(g;ij) + \sum_{n\geq 1} \sum_{1\leq i_1,\dots , i_n\leq r} \mu(g;i_1\cdots i_n i i) u_{i_1}\cdots u_{i_n}\right) \right) } & \; (i \neq j).\\
\end{array}\right.
$$
{\em Proof.} By the definition (5.1.1), we have 
$$  Q_S(g)_{ij} :=  (\theta(g) \circ \varpi \circ \psi )\left(  \frac{\partial\, y_j(g)x_j^{\chi_l(g)}y_j(g)^{-1} x_j^{-1}}{\partial x_i} \right).$$
By the basic rules of pro-$l$ Fox free derivatives, we have
$$ \frac{\partial\, y_j(g)x_j^{\chi_l(g)}y_j(g)^{-1} x_j^{-1}}{\partial x_i} = \frac{\partial\, y_j(g)x_j^{\chi_l(g)}y_j(g)^{-1}}{\partial x_i}
 - \delta_{ij} y_j(g)x_j^{\chi_l(g)}y_j(g)^{-1}x_j^{-1}.$$
By (4.3.3) and $\theta(g) \circ \varpi \circ \psi = \phi_g \circ \theta \circ \pi$, we have
$$ (\theta(g) \circ \varpi \circ \psi) \left(\frac{\partial\, y_j(g)x_j^{\chi_l(g)}y_j(g)^{-1}}{\partial x_i}\right) = \phi_g( {\rm Gass}_S(g)_{ij}),$$
and we also have
$$ (\theta(g) \circ \varpi \circ \psi)(y_j(g)x_j^{\chi_l(g)}y_j(g)^{-1}x_j^{-1}) = \theta(g)(\gamma_j^{\chi_l(g)-1}) = (1+u_j)^{\chi_l(g)-1} = 1.$$
Therefore we have
$$ Q_S(g)_{ij} = \phi_g( {\rm Gass}_S(g)_{ij} - \delta_{ij}).$$
The second assetion follows from Proposition 4.3.4 and 
$$ \phi_g(\chi_{\Lambda_r}(g)(u_j)) = \phi_g( (1+u_j)^{\chi_l(g)} - 1) = \phi_g(u_j). \;\; \Box$$
\\
{\bf Corollary 5.1.4.} {\em For $g, h \in {\rm Gal}_k[1]$, we have the following isomorphisms of $\Lambda_r$-modules}
$$ \frak{A}_S(hgh^{-1}) \simeq \frak{A}_S(g), \;\; \frak{L}_S(hgh^{-1}) \simeq \frak{L}_S(g).$$
{\em Proof.} Since ${\rm Gass}_S : {\rm Gal}_k \rightarrow {\rm GL}(r;\Lambda_r)$ is a representation, we have 
$$ Q_S(hgh^{-1}) = \phi_g({\rm Gass}_S(hgh^{-1}) - I )= \phi_g({\rm Gass}_S(h))Q_S(g)\phi_g({\rm Gass}_S(h))^{-1}$$
by Proposition 5.1.3. Then the first assertion follows from Proposition 5.1.2 (1). The second assertion follows from Proposition 5.1.2 (2). $\;\; \Box$
\\
\\
{\bf 5.2. $l$-adic Alexander invariants.} For $n \geq 0$, we define the {\em $n$-th $l$-adic Alexander ideal} $\frak{E}_S(g)^{(n)}$ of $g \in {\rm Gal}_k$ associated to ${\rm Ih}_S$ by the $n$-th Fitting ideal of the pro-$l$ Alexander module $\frak{A}_S(g)$ over $\Lambda_r(g)$. The {\em $n$-th $l$-adic  Alexander invariant} $A_S(g)^{(n)}$  is then defined by a generator of the divisorial hull of $\frak{E}_S(g)^{(n)}$. By Proposition 5.1.2 (1), $\frak{E}_S(g)^{(n)}$ is the ideal generated by all $(r-n)$-minors of $Q_S(g)$ if $r-n > 0$ and $\frak{E}_S(g)^{(n)} := \Lambda_r(g)$ if $r-n \leq 0$,  and $A_S(g)^{(n)}$  is the greatest common divisor of all $(r-n)$-minors of $Q_S(g)$ if $r-n > 0$ and  $A_S(g)^{(n)} := 1$ if $r-n \leq 0$:
$$A_S(g)^{(n)} := \left\{ \begin{array}{ll} 
\mbox{g.c.d of all}\; (r-n)\mbox{-minors of}\; Q_S(g)  & \; (r-n >0),\\
1 & \; (r-n \leq 0).
\end{array}\right.
$$
We note that $A_S(g)^{(n)}$ is defined up to multiplication of a unit of $\Lambda_r(g)$. We write $\frak{E}_S(g)$ (resp. $A_S(g)$) for $\frak{E}_S(g)^{(0)}$ (resp. $A_S(g)^{(0)}$) and call $\frak{E}_S(g)$ (resp. $A_S(g)$) the $l$-adic Alexander ideal (resp. $l$-adic Alexander invariant) of $g$ associated to ${\rm Ih}_S$.
From Proposition 5.1.3, the following proposition is immediate.\\
\\
{\bf Proposition 5.2.1.} {\em For $g \in {\rm Gal}_k$, we have}
$$ A_S(g) = \phi_g(\det({\rm Gass}_S(g) - I)).$$
{\em When $g \in {\rm Gal}_k[1]$, $A_S(g) = 0$ if and only if ${\rm Gass}_S(g)$ has the eigenvalue $1$.}\\
\\
Moreover, since the $l$-adic Alexander matrix $Q_S(g)$ is described by $l$-adic Milnor numbers as in Proposition 5.1.3, $n$-th $l$-adic Alexander invariants are also described by  $l$-adic Milnor numbers (cf. [Ms2; Chapter 10], [Mu]).
\\

\begin{center}
{\bf 6. The Ihara power series} 
\end{center}

In this section, we suppose that  $S =\{ 0, 1, \infty \}$ and so $k = \mathbb{Q}$. 
In the following, we will omit $S$ in the notations. The Ihara representation in this case is 
$$ {\rm Ih} : {\rm Gal}_{\mathbb{Q}} \longrightarrow P(\frak{F}_2),$$
which factors through the Galois group ${\rm Gal}(\Omega_l/\mathbb{Q})$ by Theorem 1.2.6 (2), where  $\Omega_l$ denotes the maximal pro-$l$ extension of $\mathbb{Q}[1] = \mathbb{Q}(\zeta_{l^{\infty}})$ unramified outside $l$.
\vspace{.2cm}\\
{\bf 6.1. The Ihara power series.} The following lemma is a restatement of [Ih1; Theorem 2 (i)]. See also [Ih2; $\S 1$ (D) Example 2]. For the sake of readers, we give a proof using Theorem 4.3.9. \\
\\
{\bf Lemma 6.1.1.} {\em We have $\frak{L}_2 = \frak{L}_2^{\rm prim}$ with basis $^t(-u_2, u_1)$ over $\Lambda_2$, and $^t(-u_2,u_1) = \nu_1([x_1,x_2])$.}\\
\\
{\em Proof.} By Theorem 4.3.9 (1), $\frak{L}_2^{\rm prim}$ is the free $\Lambda_2$-module with basis $^t(-u_2,u_1)$. On the other hand,
we note that $\lambda_1 u_1 + \lambda_2 u_2 = 0$ implies $\lambda_1 = -a u_2, \lambda_2 = au_1$ for some $a \in \Lambda_2$, because $\Lambda_2$ is U.F.D.  Therefore $\frak{L}_2$ is also the free $\Lambda_2$-module with basis $^t(-u_2,u_1)$ by (4.3.6). Hence $\frak{L}_2 = \frak{L}_2^{\rm prim}$. The second assertion follows from 
$$   \displaystyle{(\theta \circ \pi)\left( \frac{\partial [x_1,x_2]}{\partial x_1} \right) = -u_2}, \; \displaystyle{ (\theta \circ \pi)\left( \frac{\partial [x_1,x_2]}{\partial x_2} \right) = u_1.} \;\;\Box$$
\\
Thanks to Lemma 6.1.1, Ihara introduced a power series $F_g(u_1,u_2) \in \Lambda_2$, called the {\em Ihara power series}, by the following equality in $\frak{L}_2$
$$ {\rm Ih}_S(g)([x_1,x_2]) \equiv F_g(u_1,u_2)[x_1, x_2] \;\; \mbox{mod} \; \frak{F}_2''. \leqno{(6.1.2)}$$
The following theorem gives an arithmetic topological interpretaion of  the Ihara power series $F_g(u_1,u_2)$. For a multi-index $I = (i_1\cdots i_n)$ with $i_j = 1$ or $2$, we denote by $|I|_1$ (resp. $|I|_2$) the number of $j$'s $(1\leq j \leq n)$ such that $i_j = 1$ (resp. $i_j = 2$). For integers $n_1, n_2\geq 0$ with $n_1 + n_2 \geq 1$ and $g \in {\rm Gal}_\mathbb{Q}$, we let
$$ \displaystyle{ \mu(g; n_1,n_2) := \sum_{|I|_1 = n_1-1, |I|_2 = n_2} \mu(g; I12) + \sum_{|I|_1 = n_1, |I|_2 = n_2-1} \mu(g; I21). } $$
We recall the pro-$l$ Gassner and the pro-$l$ reduced Gassner cocycles in (4.3.3) and (4.3.9.1):
$$ {\rm Gass} : {\rm Gal}_{\mathbb{Q}} \longrightarrow {\rm GL}(2; \Lambda_2); \; \; 
{\rm Gass}^{\rm red} : {\rm Gal}_{\mathbb{Q}} \longrightarrow \Lambda_2^{\times}. $$
\\
{\bf Theorem 6.1.3.} {\em Notations being as above,  we have, for $g \in {\rm Gal}_{\mathbb{Q}}$,}
$$  \begin{array}{ll} F_g(u_1,u_2)  & = {\rm Gass}^{\rm red}(g) \\
& = \displaystyle{\frac{\chi_{\Lambda_2}(g)(u_1u_2)}{u_1u_2} \left(1 + \sum_{n\geq 1}\sum_{{\scriptstyle 1\leq i_1,\dots , i_n\leq 2}\atop{\scriptstyle  i_n \neq i_{n+1} }} \mu(g; i_1\cdots i_n i_{n+1}) u_{i_1} \cdots u_{i_n}\right) }\\
 & =\displaystyle{\frac{\chi_{\Lambda_2}(g)(u_1u_2)}{u_1u_2} \left( 1 + \sum_{{\scriptstyle n_1, n_2 \geq 0}\atop{\scriptstyle n_1+n_1 \geq 1}}  \mu(g; n_1,n_2) u_1^{n_1}u_2^{n_2}\right). }
\end{array}\\
$$ 
{\em Proof.} Applying the $\Lambda_2$-homomorphism $\nu_1$ to (6.1.2), we have, for $g \in {\rm Gal}_k$,
$$ \displaystyle{\nu_1({\rm Ih}(g)([x_1,x_2])) = F_g(u_1,u_2)\nu_1([x_1,x_2]) = F_g(u_1,u_2)\left( \begin{array}{c} - u_2 \\ u_1 \end{array} \right).}$$
On the other hand, by the definition of ${\rm Gass}_S^{\rm red}(g)$ (cf. (4.3.9.1)), we have
$$ \nu_1({\rm Ih}(g)([x_1,x_2])) = {\rm Gass}^{\rm red}(g) \left( \begin{array}{c} - u_2 \\ u_1 \end{array} \right).$$
Hence we have 
$$ F_g(u_1,u_2)  = {\rm Gass}^{\rm red}(g). $$
By Proposition 4.3.7 and Lemma 6.1.1, we have
$$ \begin{array}{ll}
\nu_1({\rm Ih}(g)([x_1,x_2])) & = {\rm Gass}(g) \chi_{\Lambda_2}(g)(\nu_1([x_1,x_2]))\\
 & = {\rm Gass}(g) \left( \begin{array}{c} - \chi_{\Lambda_2}(g)(u_2) \\ \chi_{\Lambda_2}(g)(u_1) \end{array} \right).
 \end{array}$$
A straightforward calculation using Proposition 4.3.4 yields
 $$ \begin{array}{l} \displaystyle{ {\rm Gass}(g) \left( \begin{array}{c} - \chi_{\Lambda_2}(g)(u_2) \\ \chi_{\Lambda_2}(g)(u_1) \end{array} \right)} \\
 = \displaystyle{ \frac{\chi_{\Lambda_2}(g)(u_1u_2)}{u_1u_2} \left(1 + \sum_{n\geq 1}\sum_{{\scriptstyle1\leq i_1,\dots , i_n\leq 2}\atop{\scriptstyle  i_n \neq i_{n+1} }} \mu(g; i_1\cdots i_n i_{n+1}) u_{i_1} \cdots u_{i_n}\right) \left( \begin{array}{c} - u_2 \\ u_1 \end{array} \right) } \\
 = \displaystyle{ \frac{\chi_{\Lambda_2}(g)(u_1u_2)}{u_1u_2} \left( 1 + \sum_{{\scriptstyle n_1, n_2 \geq 0}\atop{\scriptstyle n_1+n_1 \geq 1}}  \mu(g; n_1,n_2) u_1^{n_1}u_2^{n_2}\right) \left( \begin{array}{c} - u_2 \\ u_1 \end{array} \right).} 
\end{array}
$$
Getting these together, we obtain the assertion. $\;\; \Box$
\\

Ihara also interpreted $\frak{L}_2$ in terms of Fermat Jacobians. For a positive integer $n$, let $C_n$ be the non-singular, projective curve over $\mathbb{Q}$ defined by
$$ X^{l^n} + Y^{l^n} = Z^{l^n}$$
and let ${\rm Jac}_n$ be the Jacobian variety of $C_n$. Let ${\rm T}({\rm Jac}_n)$ be the $l$-adic Tate module of ${\rm Jac}_n$:
$$ {\rm T}({\rm Jac}_n):= {\rm Hom}(\mathbb{Q}_l/\mathbb{Z}_l, {\rm Jac}_n(\overline{\mathbb{Q}})) \simeq  H_1^{\rm sing}(C_n(\mathbb{C}), \mathbb{Z}) \otimes \mathbb{Z}_l,$$
 and let 
$$ \mathbb{T} := \lim_{{\scriptstyle \longleftarrow}\atop{\scriptstyle  n }} {\rm T}({\rm Jac}_n),$$
where the inverse limit is taken with respect to the maps  ${\rm T}({\rm Jac}_{n+1}) \rightarrow {\rm T}({\rm Jac}_n)$ induced by the morphisms $C_{n+1} \rightarrow C_n$;  $(X,Y,Z) \mapsto (X^l, Y^l,Z^l)$. 
Let $g_{X,n}, g_{Y,n}$ be the automorphisms of $\overline{C_n} := C_n \times_{{\rm Spec}\, \mathbb{Q}} {\rm Spec} \,\overline{\mathbb{Q}}$ over $\mathbb{P}^1_{\overline{\mathbb{Q}}}$ defined by
$$ g_{X,n} : (X,Y,Z) \mapsto (\zeta_{l^n}X, Y, Z), \;\; g_{Y,n} : (X,Y,Z) \mapsto (X, \zeta_{l^n}Y, Z)$$
and set $\displaystyle{ g_X := \lim_{\longleftarrow} g_{X,n}, g_Y := \lim_{\longleftarrow} g_{Y,n}.}$ Then ${\rm Gal}(\overline{C_n}/\mathbb{P}^1_{\overline{\mathbb{Q}}}) = (\mathbb{Z}/l^n\mathbb{Z}_l)g_{X,n} \oplus (\mathbb{Z}/l^n\mathbb{Z})g_{Y,n}$ and so $\displaystyle{ \lim_{\longleftarrow} \mathbb{Z}_l[{\rm Gal}(\overline{C_n}/\mathbb{P}^1_{\overline{\mathbb{Q}}})] \simeq \Lambda_2} $ by the correspondence $ g_X \mapsto 1+u_1, g_Y \mapsto 1+ u_2$. 
Thus $\mathbb{T}$ is regarded as a $\Lambda_2$-module.  Then we have the isomorphism of $\Lambda_2$-modules
$$ \frak{L}_2 \simeq \mathbb{T}.$$
For an explicit construction of the basis of $\mathbb{T}$ corresponding to $[x_1,x_2]$, we consult [Ae; $\S 13$].\\

Now, the main results in [Ih1] are arithmetic descriptions of\\
$\bullet$ values of $F_g(u_1,u_2)$ at $l$-powerth roots of unity in terms of the Jacobi sums which arise from the Galois action on ${\rm T}({\rm Jac}_n)$, and\\
$\bullet$ coefficients of $F_g(u_1,u_2)$ in terms of $l$-adic Soul\'{e} cocycles which are defined by the Galois action on higher cyclotomic $l$-units.\\

We will describe these, using Theorem 6.1.3, from the view point of arithmetic topology.\\
\\
{\bf 6.2. Values of the Ihara power series.}  Let $p$ be a rational prime which is in $R_S$ of (1.2.5) and let $\overline{p}$ be a prime of $\overline{\mathbb{Q}}$ lying over $p$. By Theorem 1.2.6 (2), $\overline{p}$ is unramified in $\Omega_S/\mathbb{Q}$ and so we have the Frobenius automorphism $\sigma_{\overline{p}} \in {\rm Gal}(\Omega_S/\mathbb{Q})$. Let $n$ be a fixed positive integer.  Let $\frak{p}_n$ be the prime of $\mathbb{Q}(\zeta_{l^n})$ lying below $\overline{p}$ and let $\displaystyle{\left( \frac{x}{\frak{p}_n}\right)_{l^n}}$ denote the $l^n$-th power residue symbol at $\frak{p}_n$ for $x \in (\mathbb{Z}[\zeta_{l^n}]/\frak{p}_n)^{\times}$. For $a, b \in \mathbb{Z}/l^n\mathbb{Z} \setminus \{ 0 \}$ with $(a,b,l)=1$, we define the {\em Jacobi sum} by
$$ \displaystyle{ J_{l^n}(\frak{p}_n)^{(a,b)} = \sum_{{\scriptstyle x, y \in (\mathbb{Z}[\zeta_{l^n}]/\frak{p}_n)^{\times}}\atop{\scriptstyle x + y = -1}}  \left( \frac{x}{\frak{p}_n}\right)_{l^n}^a \left( \frac{y}{\frak{p}_n}\right)_{l^n}^b. }$$
For $l=2$, $J_{l^n}(\frak{p}_n)^{(a,b)}$ must be multiplied by $\displaystyle{ \left( \frac{-1}{\frak{p}_n} \right)^a}$. Let $f$ be the order of $p$ in $(\mathbb{Z}/l^n\mathbb{Z})^{\times}$. We note that $\sigma_{\overline{p}}^f \in {\rm Gal}(\Omega_S/\mathbb{Q}(\zeta_{l^n}))$. By using Weil's theorem, Ihara showed the following \\
\\
{\bf Theorem 6.2.1} ([Ih1; Theorem 7]). {\em Let $a, b\in \mathbb{Z}/l^n\mathbb{Z} \setminus \{ 0 \}$ such that $a+b \neq 0$ and $(a,b, a+b,l) = 1$. Then we have}
$$ F_{\sigma_{\overline{p}}^f}(\zeta_{l^n}^a -1, \zeta_{l^n}^b -1) = J_{l^n}(\frak{p}_n)^{(a,b)}.$$
\\
Combining Theorem 6.1.3 and Theorem 6.2.1, we obtain the following $l$-adic expansion of the Jacobi sum $J_{l^n}({\frak{p}_n})^{(a,b)}$ with coefficients $l$-adic Milnor numbers.\\
\\
{\bf Theorem 6.2.2.} {\em Notations being as above, we have}
$$ J_{l^n}(\frak{p}_n)^{(a,b)} =\displaystyle{1 + \sum_{{\scriptstyle  n_1, n_2 \geq 0}\atop{\scriptstyle  n_1 + n_2 \geq 1}} \mu(\sigma_{\overline{p}}^f; n_1, n_2) (\zeta_{l^n}^a - 1)^{n_1} (\zeta_{l^n}^b - 1)^{n_2}.}$$
{\em Proof.} Since we have $\zeta_{l^n}^{\chi_l(\sigma_{\overline{p}}^f)} = \zeta_{l^n}^{p^f} = \zeta_{l^n}$ by $p^f \equiv 1$ mod $l^n$, the formula follows from Theorem 6.1.3 and Theorem 6.2.1. $\;\; \Box$
\\
\\
{\bf 6.3. Coefficients of the Ihara power series.} We will combine Theorem 6.1.3 with the result of Ihara, Kaneko and Yukinari on the Ihara power series ([IKY]) and deduce some formulas relating  our $l$-adic Milnor numbers with the Soul\'{e} cocycles ([So]).  As in 1.2,  let $\zeta_{l^n}$ be a primitive $l^n$-th root of unity for a positive integer $n$ such that $(\zeta_{l^{n+1}})^l = \zeta_{l^n}$ for $n \geq 1$.
 For $a \in \mathbb{Z}/l^n\mathbb{Z}$, let $\langle a \rangle_{l^n}$ denote the integer such that $0\leq \langle a \rangle_{l^n} < l^n$ and $a = \langle a \rangle_{l^n}$ mod $l^n$. For a positive interger $m$, we let
$$ \varepsilon_{l^n}^{(m)} := \prod_{ a \in (\mathbb{Z}/l^n\mathbb{Z})^{\times}} (\zeta_{l^n}-1)^{\langle a^{m-1} \rangle_{l^n}},$$
which is an $l$-unit in $\mathbb{Q}(\zeta_{l^n})$, called a {\em cyclotomic $l$-unit}. Then we define the {\em $m$-th $l$-adic Soul\'{e} cocycle} $\chi^{(m)} : {\rm Gal}_{\mathbb{Q}} \rightarrow \mathbb{Z}_l$ by the Kummer cocycle attached to the system of cyclotomic $l$-units $\{\varepsilon_{l^n}^{(m)}\}_{n \geq 1}$ 
$$ \zeta_{l^n}^{\chi^{(m)}(g)} = \{(\varepsilon_{l^n}^{(m)})^{1/l^n}\}^{g-1} \;\; (n \geq 1, g \in {\rm Gal}_{\mathbb{Q}}).$$
It is easy to see the cocycle relation
$$ \chi^{(m)}(gh) = \chi^{(m)}(g) + \chi_l(g)\chi^{(m)}(h) \;\;\;\; (g, h \in {\rm Gal}_{\mathbb{Q}})$$
and hence the restriction of $\chi^{(m)}|_{{\rm Gal}_{\mathbb{Q}[1]}}$ is a character. Let $\Omega_l^{\rm ab}$ be the maximal abelian subextension of $\Omega_l/\mathbb{Q}[1]$. Since $\mathbb{Q}(\zeta_{l^n}, (\varepsilon_{l^n}^{(m)})^{1/l^n})$ is a cyclic extension of $\mathbb{Q}(\zeta_{l^n})$ unramified outside $l$,  we have  $(\varepsilon_{l^n}^{(m)})^{1/l^n} \in \Omega_l^{\rm ab}$ and so the Soul\'{e} character  $\chi^{(m)}|_{{\rm Gal}_{\mathbb{Q}}[1]}$ factors through the Galois group ${\rm Gal}(\Omega_l^{\rm ab}/\mathbb{Q}[1])$. We note by Theorem 4.2.8 (3) that the pro-$l$ reduced Gassner representation ${\rm Gass}^{\rm red}$ also factors through  ${\rm Gal}(\Omega_l^{\rm ab}/\mathbb{Q}[1])$.

We set
$$ \displaystyle{ \kappa_m(g) := \frac{\chi^{(m)}(g)}{1- l^{m-1}}} \;\;\;\; (g \in {\rm Gal}_{\mathbb{Q}}),$$
and introduce new variables $U_1, U_2$  defined by
$$ 1 + u_i  = \exp(U_i) = \sum_{n=0}^{\infty} \frac{U_i^n}{n!} \; \in \mathbb{Q}_l[[U_i]] \;\;\;\; (i=1, 2)$$
and set
$$ {\cal F}_g(U_1,U_2) := F_g(u_1,u_2)|_{u_i = \exp(U_i)-1}.$$
\\
{\bf Theorem 6.3.1.} ([IKY; Theorem ${\rm A}_2$]). {\em Notations being as above, we have, for $g \in {\rm Gal}(\Omega_l^{\rm ab}/\mathbb{Q}[1])$,}
$$ \displaystyle{ {\cal F}_g(U_1,U_2) = \exp \left\{ - \sum_{{\scriptstyle m \geq 3}\atop{\scriptstyle  {\rm odd} }}  \kappa_m(g) \left( \sum_{{\scriptstyle m_1, m_2 \geq 1}\atop{\scriptstyle m_1+m_2=m}} \frac{U_1^{m_1} U_2^{m_2}}{m_1 ! \, m_2!} \right) \right\}. } $$
\\
Combining Theorem 6.1.3 and Theorem 6.3.1, we can deduce relations between $l$-adic Milnor numbers and $l$-adic Soul\'{e} characters. For this, we prepare
the following \\
\\
{\bf Lemma 6.3.2.} {\em Let $a(n_1,n_2)$ and $c(m_1,m_2)$ be given $l$-adic numbers for integers $m_1,m_2,n_1,n_2 \geq 0$ with $m_1+m_2, n_1+n_2 \geq 1$. Let 
$$ A(u_1,u_2) := 1 + \sum_{{\scriptstyle n_1, n_2 \geq 0}\atop{\scriptstyle n_1+n_2 \geq 1}} a(n_1, n_2)u_1^{n_1} u_2^{n_2} \;\;  \in \mathbb{Q}_l[[u_1,u_2]]$$
and set
$$ \begin{array}{ll} 
B(U_1,U_2) & := A(u_1,u_2)|_{u_i = \exp(U_i)-1}\\
                 & = \displaystyle{ 1 + \sum_{{\scriptstyle N_1, N_2 \geq 0}\atop{\scriptstyle N_1+N_2 \geq 1}} b(N_1,N_2) U_1^{N_1} U_2^{N_2} } \;\; \in \mathbb{Q}_l[[U_1,U_2]].
\end{array}$$
Then we have
$$ \displaystyle{ b(N_1,N_2) = \sum_{{\scriptstyle n_1 + n_2 \geq 1}\atop{\scriptstyle 0\leq n_1 \leq N_1, 0\leq n_2 \leq N_2}}  a(n_1,n_2)a_{n_1}(N_1)a_{n_2}(N_2),}$$
where for $j = 1,2$, 
$$ a_{n_j}(N_j) := \left\{ 
\begin{array}{ll}
1 & \;\;\;\; (n_j = 0),\\
\displaystyle{ \sum_{{\scriptstyle e_1, \dots, e_{n_j} \geq 1}\atop{\scriptstyle e_1 + \cdots + e_{n_j} = N_j}}    \frac{1}{e_1 ! \cdots e_{n_j} !}   }       & \;\;\;\;  (n_j \geq 1).
\end{array}\right.
$$}
{\em Let 
$$ \displaystyle{ C(U_1,U_2) := \sum_{{\scriptstyle m_1, m_2 \geq 0}\atop{\scriptstyle m_1+m_2 \geq 1}} c(m_1,m_2) U_1^{m_1} U_2^{m_2} \;\; \in \mathbb{Q}_l[[U_1,U_2]] }$$
and set
$$ \begin{array}{ll} 
D(U_1,U_2) & := \exp(C(U_1,U_2))\\
                 & = \displaystyle{ 1 + \sum_{{\scriptstyle N_1, N_2 \geq 0}\atop{\scriptstyle N_1+N_2 \geq 1}} d(N_1,N_2) U_1^{N_1} U_2^{N_2} } \;\; \in \mathbb{Q}_l[[U_1,U_2]].
\end{array}$$
Then we have
$$ \displaystyle{  d(N_1,N_2) = \sum_{1\leq n \leq N_1+N_2} \frac{1}{n !} \sum c(m_1^{(1)},m_2^{(1)})\cdots c(m_1^{(n)},m_2^{(n)}),}$$
where the second sum ranges over integers $m_1^{(1)}, \dots, m_1^{(n)}, m_2^{(1)},\dots , m_2^{(n)} \geq 0$ satisfying $m_1^{(i)}+ m_2^{(i)} \geq 1$ $(1\leq i \leq n)$,  $m_1^{(1)}+ \cdots + m_1^{(n)} = N_1$ and $m_2^{(1)} + \cdots + m_2^{(n)} = N_2.$
}\\
\\
{\em Proof.} Both formulas for $b(N_1,N_2)$ and $d(N_1,N_2)$ follow from straightforward computations. $\;\;\Box$ \\
\\
We apply Lemma 6.3.2 to the case that $A(u_1,u_2) = {\rm Gass}^{\rm red}(g)$, where
$$ a(n_1,n_2) = \mu(g; n_1, n_2)$$
and $\displaystyle{ C(U_1,U_2) = \log ({\cal F}_g(U_1,U_2))}$, where
$$ c(m_1,m_2) = \left\{ 
\begin{array}{ll}
\displaystyle{ - \frac{\kappa_{m_1+m_2}(g)}{m_1 ! \,m_2 !}} & \;\; (m_1+m_2 \geq 3, \; {\rm odd}),\\
0 & \;\; {\rm otherwise}.
\end{array}\right.
$$
Then, by comparing coefficients of $U_1^{N_1}U_2^{N_2}$ in ${\rm Gass}^{\rm red}(g)|_{u_i = \exp(U_i)-1} = {\cal F}_g(U_1,U_2)$,  we obtain the following\\
\\
{\bf Theorem 6.3.3.} {\em Notations being as above, we have the following equality for $g \in {\rm Gal}_{\mathbb{Q}}[1]$}:
$$ \begin{array}{l}
\displaystyle{ \sum_{{\scriptstyle n_1 + n_2 \geq 1}\atop{\scriptstyle 0\leq n_1 \leq N_1, 0\leq n_2 \leq N_2}}  \mu(g; n_1,n_2)a_{n_1}(N_1)a_{n_2}(N_2) }\\
\;\;\;\; \displaystyle{ = \sum_{1\leq n \leq N_1+N_2} \frac{(-1)^n}{n !} \sum  \frac{\kappa_{m_1^{(1)}+m_2^{(1)}}(g)}{m_1^{(1)}! \, m_2^{(1)}!}    \cdots \frac{\kappa_{m_1^{(n)}+m_2^{(n)}}(g)}{m_1^{(n)}! \,m_2^{(n)}!}, }
\end{array}
$$
{\em where the last sum ranges over integers $m_1^{(1)}, \dots, m_1^{(n)}, m_2^{(1)},\dots , m_2^{(n)} \geq 0$ satisfying $m_1^{(i)}+ m_2^{(i)} \geq 3;$ odd  $\;(1\leq i \leq n)$,  $m_1^{(1)}+ \cdots + m_1^{(n)} = N_1$ and $m_2^{(1)} + \cdots + m_2^{(n)} = N_2.$}
\\
\\
For example, lower terms are given by
$$ \begin{array}{ll}
\mu(g; (12)) = \mu(g; (21)) = 0, \; \mu(g; (212)) + \mu(g; (121)) = 0,\\
\displaystyle{ \mu(g;(221)) + \mu(g; (2212)) + \mu(g; (1221)) + \mu(g; (2121)) = - \frac{\kappa_3(g)}{2} },\\
\displaystyle{ \mu(g;(112)) + \mu(g; (1121)) + \mu(g; (2112)) + \mu(g; (1212)) = - \frac{\kappa_3(g)}{2} }.
\end{array}
$$
\vspace{.4cm}\\

{\em Acknowledgement.} We would like to thank Hidekazu Furusho, Masanobu Kaneko, Nariya Kawazumi, Hiroaki Nakamura, Takayuki Oda, Kenji Sakugawa and Takao Satoh for communications. We also thank the referees for comments.\\
\begin{flushleft}
{\bf References}\\
{[Am]} F. Amano, On a certain nilpotent extension over $\mathbb{Q}$  of degree 64 and the 4-th multiple residue symbol, Tohoku Math. J. (2)  {\bf 66}  (2014),  no. 4, 501--522.\\
{[Ae]} G. Anderson, The hyperadelic gamma function, Invent. Math.  {\bf 95},  (1989),  no. 1, 63--131.\\
{[AI]} G. Anderson, Y. Ihara, Pro-$l$  branched coverings of $\mathbb{P}^1$ and higher circular $l$-units, Ann. of Math. (2)  {\bf 128}  (1988),  no. 2, 271--293.\\
{[Aa]} S. Andreadakis, On the automorphisms of free groups and free nilpotent groups, Proc. London Math. Soc. {\bf 15} (1965), 239--268.\\
{[Ar]} E. Artin, Theorie der Z\"{o}pfe, Abh. Math. Sem. Univ. Hamburg, {\bf 4}  (1925),  no. 1, 47--72. \\
{[B]} J. S. Birman, Braids, links, and mapping class groups, Annals of Mathematics Studies, {\bf 82}. Princeton Univ. Press, Princeton, N.J.; Univ. of Tokyo Press, Tokyo, 1974.\\
{ [CFL] } K. T. Chen, R. H. Fox and R. C. Lyndon, Free differential calculus. IV. The quotient groups of the lower central series, Ann. of Math. (2) {\bf 68} (1958), 81--95.\\
{[Da]} M. Day, Nilpotence invariants of automorphism groups, Lecture note. Available at {\small http://www.math.caltech.edu/~2010-11/1term/ma191a/}\\
{[DDMS]} J. D. Dixon, M. P. F. du Sautoy, A. Mann, D. Segal, Analytic pro-$p$  groups, Second edition. Cambridge Studies in Advanced Mathematics, {\bf 61}, Cambridge University Press, Cambridge, 1999. \\
{[Dw]} W. G.  Dwyer, Homology, Massey products and maps between groups, J. Pure Appl. Algebra  {\bf 6}  (1975), no. 2, 177--190.\\
{[F]} R. H. Fox, Free differential calculus. I. Derivation in the free group ring, Ann. of Math. (2)  {\bf 57}  (1953). 547--560. \\
{[G]} A. Grothendieck,  Rev\^{e}tements \'{E}tales et Groupe Fondamental, Lecture Notes in Mathematics, {\bf 224}, Springer-Verlag, 1971.\\
{[Ha]} N. Habegger, Milnor, Johnson, and Tree Level Perturbative Invariants, preprint, 2000.\\
{[HM]} N. Habegger,  G. Masbaum, The Kontsevich integral and Milnor's invariants, Topology  {\bf 39}  (2000),  no. 6, 1253--1289.\\
{[Hi]} J. Hillman, Algebraic invariants of links, Second edition. Series on Knots and Everything, {\bf 52}. World Scientific Publishing Co. Pte. Ltd., Hackensack, NJ, 2012.\\
{[Ih1]} Y. Ihara, Profinite braid groups, Galois representations and complex multiplications, Ann. of Math. (2)  {\bf 123}  (1986),  no. 1, 43--106.\\
{[Ih2]} Y. Ihara, On Galois representations arising from towers of coverings of $\mathbb{P}^1 \setminus \{0,1,\infty \}$, Invent. Math.  86  (1986),  no. 3, 427--459.\\
{[Ih3]} Y. Ihara, Arithmetic analogues of braid groups and Galois representations,  Braids (Santa Cruz, CA, 1986),  245--257, Contemp. Math., {\bf 78}, Amer. Math. Soc., Providence, RI, 1988.\\
{[Ih4]} Y. Ihara, The Galois representation arising from  $\mathbb{P}^1 \setminus \{0,1,\infty \}$ and Tate twists of even degree,  Galois groups over $\mathbb{Q}$  (Berkeley, CA, 1987),  299--313, Math. Sci. Res. Inst. Publ., 16, Springer, New York, 1989. \\
{[IKY]} Y. Ihara, M. Kaneko, A. Yukinari, On some properties of the universal power series for Jacobi sums,  in Galois representations and arithmetic algebraic geometry (Kyoto, 1985/Tokyo, 1986),  Adv. Stud. Pure Math., {\bf 12} North-Holland, Amsterdam, (1987), 65--86. \\
{[J1]} D. Johnson, An abelian quotient of the mapping class group ${\cal T}_g$, Math. Ann. {\bf 249} (1980), 225--242.\\
{[J2]} D. Johnson, A survey of Torelli group, Low-dimensional topology (San Francisco, Calif., 1981),  165--179, Contemp. Math., {\bf 20}, Amer. Math. Soc., Providence, RI, 1983. \\
{[Ka]} N. Kawazumi,  Cohomological aspects of Magnus expansions, arXiv:0505497, 2006.\\
{[Ki]} T. Kitano, Johnson's homomorphisms of subgroups of the mapping class group, the Magnus expansion and Massey higher products of mapping tori. Topology and its application, {\bf 69} 1996, 165--172.\\
{[Ko1]} H. Kodani, On Johnson homomorphisms and Milnor invariants for pure braids, preprint, 2015.\\
{[Ko2]} H. Kodani, On Gassner representation, the Johnson homomorphisms and the Milnor invariants  for pure braids, preprint, 2016.\\
{[Ko3]} H. Kodani, Arithmetic topology on braid and absolute Galois groups, Ph.D Thesis, Kyushu University, 2017.\\
{[Kr] } D. Kraines, Massey higher products, Trans. Amer. Math. Soc.  {\bf 124}  (1966),  431--449.\\
{[MKS]}  W. Magnus, A. Karrass, D. Solitar, Combinatorial group theory. Presentations of groups in terms of generators and relations, Reprint of the 1976 second edition. Dover Publications, Inc., Mineola, NY, 2004. \\
{[Ma]} J. P. May, Matric Massey products, J. Algebra {\bf 12}  (1969), 533--568.\\
{ [Mi] } J. Milnor, Isotopy of links, in Algebraic Geometry and Topology, A symposium in honor of S. Lefschetz ( edited by R.H. Fox, D.C. Spencer and A.W. Tucker), 280--306 Princeton University Press, Princeton, N.J., 1957.\\
{[Ms1]} M. Morishita, Milnor invariants and Massey products for prime numbers, Compos. Math.  {\bf 140}  (2004),  no. 1, 69--83.\\
{[Ms2]} M. Morishita, Knots and Primes -- An Introduction to Arithmetic Topology, Universitext, Springer, 2011.\\
{[MT]} M. Morishita, Y. Terashima, $p$-Johnson homomorphisms and pro-$p$ groups, J. of Algebra. {\bf 479} (2017), 102--136.\\
{[Mt1]} S. Morita, Abelian quotients of subgroups of the mapping class group of surfaces, Duke Math. J. {\bf 70} (1993), 699--726.\\
{[Mt2]} S. Morita, The extension of Johnson's homomorphism from the Torelli group to the mapping class group, Invent. Math.  {\bf 111}  (1993),  no. 1, 197--224.\\
{[Mu]} K. Murasugi, On Milnor's invariant for links, Trans. Amer. Math. Soc.  {\bf 124}  (1966), 94--110. \\
{[MK]}  K. Murasugi, B.I. Kurpita, A study of braids, Mathematics and its Applications, {\bf 484}, Kluwer Academic Publishers, Dordrecht, 1999. \\
{[N]} H. Nakamura, Tangential base points and Eisenstein power series,  Aspects of Galois theory (Gainesville, FL, 1996),  202--217, 
London Math. Soc. Lecture Note Ser., {\bf 256}, Cambridge Univ. Press, Cambridge, 1999. \\
{[NW]} H. Nakamura, Z. Wojtkowiak,  On explicit formulae for $l$-adic polylogarithms,  In: Arithmetic fundamental groups and noncommutative algebra (Berkeley, CA, 1999),  285--294, Proc. Sympos. Pure Math., 70, Amer. Math. Soc., Providence, RI, 2002.\\
{[O1]} T. Oda, Two propositions on pro-$l$ braid groups, unpublished note, 1985.\\
{[O2]} T. Oda, Note on meta-abelian quotients of pro-$l$ free groups, unpublished note, 1985.\\
{ [R] } L. R\'{e}dei, Ein neues zahlentheoretisches Symbol mit Anwendungen auf die Theorie der quadratischen Zahlk\"{o}rper I, J. Reine Angew. Math., {\bf 180} (1939), 1-43.\\
{[Sa]} T. Satoh, A survey of the Johnson homomorphisms of the automorphism groups of free groups and related topics, Handbook of Teichm\"{u}ller theory, volume V. (editor: A. Papadopoulos), 2016, 167-209.\\
{[Se]} J.-P. Serre, Lie algebras and Lie groups, Lecture Notes in Mathematics, {\bf 1500}, Springer-Verlag, Berlin, 2006.\\
{[So]} C. Soul\'{e}, On higher $p$-adic regulators, Lecture Notes in Mathematics, {\bf 854}, Springer-Verlag, 1981, 372--401.\\
{[St] } D. Stein, Massey products in the cohomology of groups with applications to link theory, Trans.
Amer. Math. Soc. {\bf 318}  (1990), 301--325.\\
{[T] } V. Turaev, The Milnor invariants and Massey products, (Russian) Studies in topology, II. Zap. Nau\v{c}n. Sem. Leningrad. Otdel. Mat. Inst. Steklov. (LOMI)  {\bf 66},  (1976), 189--203, 209--210.\\ 
{[W1]} Z. Wojtkowiak, On $l$-adic iterated integrals. I. Analog of Zagier conjecture, Nagoya Math. J.  {\bf 176}  (2004), 113--158.\\
{[W2]}  Z. Wojtkowiak, On $l$-adic iterated integrals. II. Functional equations and $l$-adic polylogarithms, Nagoya Math. J.  {\bf 177}  (2005), 117--153.\\
{[W3]}  Z. Wojtkowiak, On $l$-adic iterated integrals. III. Galois actions on fundamental groups, Nagoya Math. J.  178  (2005), 1--36.  \\
{[W4]} Z. Wojtkowiak, A remark on nilpotent polylogarithmic extensions of the field of rational functions of one variable over $\mathbb{C}$, Tokyo J. Math.  {\bf 30}  (2007),  no. 2, 373--382.\\
\end{flushleft}
\vspace{.2cm}
{\small
H. Kodani:\\
Graduate School of Mathematics, Kyushu University, 744, Motooka, Nishi-ku, Fukuoka, 819-0395, Japan.\\
e-mail: h-kodani@math.kyushu-u.ac.jp\\
\\
M. Morishita:\\
Graduate School of Mathematics, Kyushu University, 744, Motooka, Nishi-ku, Fukuoka, 819-0395, Japan.\\
e-mail: morisita@math.kyushu-u.ac.jp \\
\\
Y. Terashima:\\
Department of Mathematical and Computing Sciences, Tokyo Institute of Technology, 2-12-1 Oh-okayama, Meguro-ku, Tokyo 152-8551, Japan.\\
e-mail: tera@is.titech.ac.jp
}

\end{document}